\documentclass[10pt]{elsarticle}
\usepackage{layout}
\usepackage{amsmath,amsfonts,amssymb}
\usepackage{xcolor}
\usepackage{hyperref}
\usepackage{epsfig}
\usepackage{todonotes}
\usepackage{algorithm}
\usepackage{algpseudocode}
\usepackage{lscape}
\usepackage{subcaption}
\usepackage{multirow}

\usepackage[font=sl,labelfont=normal,size=small]{caption}
\usepackage{rotating}
\definecolor{lightgray}{gray}{0.95}
\usepackage[a4paper,text={6.0in,9.0in},centering,includefoot,foot=0.6in]{geometry} 

\newcommand{\dix}{{\color{black}\delta_{i+1}^x}\color{black}}
\newcommand{\dig}{{\color{black}\delta_i^g}\color{black}}
\newcommand{\dis}{{\color{black}\delta_i^s}\color{black}}
\newcommand{\dil}{{\color{black}\delta_i^{\ell}}\color{black}}
\newcommand{\diz}{{\color{black}\delta_i^z}\color{black}}
\newcommand{\diq}{{\color{black}\delta_i^q}\color{black}}
\newcommand{\diu}{{\color{black}\delta_i^u}\color{black}}
\newcommand{\diy}{{\color{black}\delta_i^y}\color{black}}
\newcommand{\dir}{{\color{black}\delta_{i+1}^r}\color{black}}
\newcommand{\dik}{{\color{black}\delta_{i+1}^k}\color{black}}
\newcommand{\diw}{{\color{black}\delta_{i+1}^w}\color{black}}

\newcommand{\djx}{{\color{black}\delta_{j+1}^x}\color{black}}

\newcommand{\djs}{{\color{black}\delta_j^s}\color{black}}
\newcommand{\djl}{{\color{black}\delta_j^{\ell}}\color{black}}
\newcommand{\djz}{{\color{black}\delta_j^z}\color{black}}
\newcommand{\djq}{{\color{black}\delta_j^q}\color{black}}
\newcommand{\dju}{{\color{black}\delta_j^u}\color{black}}
\newcommand{\djy}{{\color{black}\delta_j^y}\color{black}}
\newcommand{\djr}{{\color{black}\delta_{j+1}^r}\color{black}}
\newcommand{\djk}{{\color{black}\delta_{j+1}^k}\color{black}}
\newcommand{\djw}{{\color{black}\delta_{j+1}^w}\color{black}}

\newcommand{\lef}{\left\|}
\newcommand{\rig}{\right\|}

\newcommand{\msn}{\mu\sqrt{N}}
\newcommand{\tmsn}{\tilde{\mu}\sqrt{N}}

\newcommand{\Dr}{\Delta^r}
\newcommand{\Ds}{\Delta^s}
\newcommand{\Dw}{\Delta^w}
\newcommand{\Dz}{\Delta^z}
\newcommand{\Dk}{\Delta^k}
\newcommand{\Dl}{\Delta^{\ell}}

\newcommand{\glred}{\textsc{glred}}
\newcommand{\axpy}{\textsc{axpy}}
\newcommand{\dotpr}{\textsc{dotpr}}
\newcommand{\spmv}{\textsc{spmv}}
\newcommand{\preco}{\textsc{prec}}

\date{\today}

\title{\textbf{Analyzing and improving maximal attainable accuracy in the communication hiding pipelined BiCGStab method}}

\author[ua]{\underline{S. Cools}} \ead{siegfried.cools@uantwerpen.be}
\address[ua]{Applied Mathematics Research Group, Department of Mathematics and Computer Science, University of Antwerp, Middelheimlaan 1, B-2020 Antwerp, Belgium. Contact: \underline{siegfried.cools@uantwerp.be}} 

\begin{document}

\begin{abstract} 
  Pipelined Krylov subspace methods avoid communication latency by reducing the number of global synchronization bottlenecks and by hiding global communication behind useful computational work. In exact arithmetic pipelined Krylov subspace algorithms are equivalent to classic Krylov subspace methods and generate identical series of iterates. However, as a consequence of the reformulation of the algorithm to improve parallelism, pipelined methods may suffer from severely reduced attainable accuracy in a practical finite precision setting. This work presents a numerical stability analysis that describes and quantifies the impact of local rounding error propagation on the maximal attainable accuracy of the multi-term recurrences in the preconditioned pipelined BiCGStab method. Theoretical expressions for the gaps between the true and computed residual as well as other auxiliary variables used in the algorithm are derived, and the elementary dependencies between the gaps on the various recursively computed vector variables are analyzed. The norms of the corresponding propagation matrices and vectors provide insights in the possible amplification of local rounding errors throughout the algorithm. Stability of the pipelined BiCGStab method is compared numerically to that of pipelined CG on a symmetric benchmark problem. 
	Furthermore, numerical evidence supporting the effectiveness of employing a residual replacement type strategy to improve the maximal attainable accuracy for the pipelined BiCGStab method is provided.
\end{abstract}


\maketitle

\section{Introduction}

The class of iterative methods known as Krylov subspace methods \cite{greenbaum1997iterative,liesen2012krylov,meurant1999computer,saad2003iterative,van2003iterative} are often considered as the methods of choice for solving large scale linear systems due to their computational efficiency. Given an initial guess $x_0$ for the solution $x$ and an initial residual $r_0 = b-Ax_0$, Krylov subspace methods construct a series of approximate solutions $x_i$ to the algebraic linear system of the form $Ax=b$ that lie in the $i$-th Krylov subspace
\[
	x_i \in x_0 + \mathcal{K}_i(A,r_0) = x_0 + \text{span}\{r_0, Ar_0, \ldots , A^{i-1} r_0\}, \qquad i = 1,2,\ldots,
\]
where $A$ is a real or complex non-singular square $N \times N$ matrix which is assumed to be sparse, and the corresponding right-hand side vector $b$ has length $N$. The choice of the orthogonality constraint differentiates the various Krylov subspace methods, see e.g.~\cite{saad2003iterative}. For problems with symmetric (or Hermitian) positive definite (SPD) matrices $A$ one of the most widely used Krylov subspace methods is the (preconditioned) Conjugate Gradient (CG) method, that dates back to the 1952 paper \cite{hestenes1952methods}. For solving linear systems with unsymmetric and/or indefinite linear systems -- which are the main focus of this work -- several Krylov subspace methods are eligible. Employing the popular Generalized Minimal Residual method (GMRES) is often not advisable in practical applications due to memory constraints. Methods based on Lanczos bidiagonalization \cite{golub1965calculating,saad1984practical,bjorck1988bidiagonalization,golub2012matrix} like Conjugate Gradients Squared (CGS) \cite{sonneveld1989cgs}, Bi-Conjugate Gradients (BiCG) \cite{fletcher1976conjugate}, or its stabilized variant Bi-Conjugate Gradients Stabilized (BiCGStab) \cite{van1992bi,sleijpen1994bicgstab} are usually preferred. This work focuses on the latter method, which is considered efficient from a numerical, computational and memory perspective.

From an algorithmic point of view Krylov subspace methods consist of the following three basic building blocks: dot-products between two vectors (\dotpr) that are used to guarantee basis vector orthogonality, sparse matrix-vector products (\spmv) to construct the next Krylov basis vector, and local vector operations such as scalar multiplications and additions $y \leftarrow \alpha x + y$ (\axpy). On large, parallel, distributed memory hardware the \axpy~operations are computed fully locally by each processor and require no communication, whereas an \spmv~typically requires a limited amount of communication between neighboring processors to transfer e.g.~boundary information. The main bottleneck for efficient parallel execution is however the dot-products, which require a global synchronization phase and thus induce significant communication latency, leading to worker idling. Indeed, the increasing gap between computational performance and memory/interconnect latency implies that on HPC hardware data movement is much more expensive than floating point operations (flops), both with respect to execution time as well as energy consumption \cite{dongarra2013toward,dongarra2015hpcg,dongarra2011international,fuller2011computing}.

Over the last decades various approaches on reducing or eliminating the synchronization bottleneck in Krylov subspace methods have been proposed \cite{strakovs1987effectivity,chronopoulos1989s,chronopoulos1991s,d1992reducing,demmel1993parallel,erhel1995parallel,de1995reducing,chronopoulos1996parallel}. Recent methods that aim to eliminate global synchronization points include improved Krylov subspace methods \cite{yang2002improved,yang2002improved2,yang2003improved}, hierarchical Krylov subspace methods \cite{mcinnes2014hierarchical}, enlarged Krylov subspace methods \cite{grigori2016enlarged}, an iteration fusing Conjugate Gradient method \cite{zhuang2017iteration}, $s$-step Krylov subspace methods \cite{chronopoulos2010block,hoemmen2010communication,carson2013avoiding,carson2014residual,carson2015communication,imberti2017varying}, and pipelined Krylov subspace methods \cite{ghysels2013hiding,ghysels2014hiding,sanan2016pipelined,eller2016scalable,yamazaki2017improving}. Pipelined Krylov subspace methods aim to avoid communication latency by reducing the number of global synchronization bottlenecks and by hiding global communication behind useful computational work. The possible performance advantages of using pipelined (and other) communication reducing Krylov subspace methods were illustrated in original works \cite{ghysels2013hiding,ghysels2014hiding,cools2017communication,cornelis2017communication}. As a subclass of these methods pipelined bidiagonalization methods such as pipelined BiCGStab (denoted as `p-BiCGStab' for short) have been shown to outperform their classic counterparts when solving large unsymmetric and/or indefinite linear systems on large-scale multi-node parallel hardware, see \cite{cools2017communication}.  

Caution is however advised when using pipelined Krylov subspace methods in practice, since reorganizing a Krylov subspace algorithm into a communication reducing variant typically affects the numerical stability of the algorithm. Although pipelined Krylov subspace methods produce a series of iterates identical to classic Krylov subspace methods in exact arithmetic, their behavior in finite precision arithmetic can differ significantly. This deviation is due to contamination of the finite precision recurrence and orthogonality relations by local rounding errors, which may lead to decreased attainable accuracy and a delay of convergence respectively. The impact of finite precision round-off errors on numerical stability of classic Krylov subspace methods has been extensively studied in the literature \cite{paige1976error,paige1980accuracy,greenbaum1989behavior,greenbaum1992predicting,notay1993convergence,sleijpen1995maintaining,greenbaum1997estimating,gutknecht2000accuracy,tong2000analysis,sleijpen2001differences,strakovs2002error,strakovs2005error,meurant2006lanczos,gergelits2014composite}. Recent publications stress the importance of a similar analysis for the emerging classes of communication reducing Krylov subspace methods \cite{carson2014residual,carson2016numerical}, including the class of pipelined Krylov subspace methods \cite{cools2018analyzing,cools2018numerical}. In the original work on pipelined BiCGStab \cite{cools2017communication} it was remarked that the numerical accuracy attainable by the p-BiCGStab method is often significantly worse than the precision obtainable by the classic BiCGStab algorithm and that unstable convergence behavior can be observed for p-BiCGStab. The current work aims to establish a theoretical framework for the analysis of local rounding error behavior in pipelined BiCGStab similar to the work in \cite{cools2018analyzing}. It is shown that the numerical analysis allows to explain (and correct for) the numerical instabilities caused by the propagation of local rounding errors that stem from the finite precision recurrence relations in p-BiCGStab.

The paper is structured as follows. Section \ref{sec:analysis_classic} gives an overview of the existing numerical analysis of rounding errors in the finite precision BiCGStab algorithm while simultaneously establishing the general framework and notations that will be used throughout the manuscript. In Section \ref{sec:analysis_pipe} the pipelined BiCGStab algorithm is analyzed by deriving expressions for the gap on the residual and related auxiliary variables. The corresponding error propagation matrices characterize the possible amplification of local rounding errors throughout the p-BiCGStab algorithm. A comparison to the numerical analysis of the pipelined CG method from \cite{cools2018numerical}, to which the analysis in this work is a generalization, is also considered. Numerical results that validate the analysis are provided in Section \ref{sec:numerical}. Accuracy and performance experiments on parallel hardware provide the reader with additional insights on the analysis. In addition the use of residual replacement type strategies \cite{sleijpen1996reliable,greenbaum1997estimating,van2000residual,sleijpen2001differences,carson2014residual,cools2018analyzing} is proposed to improve attainable accuracy. The effectiveness of this approach is supported by the numerical analysis. The paper in concluded by a short discussion in Section \ref{sec:conclusions}.

\section{Numerical stability analysis of the classic BiCGStab algorithm} \label{sec:analysis_classic}

The analysis in this section relies on the rounding error analysis that was originally established in the context of the Lanczos algorithm by Paige \cite{paige1971computation,paige1972computational,paige1976error,paige1980accuracy} and that was elaborated in a number of works among which \cite{sleijpen1995maintaining,greenbaum1997estimating,demmel1997applied,greenbaum1989behavior,meurant2006lanczos,gutknecht2000accuracy,strakovs2002error,strakovs2005error}. The analysis is framed in the following finite precision framework for floating point arithmetic with machine precision given by $\epsilon$, where the round-off errors on scalar multiplication, vector summation, \spmv~application and \dotpr~computation on an $N$-by-$N$ matrix $A$, length $N$ vectors $v$ and $w$ and scalar number $\alpha$ are respectively bounded by
\begin{align*}
	\| \alpha v - \text{fl}(\alpha v) \| &\leq \| \alpha v \| \, \epsilon =  |\alpha| \, \|v\| \, \epsilon, &\qquad
	\| v + w - \text{fl}(v + w) \| &\leq (\|v\| + \|w\|) \, \epsilon,\\
	\| Av - \text{fl}(Av) \| &\leq \mu\sqrt{N} \, \|A\| \, \|v\| \, \epsilon, &\qquad
	| \left( v,w \right) - \text{fl}(\,\left(v,w \right)\,) | &\leq N \, \|v\| \, \|w\| \, \epsilon,
\end{align*}
where $\text{fl}(\cdot)$ indicates the finite precision floating point representation of the variable, $\mu$ is the maximum number of nonzeros in any row of $A$, and the norm $\|\cdot\|$ is the Euclidean 2-norm in this manuscript. To avoid confusion the maximum norm will be denoted by $\|\cdot\|_{\max}$ where applicable.

\subsection{Stability analysis for the classic BiCGStab method}

The classic BiCGStab method for the solution of the preconditioned linear system $M^{-1}A x = M^{-1}b$ is given in Algorithm \ref{algo::bicgstab}. The algorithm uses the current guess for the solution $x_i$, the residual $r_i = b-Ax_i$ and the pair of 
search directions $p_i$ and $q_i$ to construct an approximate solution $x_{i+1}$ in iteration $i$. Furthermore, the following auxiliary variables are defined:
\begin{align}
	g_i &:= M^{-1} p_i, & s_i &:= A g_i, & u_i &:= M^{-1} q_i, & y_i &:= A u_i. \label{eq:auxvars0}
\end{align}
These auxiliary vectors are computed \emph{explicitly} in Algorithm \ref{algo::bicgstab} using the relations \eqref{eq:auxvars0}. This implies that in practice, i.e.~in a finite precision framework, one computes
\begin{align}
	\bar{g}_i &:= \text{fl}(M^{-1} \bar{p}_i), & \bar{s}_i &:= \text{fl}(A \bar{g}_i), & \bar{u}_i &:= \text{fl}(M^{-1} \bar{q}_i), & \bar{y}_i &:= \text{fl}(A \bar{u}_i). \label{eq:auxvars01}
\end{align}
The finite precision equivalents, or \emph{actually computed} variants, of the true vector variables occurring throughout the algorithm are denoted by a bar symbol in this work. The solution $x_i$, the residual $r_i$ and the search directions $p_i$ and $q_i$ are not computed explicitly; instead, they are computed using recurrence relations in Algorithm \ref{algo::bicgstab}. Mathematically, i.e.~in exact arithmetic, the vector quantity resulting from the recurrence relations (e.g.\,for the residual $r_{i+1} = q_i - \omega_i y_i$) is identical to the corresponding true variable (e.g.~the true residual $b-Ax_{i+1}$). However, in finite precision arithmetic the recursively computed variables may deviate from their true counterparts.

\begin{algorithm}[t]
  \caption{\textbf{: Classic preconditioned BiCGStab}\\
	\small Standard version of the BiCGStab algorithm consisting of 4 \dotpr s, 6 \axpy s, 3 \glred~phases, 2 \spmv s and 2 \preco~applications per iteration. Each \glred~phase induces a global synchronization bottleneck leading to processor idle time. Operations are executed sequentially; each line has to be completed before the next line can be executed, and there is no overlap of the \glred~phases by other arithmetic computations.}
  \label{algo::bicgstab}
  \begin{algorithmic}[1]
  	\Function{bicgstab}{$A$, $M^{-1}$, $b$, $x_0$}
    \State $r_0 := b - Ax_0$; $p_0 := r_0$  
    \For{$i = 0,1,2, \dots$}\
		\State $g_i := M^{-1} p_i$
    \State $s_i := Ag_{i}$ 
		\State \textbf{begin global reduction} $\left( r_0, s_i \right)$ \textbf{end reduction}
    \State $\alpha_i := \left( r_0, r_i \right) / \left( r_0, s_i \right)$ 
		\State $q_i := r_i - \alpha_i s_i$
		\State $u_i := M^{-1} q_i$
		\State $y_i := A u_i$
		\State \textbf{begin global reduction} $\left( q_i, y_i \right)$; $\left( y_i, y_i \right)$ \textbf{end reduction}
		\State $\omega_i := \left( q_i, y_i \right)/\left( y_i, y_i \right)$
    \State $x_{i+1} := x_i + \alpha_{i} g_i + \omega_i u_i$
    \State $r_{i+1} := q_i - \omega_{i} y_i$
		\State \textbf{begin global reduction} $\left( r_0, r_{i+1} \right)$ \textbf{end reduction}
    \State $\beta_i := \left( \alpha_i / \omega_i \right)  \left( r_0, r_{i+1} \right) / \left( r_0, r_i \right)$
    \State $p_{i+1} := r_{i+1} + \beta_i \left( p_i - \omega_i s_i \right)$ 
    \EndFor
    \EndFunction
  \end{algorithmic}
\end{algorithm}

The following perturbed recurrence relations hold in finite precision arithmetic:
\begin{equation} 
	\bar{x}_{i+1} := \left(\bar{x}_i + \bar{\alpha}_i \bar{g}_i\right) + \bar{\omega}_i \bar{u}_i + \dix, \qquad
	\bar{r}_{i+1} := \bar{q}_i - \bar{\omega}_i \bar{y}_i + \dir, \qquad
	\bar{q}_i 		:= \bar{r}_i - \bar{\alpha}_i \bar{s}_i + \diq, \label{eq:auxrec0}
\end{equation}
where $\dix$, $\dir$ and $\diq$ represent local rounding errors that may be bounded in terms of the computed quantities $\bar{x}_i$, $\bar{g}_i$, $\bar{u}_i$, $\bar{q}_i$, $\bar{y}_i$, $\bar{r}_i$ and $\bar{s}_i$ as follows:
\begin{align}
	\|\dix\| 	& \leq \left( 2 \lef \bar{x}_i\rig + 3\lef\bar{\alpha}_i \bar{g}_i \rig + \tilde{\mu}\sqrt{N} \lef M^{-1}\rig \lef \bar{\alpha}_i \bar{p}_i \rig + 2 \lef\bar{\omega}_i\bar{u}_i\rig  + \tilde{\mu}\sqrt{N} \lef M^{-1} \rig \lef \bar{\omega}_i\bar{q}_i \rig \right) \epsilon, \notag \\
	\|\dir\| 	& \leq \left( \lef \bar{q}_i \rig + 2 \lef \bar{\omega}_i \bar{y}_i \rig + \mu \sqrt{N} \lef A \rig \lef \bar{\omega}_i \bar{u}_i \rig + \tilde{\mu} \sqrt{N} \lef A \rig \lef M^{-1} \rig \lef \bar{\omega}_i \bar{q}_i \rig \right) \epsilon, \notag \\
	\|\diq\| 	& \leq \left( \lef \bar{r}_i \rig + 2 \lef \bar{\alpha}_i \bar{s}_i \rig + \mu \sqrt{N} \lef A \rig \lef \bar{\alpha}_i \bar{g}_i \rig + \tilde{\mu} \sqrt{N} \lef A \rig \lef M^{-1} \rig \lef \bar{\alpha}_i \bar{p}_i \rig \right) \epsilon. \notag 
\end{align}
It is assumed that the default order of operations in finite precision arithmetic is from left to right unless otherwise indicated by the brackets in the expressions \eqref{eq:auxrec0} (while evidently always abiding by the default order of operations for multiplication and addition).
In these expressions $\mu$ is the maximum number of nonzeros (nnz) over all rows of $A$ and $\tilde{\mu}$ denotes the maximum nnz over all rows of $M^{-1}$. 
Note that the error bounds established in this work neglect terms of $\mathcal{O}(\epsilon^2)$, cf.~\cite{gutknecht2000accuracy}.

To study the impact of the behavior of local rounding errors on the BiCGStab iteration the difference between the true residual $b-A\bar{x}_i$ and the recursively computed residual $\bar{r}_i$ is considered. This quantity is commonly called the \emph{residual gap}, and will be denoted throughout this work as
\begin{equation}
	\Dr_{i+1} := \left(b-A\bar{x}_{i+1}\right) - \bar{r}_{i+1}. \label{eq:res_gap} 
\end{equation}
The norm of the initial residual gap $\Dr_{0}$ is the round-off error from computing $\bar{r}_0$ from $b$ and $\bar{x}_0$, i.e., $\Dr_{0} := (b-A\bar{x}_0) - \text{fl}(b-A\bar{x}_0)$, and is bounded by $\lef\Dr_{0}\rig \leq ( (\msn+1) \lef A\rig \lef \bar{x}_0\rig + \lef b \rig ) \, \epsilon$. 
By substituting the recurrence relations \eqref{eq:auxrec0} into definition \eqref{eq:res_gap} for the residual gap, an expression for the residual gap in the $i$-th iteration is derived:
\begin{align} 
	\Dr_{i+1}	&= (b-A\bar{x}_{i+1})-\bar{r}_{i+1} \notag \\
						&= b-A\bar{x}_i - \bar{\alpha}_i A \bar{g}_i - \bar{\omega}_i A \bar{u}_i - \bar{q}_i + \bar{\omega}_i \bar{y}_i - A \dix - \dir \notag \\
						&= b-A\bar{x}_i - \bar{\alpha}_i A \bar{g}_i - \bar{r}_i  + \bar{\alpha}_i \bar{s}_i - A \dix - \dir -\diq \notag \\
						&= \Dr_i - A \dix - \dir  - \diq.
\end{align}
where we use that $\bar{y}_i := A\bar{u}_i$ and $\bar{s}_i := A \bar{g}_i$ are computed explicitly in Algorithm \ref{algo::bicgstab}, see \eqref{eq:auxvars01}. In classic BiCGStab the residual gap is thus simply a superposition of the local rounding errors that are compounded by the recurrence relations in each iteration:
\begin{equation} \label{eq:classic_nomatrix}
	\Dr_{i}  =  \Dr_0 - \sum_{j=0}^{i-1} (A \djx + \djr +\djq).
\end{equation}
This expression implies that no propagation of local rounding errors takes place in classic BiCGStab.

\subsection{Expressing the residual gap for classic BiCGStab in matrix notation}

Rewriting expression \eqref{eq:classic_nomatrix} in matrix notation, which describes the residual gaps up to iteration $i$, with $B = [b,b,\ldots,b]$, $\bar{X}_{i+1} = [\bar{x}_0,\bar{x}_1,\ldots, \bar{x}_i]$ and $\bar{R}_{i+1} = [\bar{r}_0,\bar{r}_1,\ldots, \bar{r}_i]$, we obtain
\begin{align}
  \mathcal{R}_{i+1} = [\Delta^r_0,\Delta^r_1,\ldots, \Delta^r_i] 
	&= (B-A\bar{X}_{i+1})-\bar{R}_{i+1} \notag \\
	&= (A\Theta_{i+1}^{x}+\Theta_{i+1}^{r}+\Theta_{i+1}^{q}) \, \mathcal{U}_{i+1}, 
\end{align}
where the local rounding error matrices $\Theta_{i+1}^{x}$, $\Theta_{i+1}^{r}$ and $\Theta_{i+1}^{q}$ are defined as
\begin{equation}
\Theta_{i+1}^{x} = [0,-\delta_1^{x},\ldots, -\delta_i^{x}], \qquad 
\Theta_{i+1}^{r} = [\Delta_0^r,-\delta_1^{r},\ldots, -\delta_i^{r}], \qquad 
\Theta_{i+1}^{q} = [0,-\delta_0^{q},\ldots, -\delta_{i-1}^{q}],
\end{equation}
and $\mathcal{U}_{i+1}$ is an $i \times i$ upper triangular matrix with all one entries.

Since (the modulus of) all entries in $\mathcal{U}_{i+1}$ is one, no amplification of local rounding errors occurs in the classic BiCGStab algorithm, and rounding errors are merely accumulated, as stated above. The algorithm will converge until the true residual norm $\|b-A\bar{x}_i\|$ stagnates at the level of the norm of the corresponding superposed local rounding errors $\|\sum_{j=0}^{i-1} (A \djx + \djr +\djq)\|$. Note that although the true residual norm does not decrease beyond this point, the recursively computed residual norm may continue to decrease, as illustrated by the numerical experiments in Section \ref{sec:numerical}.

\section{Stability analysis for the pipelined BiCGStab method} \label{sec:analysis_pipe}

The preconditioned pipelined BiCGStab method (p-BiCGStab), Algorithm \ref{algo::pipebicgstab}, uses the current guess for the solution $x_i$, the residual $r_i = b-Ax_i$ and the search direction pair $p_i$ and $q_i$ to construct an approximate solution $x_{i+1}$ in iteration $i$. A number of auxiliary vector quantities are used in the pipelined algorithm to allow for the overlap of global communication with computations. We refer to our work in \cite{cools2017communication} for details on the derivation of the pipelined BiCGStab algorithm and the auxiliary variables. The following list of (true) vector quantities is used in the pipelined BiCGStab method:
\begin{align}
	g_i &:= M^{-1} p_i, & s_i &:= A g_i, & \ell_i &:= M^{-1} s_i, & z_i &:= A \ell_i, \notag \\
	k_i &:= M^{-1} r_i, & w_i &:= A k_i, & u_i &:= M^{-1} q_i, & y_i &:= A u_i, \notag \\
	n_i &:= M^{-1} z_i, & v_i &:= A n_i, & m_i &:= M^{-1} w_i, & t_i &:= A m_i. \label{eq:auxvars}
\end{align}
To reduce the number of computationally expensive \spmv~and \preco~operations, all but the last four (i.e.~$n_i$, $v_i$, $m_i$ and $t_i$) of these variables are computed recursively as shown in Algorithm \ref{algo::pipebicgstab}. Compared to classic BiCGStab, Algorithm \ref{algo::bicgstab}, which features 4 recurrence relations (for $q_i$, $x_{i+1}$, $r_{i+1}$ and $p_{i+1}$), corresponding to a total of 6 \axpy~operations (vector operations of the form $y \leftarrow \alpha x+y$, where $x,y$ are vectors and $\alpha$ is a scalar), the pipelined BiCGStab method, Algorithm \ref{algo::pipebicgstab}, uses 11 recursively defined variables, resulting in a total of 18 \axpy~operations.

In exact arithmetic, the auxiliary recurrence relations each produce a vector that is mathematically identical to the result that is obtained by explicitly applying the matrix (or preconditioner respectively) as defined by \eqref{eq:auxvars}. However, in finite precision, numerical round-off errors in the \axpy~operations cause the recursively computed quantities to deviate from their true values. We aim to characterize the gaps between the recursive and true quantities in pipelined BiCGStab.

\begin{algorithm}[t]
  \caption{\textbf{: Preconditioned pipelined BiCGStab}\\
	\small Pipelined variant of the BiCGStab algorithm consisting of 6 \dotpr s, 18 \axpy s, 2 \glred~phases, 2 \spmv s and 2 \preco~applications per iteration. Each \glred~global synchronization phase is overlapped by one corresponding \spmv~+ \preco~application as indicated by the \textbf{begin global reduction} and \textbf{end reduction} key words, representing e.g.~the MPI delimiters \texttt{MPI\_Iallreduce} and \texttt{MPI\_Wait}.}
  \label{algo::pipebicgstab}
  \begin{algorithmic}[1]
  	\Function{p-bicgstab}{$A$, $M^{-1}$, $b$, $x_0$}
    \State $r_0 := b - Ax_0$; $k_0 := M^{-1} r_0$; $w_0 := A k_0$; $m_0 := M^{-1} w_0$
		\State $t_0 := A m_0$; $\alpha_0 := \left( r_0, r_0 \right) / \left( r_0, w_0 \right)$; $\beta_0 := 0$  
    \For{$i = 0, 1, 2, \dots$}
    \State $g_i := k_i + \beta_i \left( g_{i-1} - \omega_{i-1} \ell_{i-1} \right)$ 
		\State $s_i := w_i + \beta_i \left( s_{i-1} - \omega_{i-1} z_{i-1} \right)$
		\State $\ell_i := m_i + \beta_i \left( \ell_{i-1} - \omega_{i-1} n_{i-1} \right)$
		\State $z_i := t_i + \beta_i \left( z_{i-1} - \omega_{i-1} v_{i-1} \right)$
		\State $q_i := r_i - \alpha_i s_i$
		\State $u_i := k_i - \alpha_i \ell_i$
		\State $y_i := w_i - \alpha_i z_i$
		\State \textbf{begin global reduction} $\left( q_i, y_i \right)$; $\left( y_i, y_i \right)$
		\State $n_i := M^{-1} z_i$
		\State $v_i := A n_i$
		\State \textbf{end reduction}
		\State $\omega_i := \left( q_i, y_i \right)/\left( y_i, y_i \right)$
		\State $x_{i+1} := x_i + \alpha_{i} g_i + \omega_i u_i$
    \State $r_{i+1} := q_i - \omega_{i} y_i$
		\State $k_{i+1} := u_i - \omega_i \left( m_i - \alpha_i n_i \right)$
		\State $w_{i+1} := y_i - \omega_i \left( t_i - \alpha_i v_i \right)$
		\State \textbf{begin global reduction} $\left( r_0, r_{i+1} \right)$; $\left( r_0, w_{i+1} \right)$; $\left( r_0, s_i \right)$; $\left( r_0, z_i \right)$
		\State $m_{i+1} := M^{-1} w_{i+1}$
		\State $t_{i+1} := A m_{i+1}$
		\State \textbf{end reduction}
		\State $\beta_{i+1} := \left( \alpha_i / \omega_i \right)  \left( r_0, r_{i+1} \right) / \left( r_0, r_i \right)$
		\State $\alpha_{i+1} := \left( r_0, r_{i+1} \right) / \left(  \left( r_0, w_{i+1} \right) + \beta_{i+1} \left( r_0, s_i \right) - \beta_{i+1} \omega_i \left( r_0, z_i \right) \right)$
    \EndFor
    \EndFunction
  \end{algorithmic}
\end{algorithm}

\subsection{Local rounding errors in the finite precision recurrence relations in p-BiCGStab} \label{sec:describing}

The finite precision equivalent, or \emph{actually computed} variants, of the variables defined in \eqref{eq:auxvars} will again be denoted by the bar symbol in this section.
The following perturbed recurrence relations hold for the pipelined BiCGStab algorithm in finite precision:
\begin{align}
	\bar{g}_i 		&:= \bar{k}_i + \bar{\beta}_i \left(\bar{g}_{i-1} - \bar{\omega}_{i-1} \bar{\ell}_{i-1}\right) + \dig, &
	\bar{x}_{i+1} &:= \left(\bar{x}_i + \bar{\alpha}_i \bar{g}_i\right) + \bar{\omega}_i \bar{u}_i + \dix, & 
	\bar{q}_i 		&:= \bar{r}_i - \bar{\alpha}_i \bar{s}_i + \diq, \notag \\
	\bar{s}_i 		&:= \bar{w}_i + \bar{\beta}_i \left(\bar{s}_{i-1} - \bar{\omega}_{i-1} \bar{z}_{i-1}\right) + \dis, &
	\bar{r}_{i+1} &:= \bar{q}_i - \bar{\omega}_i \bar{y}_i + \dir, &
	\bar{u}_i 		&:= \bar{k}_i - \bar{\alpha}_i \bar{\ell}_i + \diu, \notag \\
	\bar{\ell}_i 	&:= \bar{m}_i + \bar{\beta}_i \left(\bar{\ell}_{i-1} - \bar{\omega}_{i-1} \bar{n}_{i-1}\right) + \dil, &
	\bar{k}_{i+1} &:= \bar{u}_i - \bar{\omega}_i \left(\bar{m}_i - \bar{\alpha}_i \bar{n}_i\right) + \dik, &
	\bar{y}_i 		&:= \bar{w}_i - \bar{\alpha}_i \bar{z}_i + \diy, \notag \\
	\bar{z}_i 		&:= \bar{t}_i + \bar{\beta}_i \left(\bar{z}_{i-1} - \bar{\omega}_{i-1} \bar{v}_{i-1}\right) + \diz, &
	\bar{w}_{i+1} &:= \bar{y}_i - \bar{\omega}_i \left(\bar{t}_i - \bar{\alpha}_i \bar{v}_i\right) + \diw. \label{eq:auxrec} 
\end{align}
Furthermore, the following \spmv~and \preco~applications are computed explicitly in the algorithm:
\begin{align}
	\bar{n}_i &:= \text{fl}(M^{-1} \bar{z}_i), & \bar{v}_i &:= \text{fl}(A \bar{n}_i), & \bar{m}_i &:= \text{fl}(M^{-1} \bar{w}_i), & \bar{t}_i &:= \text{fl}(A \bar{m}_i). \label{eq:auxvars2}
\end{align}
It is assumed that the default order of operations in finite precision arithmetic is from left to right unless otherwise indicated by the brackets in the expressions \eqref{eq:auxrec} (while evidently always abiding by the default order of operations for multiplication and addition). 
The local rounding errors in \eqref{eq:auxrec} 
can be bounded as follows:
\begin{align}
	\|\dig\| 		&\leq \left( \lef \bar{k}_i\rig + 3\lef\bar{\beta}_i\bar{g}_{i-1}\rig + 4 \lef\bar{\beta}_i\bar{\omega}_{i-1}\bar{\ell}_{i-1}\rig\right) \epsilon, &
	\|\dix\| 		&\leq \left( 2 \lef \bar{x}_i\rig + 3\lef\bar{\alpha}_i \bar{g}_i \rig + 2 \lef\bar{\omega}_i\bar{u}_i\rig\right) \epsilon, \notag \\
	\|\dis\| 		&\leq \left( \lef \bar{w}_i\rig + 3\lef\bar{\beta}_i \bar{s}_{i-1}\rig + 4\lef\bar{\beta}_i\bar{\omega}_{i-1} \bar{z}_{i-1}\rig \right) \epsilon,  &
	\|\dir\| 		&\leq \left( \lef \bar{q}_i \rig + 2 \lef \bar{\omega}_i \bar{y}_i \rig \right) \epsilon, \label{eq:errbounds1}
\end{align}	
and
\begin{align}
	\|\dil\| 		&\leq \left( \lef \bar{m}_i \rig + \tmsn \lef M^{-1} \rig \lef \bar{w}_i \rig + 3 \lef\bar{\beta}_i \bar{\ell}_{i-1} \rig \right.  
							  + \left. 4 \lef \bar{\beta}_i\bar{\omega}_{i-1}\bar{n}_{i-1}\rig +\tmsn \lef M^{-1} \rig \lef \bar{\beta}_i\bar{\omega}_{i-1} \bar{z}_{i-1} \rig \right) \epsilon, \notag \\
	\|\diz\| 		&\leq \left( \lef \bar{t}_i \rig + \msn \lef A \rig \lef \bar{m}_i \rig +\tmsn \lef A \rig \lef M^{-1}\rig \lef \bar{w}_i\rig + 3\lef\bar{\beta}_i\bar{z}_{i-1}\rig \right. \notag \\
							& \qquad + \left. 4 \lef \bar{\beta}_i \bar{\omega}_{i-1} \bar{v}_{i-1}\rig + \msn \lef A \rig\lef \bar{\beta}_i\bar{\omega}_{i-1} \bar{n}_{i-1}\rig + \tmsn \lef A \rig \lef M^{-1}\rig \lef \bar{\beta}_i\bar{\omega}_{i-1}\bar{z}_{i-1}\rig \right) \epsilon, \label{eq:errbounds2}
\end{align}
and 
\begin{align}
	\|\dik\|	 	&\leq \left( \lef \bar{u}_i \rig + 3 \lef \bar{\omega}_i \bar{m}_i \rig + \tmsn \lef M^{-1} \rig \lef \bar{\omega}_i \bar{w}_i \rig + 4 \lef \bar{\omega}_i \bar{\alpha}_i \bar{n}_i \rig + \tmsn \lef M^{-1} \rig \lef \bar{\omega}_i \bar{\alpha}_i \bar{z}_i \rig \right) \epsilon, \notag \\
	\|\diw\| 		&\leq \left( \lef \bar{y}_i \rig + 3 \lef \bar{\omega}_i \bar{t}_i \rig + \msn \lef A \rig \lef \bar{\omega}_i \bar{m}_i \rig + \tmsn \lef A \rig \lef M^{-1} \rig \lef \bar{\omega}_i \bar{w}_i \rig \right.\notag \\
							& \qquad + \left. 4 \lef \bar{\omega}_i \bar{\alpha}_i \bar{v}_i \rig + \msn \lef A \rig \lef \bar{\omega}_i \bar{\alpha}_i \bar{n}_i \rig + \tmsn \lef A \rig \lef M^{-1} \rig \lef \bar{\omega}_i \bar{\alpha}_i \bar{z}_i \rig\right) \epsilon, \label{eq:errbounds3}
\end{align}
and
\begin{equation}						
	\|\diq\| 		\leq \left( \lef \bar{r}_i \rig + 2 \lef \bar{\alpha}_i \bar{s}_i \rig \right) \epsilon, \qquad
	\|\diu\| 		\leq \left( \lef \bar{k}_i \rig + 2 \lef \bar{\alpha}_i \bar{\ell}_i \rig \right) \epsilon, \qquad
	\|\diy\| 		\leq \left( \lef \bar{w}_i \rig + 2 \lef \bar{\alpha}_i \bar{z}_i\rig \right) \epsilon.\label{eq:errbounds4}
\end{equation}
In the expressions \eqref{eq:errbounds1}-\eqref{eq:errbounds4} $\mu$ is again the maximum number of nonzeros (nnz) over all rows of $A$ and $\tilde{\mu}$ denotes the maximum nnz over all rows of $M^{-1}$.

\subsection{The gap between the true and recursively computed residual in p-BiCGStab} \label{sec:gap}

We derive an expression for the gap between the current true residual $b-A\bar{x}_i$ and the recursively computed residual $\bar{r}_i$, see expression \eqref{eq:res_gap}, 
where $\bar{x}_{i+1}$ and $\bar{r}_{i+1}$ are defined by \eqref{eq:auxrec}.
Similar to the definition for the residual gap, the following gaps on the auxiliary variables $s_i$, $w_i$ and $z_i$ are defined: 
\begin{equation}
\Ds_i := A \bar{g}_i - \bar{s}_i, \qquad \Dw_{i+1} := A \bar{k}_{i+1} - \bar{w}_{i+1}, \qquad \Dz_i := A \bar{\ell}_i - \bar{z}_i. \label{eq:DsDwDz}
\end{equation} 
The norm of the initial residual gap, which equals the initial round-off from computing $\bar{r}_0$ from $b$ and $\bar{x}_0$, i.e., $\Dr_{0} := (b-A\bar{x}_0) - \text{fl}(b-A\bar{x}_0)$, is bounded by
$\lef\Dr_{0}\rig \leq ( (\msn+1) \lef A\rig \lef \bar{x}_0\rig + \lef b \rig ) \, \epsilon$. In iteration $i$ of the algorithm it holds by substituting $\bar{x}_{i+1}$ and $\bar{r}_{i+1}$ that
\begin{align} \label{eq:Dr}
	\Dr_{i+1}	&= (b-A\bar{x}_{i+1})-\bar{r}_{i+1} \notag \\
						&= b-A\bar{x}_i - \bar{\alpha}_i A \bar{g}_i - \bar{\omega}_i A \bar{u}_i - \bar{q}_i + \bar{\omega}_i \bar{y}_i - A \dix - \dir \notag \\
						&= b-A\bar{x}_i-\bar{\alpha}_i A \bar{g}_i -\bar{\omega}_i A \bar{k}_i + \bar{\omega}_i \bar{\alpha}_i A \bar{\ell}_i - \bar{r}_i + \bar{\alpha}_i \bar{s}_i + \bar{\omega}_i \bar{w}_i - \bar{\omega}_i \bar{\alpha}_i \bar{z}_i \notag \\
						&~~~ - A \dix - \dir  - \bar{\omega}_i A \diu - \diq + \bar{\omega}_i\diy \notag \\
						&= \Dr_i - \bar{\alpha}_i \Ds_i - \bar{\omega}_i \Dw_i +\bar{\omega}_i\bar{\alpha}_i\Dz_i - A \dix - \dir  - \bar{\omega}_i A \diu - \diq + \bar{\omega}_i\diy,
\end{align}
where the auxiliary variables $\bar{q}_i$, $\bar{u}_i$ and $\bar{y}_i$ were also substituted by their respective recurrence relations.
Note that the gap $\Dr_{i+1}$ is coupled to the gaps $\Ds_i, \Dw_i$ and $\Dz_i$ on the other auxiliary variables by the recurrence relations in the algorithm. 
We hence additionally derive expressions for the gaps on the auxiliary variables \eqref{eq:DsDwDz}.
At the start of the algorithm the gap $\Ds_i$ is bounded by 
$\lef\Ds_{0}\rig \leq \msn \lef A\rig \lef \bar{g}_0\rig \epsilon$, whereas in the $i$-th iteration it holds
\begin{align} \label{eq:Ds}
	\Ds_i	&= A\bar{g}_i - \bar{s}_i \notag \\
				&= A\bar{k}_i + \bar{\beta}_i A \bar{g}_{i-1} - \bar{\beta}_i \bar{\omega}_{i-1} A \bar{\ell}_{i-1} - \bar{w}_i - \bar{\beta}_i \bar{s}_{i-1} + \bar{\beta}_i \bar{\omega}_{i-1} \bar{z}_{i-1} + A \dig - \dis \notag \\
				&= \Dw_i + \bar{\beta}_i \Ds_{i-1} - \bar{\beta}_i \bar{\omega}_{i-1} \Dz_{i-1} + A \dig - \dis.
\end{align}
Note the coupling to the gap $\Dw_i$. For the latter gap the bound $\lef\Dw_{0}\rig \leq \msn \lef A\rig \lef \bar{k}_0\rig \epsilon$ holds when $i = 0$, while for general iteration index $i$the gap equals
\begin{align} \label{eq:Dw}
	\Dw_{i+1}	&= A \bar{k}_{i+1}-\bar{w}_{i+1} \notag \\
						&= A \bar{u}_i - \bar{\omega}_i A \bar{m}_i + \bar{\omega}_i \bar{\alpha}_i A \bar{n}_i - \bar{y}_i + \bar{\omega}_i \bar{t}_i - \bar{\omega}_i \bar{\alpha}_i \bar{v}_i + A\dik - \diw \notag \\
						&= A \bar{k}_i - \bar{\alpha}_i A \bar{\ell}_i - \bar{\omega}_i A \bar{m}_i + \bar{\omega}_i \bar{\alpha}_i A \bar{n}_i -\bar{w}_i +\bar{\alpha}_i \bar{z}_i +\bar{\omega}_i \bar{t}_i - \bar{\omega}_i \bar{\alpha}_i \bar{v}_i + A \dik - \diw + A\diu -\diy  \notag \\
						&= \Dw_i - \bar{\alpha}_i \Dz_i  + A \dik - \diw + A\diu -\diy.
\end{align}
Finally, for the auxiliary variable $z_i$ the initial gap can be bounded by $\lef\Dz_{0}\rig \leq \msn \lef A\rig \lef \bar{\ell}_0\rig \epsilon$, and the gap $\Dz_i$ in iteration $i$ is
\begin{align}  \label{eq:Dz}
	\Dz_i	&= A \bar{\ell}_i - \bar{z}_i \notag \\
				&= A \bar{m}_i + \bar{\beta}_i A \bar{\ell}_{i-1} - \bar{\beta}_i \bar{\omega}_{i-1} A \bar{n}_{i-1} - \bar{t}_i -\bar{\beta}_i \bar{z}_{i-1}+\bar{\beta}_i\bar{\omega}_{i-1}\bar{v}_{i-1} + A \dil - \diz \notag \\
				&= \bar{\beta}_i \Dz_{i-1} + A\dil - \diz.
\end{align}
By induction it follows from expressions \eqref{eq:Dr}-\eqref{eq:Dz} that the residual gap $\Dr_{i}$ can be formulated as:
\begin{equation} \label{eq:DDr}
	\Dr_{i}  =  \Dr_0 - \sum_{j=0}^{i-1} \bar{\alpha}_j \Ds_j - \sum_{j=0}^{i-1} \bar{\omega}_j \Dw_j + \sum_{j=0}^{i-1} \bar{\omega}_j \bar{\alpha}_j \Dz_j + \sum_{j=0}^{i-1} (- A \djx - \djr  - \bar{\omega}_j A \dju - \djq + \bar{\omega}_j\djy) ,
\end{equation}
where
\begin{equation} \label{eq:DDs}
	\Ds_i  =  \left(\prod_{j=1}^{i} \bar{\beta}_j\right) \Ds_0 + \sum_{j=1}^i \left(\prod_{k=j+1}^{i}\bar{\beta}_k\right) \Dw_j - \sum_{j=0}^{i-1}\left(\prod_{k=j+1}^{i}\bar{\beta}_k\right) \bar{\omega}_j \Dz_j + \sum_{j=1}^i \left(\prod_{k=j+1}^{i} \bar{\beta}_k\right) (A \djq - \djs) ,
\end{equation}
with 
\begin{equation} \label{eq:DDw} 
	\Dw_{i}  =  \Dw_0 - \sum_{j=0}^{i-1} \bar{\alpha}_j \Dz_j + \sum_{j=0}^{i-1} (A \djk - \djw + A\dju -\djy) ,
\end{equation}
and
\begin{equation} \label{eq:DDz}
  \Dz_i  =  \left( \prod_{j=1}^{i} \bar{\beta}_j \right) \Dz_0 +  \sum_{j=1}^{i} \left( \prod_{k=j+1}^{i} \bar{\beta}_k \right) (A\djl - \djz) .
\end{equation}
Rewriting the recursions \eqref{eq:Dr}-\eqref{eq:Dz} as a function of the gaps $\Dr_i$, $\Ds_{i-1}$, $\Dw_i$ and $\Dz_{i-1}$ in iteration $i-1$, we arrive alternatively at the following coupled system of equations:
\begin{equation} \label{eq:system}
\begin{bmatrix}
 \Dr_{i+1} \\
 \Ds_i \\
 \Dw_{i+1} \\ 
 \Dz_i
\end{bmatrix} = 
\begin{bmatrix}
    1 & -\bar{\alpha}_i \bar{\beta}_i & -\left(\bar{\alpha}_i+\bar{\omega}_i\right) & \bar{\alpha}_i\bar{\beta}_i\left(\bar{\omega}_i + \bar{\omega}_{i-1}\right) \\ 
		0 & \bar{\beta}_i & 1 & -\bar{\beta}_i\bar{\omega}_{i-1} \\ 
		0 & 0 & 1 & -\bar{\alpha}_i \bar{\beta}_i \\ 
		0 & 0 & 0 & \bar{\beta}_i
\end{bmatrix}
\begin{bmatrix}
 \Dr_i \\
 \Ds_{i-1} \\
 \Dw_i \\ 
 \Dz_{i-1}
\end{bmatrix} +
\begin{bmatrix}
  \epsilon_i^r \\
	\epsilon_i^s \\
	\epsilon_i^w \\
	\epsilon_i^z 
\end{bmatrix},
\end{equation}
where the local error additions in each iteration are given by
\begin{equation} \label{eq:local}
\begin{bmatrix}
  \epsilon_i^r \\
	\epsilon_i^s \\
	\epsilon_i^w \\
	\epsilon_i^z 
\end{bmatrix} = 
\begin{bmatrix}
  -A\dix -\bar{\omega}_i A\diu +\bar{\omega}_i\diy -\dir -\diq -\bar{\alpha}_i \left( A \dig - \dis \right) +\bar{\omega}_i \bar{\alpha}_i \left( A \dil - \diz \right) \\
	 A \dig - \dis\\
	-\bar{\alpha}_i \left( A \dil - \diz \right) + A \dik - \diy - \diw + A \diu \\
	 A\dil - \diz
\end{bmatrix}.
\end{equation}

\subsection{Further analysis: expressing the residual gap for p-BiCGStab in matrix notation} \label{sec:expressing}

In this section the expressions \eqref{eq:DDr}-\eqref{eq:DDz} for the gaps are reformulated in a more comprehensible matrix notation, which grants more insight into the propagation of local rounding errors throughout the algorithm.
Let $B = [b,b,\ldots,b]$, $\bar{X}_{i+1} = [\bar{x}_0,\bar{x}_1,\ldots, \bar{x}_i]$ and $\bar{P}_{i+1} = [\bar{p}_0,\bar{p}_1,\ldots, \bar{p}_i]$ be defined. Writing the expressions for the gaps in matrix notation yields
\begin{align}
  \mathcal{R}_{i+1} &= [\Delta^r_0,\Delta^r_1,\ldots, \Delta^r_i] = (B-A\bar{X}_{i+1})-\bar{R}_{i+1} , & 	&\text{~with~~}& \bar{R}_{i+1} &= [\bar{r}_0,\bar{r}_1,\ldots, \bar{r}_i], \notag \\
  \mathcal{S}_{i+1} &= [\Delta^s_0,\Delta^s_1,\ldots, \Delta^s_i] = A\bar{P}_{i+1}-\bar{S}_{i+1} , & 			&\text{~with~~}& \bar{S}_{i+1} &= [\bar{s}_0,\bar{s}_1,\ldots, \bar{s}_i], &  \notag \\
  \mathcal{W}_{i+1} &= [\Delta^w_0,\Delta^w_1,\ldots, \Delta^w_i] = A\bar{R}_{i+1}-\bar{W}_{i+1} , & 			&\text{~with~~}& \bar{W}_{i+1} &= [\bar{w}_0,\bar{w}_1,\ldots, \bar{w}_i], &  \notag \\
  \mathcal{Z}_{i+1} &= [\Delta^z_0,\Delta^z_1,\ldots, \Delta^z_i] = A\bar{S}_{i+1}-\bar{Z}_{i+1} , & 			&\text{~with~~}& \bar{Z}_{i+1} &= [\bar{z}_0,\bar{z}_1,\ldots, \bar{z}_i]. &  \notag
\end{align}
Matrix expressions for the local rounding errors on the auxiliary variables in p-BiCGStab are introduced as follows:
\begin{align}
 \Theta_{i+1}^{x} &= [0,-\delta_1^{x},\ldots, -\delta_i^{x}], & \Theta_{i+1}^{r} &= [\Delta_0^r,-\delta_1^{r},\ldots, -\delta_i^{r}], & \Theta_{i+1}^{u} &= [0,\delta_0^{u},\ldots, \delta_{i-1}^{u}],\notag \\
 \Theta_{i+1}^{q} &= [0,\delta_1^{q},\ldots, \delta_i^{q}],   & \Theta_{i+1}^{s} &= [\Delta_0^s,-\delta_1^{s},\ldots, -\delta_i^{s}], & \Theta_{i+1}^{p} &= [0,-\delta_0^{q},\ldots, -\delta_{i-1}^{q}],\notag \\
 \Theta_{i+1}^{k} &= [0,\delta_1^{k},\ldots, \delta_i^{k}],   & \Theta_{i+1}^{w} &= [\Delta_0^w,-\delta_1^{w},\ldots, -\delta_i^{w}], & \Theta_{i+1}^{y} &= [0,-\delta_0^{y},\ldots, -\delta_{i-1}^{y}],\notag \\
 \Theta_{i+1}^{l} &= [0,\delta_1^{l},\ldots, \delta_i^{l}],   & \Theta_{i+1}^{z} &= [\Delta_0^z,-\delta_1^{z},\ldots, -\delta_i^{z}]. & \label{eq:thetas}
\end{align}
By unreeling the expressions \eqref{eq:DDr}-\eqref{eq:DDz} the following matrix expressions for the local rounding errors are then obtained: 
\begin{align} 
  \mathcal{R}_{i+1} &= (A\Theta_{i+1}^{x}+\Theta_{i+1}^{r}+A\Theta_{i+1}^{u}\mathcal{O}_{i+1}+\Theta_{i+1}^{p}+\Theta_{i+1}^{y}\mathcal{O}_{i+1}) \, \mathcal{U}_{i+1} + \ldots \notag \\
										& \qquad \qquad \qquad \qquad \qquad \qquad \qquad + \mathcal{S}_{i+1} \mathcal{A}_{i+1} + \mathcal{W}_{i+1} \mathcal{C}_{i+1} + \mathcal{Z}_{i+1} \mathcal{D}_{i+1}, \notag \\ 
  \mathcal{S}_{i+1} &= (A\Theta_{i+1}^{q}+\Theta_{i+1}^{s}) \, \mathcal{B}_{i+1} + \mathcal{W}_{i+1} \mathcal{E}_{i+1} + \mathcal{Z}_{i+1} \mathcal{P}_{i+1}, \notag \\ 
  \mathcal{W}_{i+1} &= (A\Theta_{i+1}^{k}+\Theta_{i+1}^{w}+A\Theta_{i+1}^{u}+\Theta_{i+1}^{y}) \, \mathcal{U}_{i+1} + \mathcal{Z}_{i+1} \mathcal{A}_{i+1}, \notag \\ 
  \mathcal{Z}_{i+1} &= (A\Theta_{i+1}^{l}+\Theta_{i+1}^{z}) \, \mathcal{B}_{i+1}, \label{eq:matrix_expr}
\end{align}
where the matrices $\mathcal{A}_{i+1}$, $\mathcal{B}_{i+1}$, $\mathcal{E}_{i+1}$, $\mathcal{P}_{i+1}$, $\mathcal{C}_{i+1}$, $\mathcal{D}_{i+1}$, $\mathcal{O}_{i+1}$ and $\mathcal{U}_{i+1}$, which are referred to as \emph{error propagation matrices} in this work, are defined as
\begin{align}
  \mathcal{A}_{i+1} &= -
  \left(\begin{array}{ccccc} 
    0&\bar{\alpha}_0&\bar{\alpha}_0&\cdots&\bar{\alpha}_0 \\
    &0&\bar{\alpha}_1&\cdots&\bar{\alpha}_1 \\
    &&\ddots&& \vdots\\
    &&&0& \bar{\alpha}_{i-1}\\
    &&&&0
  \end{array}\right),
&
\mathcal{B}_{i+1} &= 
\left(\begin{array}{ccccc} 
1&\bar{\beta}_1&\bar{\beta}_1 \bar{\beta}_2&\cdots&\bar{\beta}_1 \bar{\beta}_2 \ldots \bar{\beta}_i\\
&1&\bar{\beta}_2& \ddots &\bar{\beta}_2 \ldots \bar{\beta}_i\\
&&\ddots&\ddots&\vdots \\
&&&\ddots&\bar{\beta}_i \\
&&&&1
\end{array}\right),
\notag \\
  \mathcal{E}_{i+1} &= 
\left(\begin{array}{ccccc} 
0&0&0&\cdots&0\\
&1&\bar{\beta}_2& \ddots &\bar{\beta}_2 \ldots \bar{\beta}_i\\
&&\ddots&\ddots&\vdots \\
&&&\ddots&\bar{\beta}_i \\
&&&&1
\end{array}\right),
&
\mathcal{P}_{i+1} &= -
\left(\begin{array}{ccccc} 
0&\bar{\omega}_0 \bar{\beta}_1&\bar{\omega}_0 \bar{\beta}_1 \bar{\beta}_2&\cdots&\bar{\omega}_0 \bar{\beta}_1 \bar{\beta}_2 \ldots \bar{\beta}_i \\
&0&\bar{\omega}_1 \bar{\beta}_2& \ddots &\bar{\omega}_1 \bar{\beta}_2 \ldots \bar{\beta}_i\\
&&\ddots&\ddots&\vdots \\
&&&\ddots&\bar{\omega}_{i-1} \bar{\beta}_i \\
&&&&0
\end{array}\right),
\notag \\
  \mathcal{C}_{i+1} &= -
\left(\begin{array}{ccccc} 
0&\bar{\omega}_0 &\bar{\omega}_0 &\cdots&\bar{\omega}_0 \\
&0&\bar{\omega}_1 & \ddots &\bar{\omega}_1\\
&&\ddots&\ddots&\vdots \\
&&&\ddots&\bar{\omega}_{i-1} \\
&&&&0
\end{array}\right),
&
\mathcal{D}_{i+1} &= 
\left(\begin{array}{ccccc} 
0&\bar{\omega}_0 \bar{\alpha}_0 & \bar{\omega}_0 \bar{\alpha}_0 &\cdots&\bar{\omega}_0 \bar{\alpha}_0 \\
&0& \bar{\omega}_1 \bar{\alpha}_1 & \ddots & \bar{\omega}_1 \bar{\alpha}_1\\
&&\ddots&\ddots&\vdots \\
&&&\ddots&\bar{\omega}_{i-1} \bar{\alpha}_{i-1} \\
&&&&0
\end{array}\right),
\notag \\
  \mathcal{O}_{i+1} &= -
\left(\begin{array}{ccccc} 
0& & & & \\
&\bar{\omega}_0& & & \\
&&\bar{\omega}_1& & \\
&&&\ddots& \\
&&&&\bar{\omega}_{i-1}
\end{array}\right),
&
\mathcal{U}_{i+1} &= 
\left(\begin{array}{ccccc} 
 1 & 1 &\cdots&1&1 \\
   & 1 &&&1 \\
& & \ddots&&\vdots \\
&&&1&1 \\
&&&&1
\end{array}\right).
\label{eq:all_matrices}
\end{align}
The entries of these error propagation matrices determine the amplification of local rounding errors throughout the p-BiCGStab algorithm. Indeed, when the modulus of an entry is large, the corresponding error(s) in expression \eqref{eq:matrix_expr} can be expected to be amplified. Note that, as indicated by the numerical results in Section \ref{sec:numerical}, the entries of the error propagation matrices can easily be several orders of magnitude larger than one in practice.

Breaking down expression \eqref{eq:matrix_expr}, all of the following products of error propagation matrices are multiplied to the left by matrices $\Theta^{\cdot}_{i+1}$ containing local rounding errors:
\begin{equation} \label{eq:all_error_matrices}
		\mathcal{U}_i, \quad \mathcal{O}_i \mathcal{U}_i, \quad \mathcal{B}_i \mathcal{A}_i, \quad \mathcal{U}_i \mathcal{E}_i \mathcal{A}_i, \quad  \mathcal{B}_i \mathcal{A}_i \mathcal{E}_i \mathcal{A}_i, \quad \mathcal{B}_i \mathcal{P}_i \mathcal{A}_i, \quad \mathcal{U}_i \mathcal{C}_i, \quad \mathcal{B}_i \mathcal{A}_i \mathcal{C}_i, \quad \mathcal{B}_i \mathcal{D}_i.
\end{equation}
When the modulus of an entry in one of these matrix products becomes excessively large, numerical stability is typically lost. Hence, the maximum norms of the columns of these matrices provide an important indication on the numerical stability of the algorithm, see Section \ref{sec:numerical}, Fig.~\ref{fig:comparison}--\ref{fig:comparison_ichol_rr}.

\subsection{Comparison to the stability analysis for the pipelined Conjugate Gradient method} \label{sec:comparison_to_CG}

Some comments on the comparison of the analysis in Sections \ref{sec:gap} and \ref{sec:expressing} to the numerical analysis of the closely related pipelined Conjugate Gradients method (p-CG) published in \cite{cools2018analyzing} are in order here. Expressions \eqref{eq:system}-\eqref{eq:local} are the analogue to the relations that constitute the numerical analysis of the pipelined CG algorithm, see \cite{cools2018analyzing}, Section 2.3, p.\,434, expression (2.37). 

Note that when setting $\omega_0 = \omega_1 = \ldots = \omega_i = 0$, the analysis of the p-BiCGStab method simplifies to that of the p-CG method. Expression \eqref{eq:system} then simplifies to the aforementioned relation (2.37) from \cite{cools2018analyzing} that was established for the p-CG method, conforming to intuition. The numerical analysis in this study can thus be considered as a generalization of the earlier work in \cite{cools2018analyzing}.

However, the $4$-by-$4$ operator that describes the iteration-wise propagation of local rounding errors in \eqref{eq:system} for p-BiCGStab contains significantly more non-zero entries compared to the numerical analysis of p-CG. In addition, the BiCGStab algorithm propagates more contributions of local error terms as indicated by the error expressions \eqref{eq:local}. Hence, it is expected that the pipelined BiCGStab method is significantly more sensitive to rounding error propagation compared to pipelined CG. 
This conjecture is validated by the matrix expressions in Section \ref{sec:expressing}. In the numerical stability analysis of the pipelined CG method only the matrices $\mathcal{U}_i$, $\mathcal{B}_i \mathcal{A}_i$, $\mathcal{U}_i \mathcal{E}_i \mathcal{A}_i$ and $\mathcal{B}_i \mathcal{A}_i \mathcal{E}_i \mathcal{A}_i$ (using the notations used in the current manuscript\footnote{In the numerical stability analysis of the pipelined CG algorithm in \cite{cools2018analyzing} the matrices $\mathcal{A}_i$ and $\mathcal{B}_i$ are defined identically to the definitions given in the current manuscript, see \eqref{eq:matrix_expr}. However, to avoid notational confusion it should be noted that the scalar entries $\alpha_i$ and $\beta_i$ are defined differently in the CG and BiCGStab algorithms.}) appear in the expression for the residual gap, see  \cite{cools2018analyzing}.

\subsection{The gap on the preconditioned residual in preconditioned p-BiCGStab}

Apart from the residual $\bar{r}_i$ the pipelined BiCGStab algorithm also computes the preconditioned residual $\bar{k}_i$ recursively.
However, as shown by \eqref{eq:DDr}-\eqref{eq:DDz}, the residual gap $\Dr_{i+1}$ is not coupled to the preconditioned residual gap $\Dk_{i+1} := M^{-1}\bar{r}_{i+1} - \bar{k}_{i+1}$. The gap $\Dk_{i+1}$ satisfies the following recurrence in iteration $i$:
\begin{align} \label{eq:Dk}
	\Dk_{i+1}&= M^{-1} \bar{r}_{i+1} - \bar{k}_{i+1} \notag \\
						&= M^{-1} \bar{q}_i - \bar{\omega}_i M^{-1} \bar{y}_i - \bar{u}_i + \bar{\omega}_i \bar{m}_i - \bar{\omega}_i\bar{\alpha}_i \bar{n}_i + M^{-1} \dir - \dik \notag \\
						&= M^{-1} \bar{r}_i -\bar{\alpha}_i M^{-1} \bar{s}_i - \bar{\omega}_i M^{-1} \bar{w}_i + \bar{\omega}_i \bar{\alpha}_i M^{-1} \bar{z}_i - \bar{k}_i + \bar{\alpha}_i \bar{\ell}_i + \bar{\omega}_i \bar{m}_i - \bar{\omega}_i \bar{\alpha}_i \bar{n}_i \notag \\
						&  ~~~ + M^{-1} \dir - \dik + M^{-1} \diq - \diu - \bar{\omega}_i M^{-1} \diy \notag \\
						&= \Dk_i - \bar{\alpha}_i \Dl_i + M^{-1} \dir - \dik + M^{-1} \diq - \diu - \bar{\omega}_i M^{-1} \diy,
\end{align}
where for $\Dl_i = M^{-1} \bar{s}_i - \bar{\ell}_i$ the following equation holds:
\begin{align} \label{eq:Dl}
	\Dl_i &= M^{-1}\bar{s}_i - \bar{\ell}_i \notag \\
				 &= M^{-1} \bar{w}_i + \bar{\beta}_i M^{-1} \bar{s}_{i-1} -\bar{\beta}_i \bar{\omega}_{i-1} M^{-1} \bar{z}_{i-1} - \bar{m}_i - \bar{\beta}_i \bar{\ell}_{i-1} + \bar{\beta}_i \bar{\omega}_{i-1} \bar{n}_{i-1} + M^{-1} \dis - \dil \notag \\
				 &= \bar{\beta}_i \Dl_{i-1} + M^{-1}  \dis - \dil.
\end{align}
It follows from \eqref{eq:Dk}-\eqref{eq:Dl} using induction that the preconditioned residual gap $\Dk_{i}$ can be formulated as:
\begin{equation} \label{eq:DDk}
	\Dk_{i}  =  \Dk_0 - \sum_{j=0}^{i-1} \bar{\alpha}_j \Dl_j + \sum_{j=0}^{i-1} (M^{-1} \djr - \djk + M^{-1} \djq -\dju - \bar{\omega}_j M^{-1} \djy) ,
\end{equation}
where
\begin{equation} \label{eq:DDl}
	\Dl_i  =  \left(\prod_{j=1}^{i} \bar{\beta}_j\right) \Dl_0 + \sum_{j=1}^i \left(\prod_{k=j+1}^{i} \bar{\beta}_k\right) (M^{-1} \djs - \djl) .
\end{equation}
The expressions \eqref{eq:Dk}-\eqref{eq:Dl} alternatively yield the system of coupled equations
\begin{equation} \label{eq:system2}
\begin{bmatrix}
 \Dk_{i+1} \\
 \Dl_i 
\end{bmatrix} = 
\begin{bmatrix}
    1 & -\bar{\alpha}_i \bar{\beta}_i  \\ 
		0 & \bar{\beta}_i   
\end{bmatrix}
\begin{bmatrix}
 \Dk_i \\
 \Dl_{i-1} 
\end{bmatrix} +
\begin{bmatrix}
  \epsilon_i^k \\
	\epsilon_i^\ell
\end{bmatrix},
\end{equation}
where the local error additions are
\begin{equation} \label{eq:local2}
\begin{bmatrix}
  \epsilon_i^k \\
	\epsilon_i^\ell
\end{bmatrix} = 
\begin{bmatrix}
  M^{-1} \dir - \dik + M^{-1} \diq - \diu - \bar{\omega}_i M^{-1} \diy - \bar{\alpha}_i M^{-1} \dis + \bar{\alpha}_i \dil \\
	M^{-1} \dis - \dil
\end{bmatrix}.
\end{equation}
Writing the gaps \eqref{eq:DDk}-\eqref{eq:DDl} in matrix notation, it holds that
\begin{align}
  \mathcal{K}_{i+1} &= [\Delta^k_0,\Delta^k_1,\ldots, \Delta^k_i] = M^{-1}\bar{R}_{i+1}-\bar{K}_{i+1} , & 	&\text{~with~~}& \bar{K}_{i+1} &= [\bar{k}_0,\bar{k}_1,\ldots, \bar{k}_i], \notag \\
  \mathcal{L}_{i+1} &= [\Delta^\ell_0,\Delta^\ell_1,\ldots, \Delta^\ell_i] = M^{-1}\bar{S}_{i+1}-\bar{L}_{i+1} , & 	&\text{~with~~}&  \bar{L}_{i+1} &= [\bar{\ell}_0,\bar{\ell}_1,\ldots, \bar{\ell}_i]. &  \notag
\end{align}
Let the local rounding error matrices now be redefined in analogy to \eqref{eq:thetas} (replacing the earlier definitions) as follows:
\begin{align}
 \Theta_{i+1}^{r} &= [0,\delta_1^{r},\ldots, \delta_i^{r}], & \Theta_{i+1}^{k} &= [\Delta_0^k,-\delta_1^{k},\ldots, -\delta_i^{k}], & \Theta_{i+1}^{u} &= [0,-\delta_0^{u},\ldots, -\delta_{i-1}^{u}],\notag \\
 \Theta_{i+1}^{s} &= [0,\delta_1^{s},\ldots, \delta_i^{s}], & \Theta_{i+1}^{l} &= [\Delta_0^l,-\delta_1^{l},\ldots, -\delta_i^{l}], & \Theta_{i+1}^{q} &= [0,\delta_0^{q},\ldots, \delta_{i-1}^{q}],\notag \\
  \Theta_{i+1}^{y} &= [0,\delta_0^{y},\ldots, \delta_{i-1}^{y}]. & &  & &  \notag 
\end{align}
Then the preconditioned residual gap can be expressed in matrix notation as
\begin{align} 
  \mathcal{K}_{i+1} &= (M^{-1}\Theta_{i+1}^{r}+\Theta_{i+1}^{k}+M^{-1}\Theta_{i+1}^{q}+\Theta_{i+1}^{u}+M^{-1}\Theta_{i+1}^{y}\mathcal{O}_{i+1}) \, \mathcal{U}_{i+1} + \mathcal{L}_{i+1} \mathcal{A}_{i+1} , \notag \\ 
   \mathcal{L}_{i+1} &= (M^{-1}\Theta_{i+1}^{s}+\Theta_{i+1}^{l}) \, \mathcal{B}_{i+1}, \label{eq:matrix_expr2}
\end{align}
where the error propagation matrices $\mathcal{A}_{i+1}$, $\mathcal{B}_{i+1}$, $\mathcal{O}_{i+1}$ and $\mathcal{U}_{i+1}$ are defined by \eqref{eq:all_matrices}.

Comparing to the unpreconditioned residual gap \eqref{eq:system}, the system for the preconditioned residual gap is coupled to only one auxiliary variable. 
The impact of rounding errors due to the recurrence for the preconditioned residual $k_i$ can thus be expected to be far less pronounced than for the unpreconditioned residual $r_i$.
A small subset of the matrices listed in \eqref{eq:all_error_matrices}, namely $\mathcal{U}_i, \mathcal{O}_i \mathcal{U}_i$, and $\mathcal{B}_i \mathcal{A}_i$, governs the behavior of the error on the preconditioned residual, as indicated by \eqref{eq:matrix_expr2}.

Note that in pipelined BiCGStab the norm of the preconditioned residual $\|\bar{k}_i\|$ could be used to formulate an alternative stopping criterion to comparing a given tolerance to the residual norm $\|\bar{r}_i\|$ (`default' stopping criterion). We opt to use the residual norm as the error measure to define a stopping criterion in this work as it allows for an easy comparison to the classic BiCGStab algorithm in which the preconditioned residual norm is by default not computed.

\subsection{Residual replacement strategies for pipelined BiCGStab} \label{sec:rr}

The principal idea behind residual replacement type approaches to improve maximal attainable accuracy has been introduced by Sleijpen and Van der Vorst \cite{sleijpen1996reliable}. We summarize the general technique below and refer the reader to the related literature \cite{greenbaum1997estimating,van2000residual,sleijpen2001differences,carson2014residual} for additional information. 

The replacement technique 
for pipelined BiCGStab resets the residuals $r_i$ and $k_i$, as well as the auxiliary variables $w_i$, $s_i$, $l_i$ and $z_i$, to their true values in certain iterations of the algorithm. This means the following quantities are explicitly (re-)computed in selected iterations:
\begin{equation} 
	\bar{r}_i :=  b-A\bar{x}_i,	\quad	\bar{k}_i := M^{-1}\bar{r}_i, \quad	 \bar{w}_i := A\bar{k}_i, \quad
	\bar{s}_i :=  A \bar{g}_i, 	\quad	\bar{l}_i := M^{-1}\bar{s}_i, \quad	 \bar{z}_i := A\bar{l}_i. \label{eq:rr}
\end{equation}
In practice the replacements are performed in p-BiCGStab 
after executing line 20 in Algorithm \ref{algo::pipebicgstab}.
Since these replacements induce an additional computational cost (i.e.~the \textsc{spmv}s and preconditioner applications in the expressions \eqref{eq:rr}), the number of iterations in which replacements are performed should ideally be small with respect to the total number of iterations. In addition, note that in finite precision the explicitly recomputed variables may not be orthogonal to previously (recursively) constructed vectors. This may induce unwanted loss of basis orthogonality and could lead to delayed convergence \cite{strakovs2002error}. 
We refer to the literature \cite{greenbaum1992predicting,gutknecht2000accuracy} for a more detailed discussion on this topic and redirect the reader to the numerical evidence in Section \ref{sec:numerical} of the current work for additional insights.

In this work we propose two approaches to residual replacement in pipelined BiCGStab. The first methodology consists of a periodical explicit computation of the variables in expressions \eqref{eq:rr} based on a preset period. In this `periodic setting' replacements are done until the criterion $\|\bar{r}_i\| < \sqrt{\epsilon} \, \|\bar{r}_0\|$ is satisfied, at which point performing additional replacements would be detrimental for convergence \cite{strakovs2002error}. However, choosing a suitable fixed replacement period may prove challenging in practice when no prior information on the problem is available. 

Alternatively, a criterion that allows for fully automated residual replacement in p-BiCGStab, similar to the criterion suggested in \cite{cools2018analyzing} for the pipelined CG method, can be derived. From the expressions for the gaps \eqref{eq:DDr}--\eqref{eq:DDz}, or the equivalent matrix expression \eqref{eq:matrix_expr}, an upper bound on the norm of the residual gap $\Dr_{i}$ is established. We propose the following bound $f^r_i$ on the residual gap:
\begin{align} \label{eq:DDr_bound}
	\|\Dr_{i}\|  ~\leq~ f^r_i ~:=~  \|\Dr_0\| &+ \sum_{j=0}^{i-1} |\bar{\alpha}_j| \|\Ds_j\| + \sum_{j=0}^{i-1} |\bar{\omega}_j| \|\Dw_j\| + \sum_{j=0}^{i-1} |\bar{\omega}_j| |\bar{\alpha}_j| \|\Dz_j\| + \ldots \notag \\
	&+ \sum_{j=0}^{i-1} ( \| A \| \|\djx\| + \|\djr\|  + |\bar{\omega}_j| \|A\| \|\dju\| + \|\djq\| + |\bar{\omega}_j|\|\djy\| ),
\end{align}
where the quantities $\|\Ds_j\|$, $\|\Dw_j\|$ and $\|\Dz_j\|$ can be bounded in a similar fashion based on expression \eqref{eq:DDs}, \eqref{eq:DDw} and \eqref{eq:DDz} respectively. The bound \eqref{eq:DDr_bound} is obviously not tight, but captures the behavior of the residual norm fairly well as demonstrated by the numerical results in Section \ref{sec:numerical}. Alternatively, an estimate for the residual gap could be computed by adding a square root to specific terms in the upper bound \eqref{eq:DDr_bound} in analogy to the estimates proposed in \cite{cools2018analyzing} based on the typical ratio between worst-case vs.~actual floating point errors \cite{higham2002accuracy}.

The computation of the bound(s) pours down to summing the moduli of the entries in the $i$-th column of the matrices in \eqref{eq:all_error_matrices} and multiplying by the upper bounds on the corresponding local rounding error vectors given by the expressions \eqref{eq:errbounds1}-\eqref{eq:errbounds4}.
This implies the bound $f^r_i$ can easily be computed at run-time in the p-BiCGStab algorithm, since the entries of the error propagation matrices are known and the global reduction required for calculating the norms of the local rounding error vectors can be combined with the second global reduction phase in Algorithm \ref{algo::pipebicgstab}. Thus, adding the computation of the residual gap bound $f^r_i$ to the p-BiCGStab algorithm causes no overhead (apart from computing the additional local dot-product contributions). 

Given that the bound $f^r_i$ is computed in each iteration of Algorithm \ref{algo::pipebicgstab} the following heuristic criterion for automated replacements is proposed. A replacement is performed in step $i$ of Alg.~\ref{algo::pipebicgstab} 
\begin{equation}\label{eq:criterion}
	 f^r_{i-1}  \leq \tau \|\bar{r}_{i-1}\|  \quad  \text{and}  \quad   f^r_{i}  > \tau \|\bar{r}_{i}\|,
\end{equation}
where typically $\tau = \sqrt{\epsilon}$ is used. This criterion ensures that replacements are only performed when $\|\bar{r}_i\|$ is sufficiently large with respect to $\|f_i\|$ and that no excess replacements are performed. For more details on automated residual replacement strategies we refer to the closely related work \cite{cools2018analyzing}, 
specifically Section 3 and the references therein.

\section{Numerical results} \label{sec:numerical}


A set of five different test problems with various spectral properties is introduced to validate the analysis and exemplify the convergence of the BiCGStab algorithms considered in this work. The first test problem (TP1) is a small academic test case that is included here for the purpose of numerical validation and model reference only. Test problems (TP2), (TP3), (TP4) and (TP5) represent larger and/or more realistic application-driven model problems for which performance experiments on parallel hardware are included, see Section \ref{sec:performance}. Details on all test problems are given below and are summarized in Table \ref{tab:specs}.
For all numerical experiments presented in this work the system matrix $A$ is rescaled to have a 2-norm of one. This common manipulation is performed to standardize the size of the coefficients $\alpha_i$, $\beta_i$ and $\omega_i$ in the algorithms over different experiments, since these coefficients are scaled by the matrix norm by default. The solution to the original (unscaled) system can be obtained through division 
by the original matrix norm $\|A\|_2$.

\vspace{0.2cm}

\textbf{Test problem 1: Symmetric positive definite 5-point 2D stencil.} [Results shown in Figs.~\ref{fig:comparison_ichol}, \ref{fig:comparison_ichol_rr} \& \ref{fig:symmetric_and_unsymmetric}.] The first model problem is a linear system $Ax = b$ that stems from a second order central uniform finite difference discretization (FDD) of the 2D Poisson equation on the unit square with homogeneous Dirichlet boundary conditions and right-hand side $b = Ax_{ex}$, where the exact solution is $x_{ex} = \bar{\bold{1}}/\sqrt{N}$. The number of grid points in each spatial direction is $n_x = n_y = 200$, leading to a system matrix $A$ of size $N = 40,000$. The initial guess for all methods is all-zero $\bar{x}_0 = \bar{\bold{0}}$. The matrix stencil for (TP1) is shown in Table \ref{tab:specs}.
Note that solving this simple system with symmetric matrix $A_1$ does not require the use of the BiCGStab method. Conjugate Gradients would likely be the Krylov subspace method of choice for solving (TP1). The problem is considered here as a purely demonstrative benchmark to compare the numerical behavior of rounding errors in the (pipelined) CG and BiCGStab methods.

\begin{table}
\begin{center}
\begin{tabular}{|c|c|c|c|c|c|c|c|c|}
	\hline
	ID & Type & Sym. & Pos. & $n_x \times n_y \left( \times n_z \right)$ & Stencil & $\varepsilon$ & Precond. 
	\\
  \hline
  (TP1) & 2D 5pt & Yes & Yes & 200 $\times$ 200 & 
			{\small $\begin{matrix} 
			& -1  &  \\
			-1 & ~~4 & -1 \\
			& -1  &  
			\end{matrix}$ }
			& -- & ICC(0) 
	\\ 
	\hline
  (TP2) & 2D 5pt & No & Yes & 1,000 $\times$ 1,000 & 
			{\small $\begin{matrix}
			& -1  &  \\
			-1 & ~~4 & -1+\varepsilon \\
			& -1+\varepsilon  &  
			\end{matrix}$}   
			& $1$e-$3$ & -- 
	\\
  \hline
	  (TP3) & 2D 5pt & Yes & No & 500 $\times$ 500 & 
			{\small $\begin{matrix}
			& -1  &  \\
			-1 & ~~4-\varepsilon & -1 \\
			& -1  &  
			\end{matrix}$} 
			& $5$e-$4$ & ICC(0)  
	\\
	\hline 
		  (TP4) & 2D 9pt & Yes & Yes & 200 $\times$ 200 & 
			{\small $\begin{matrix}
			-1& -4  &-1  \\
			-4 & ~~20 & -1 \\
			-1& -4  &-1  
			\end{matrix}$}
			& -- & ICC(0) 
	\\
	\hline 
		  (TP5) & 3D 7pt & Yes & No & 50 $\times$ 50 $\times$ 50 & 
			{\small $\begin{matrix}
			& -1  & \hspace{-0.3cm}-1~~~~\\
			-1 & ~~6-\varepsilon & -1 \\
			~~~~{-1}\hspace{-0.3cm} & -1  &  
			\end{matrix}$} 
			& $1$e-$2$ & ICC(0) 
	\\
	\hline 
\end{tabular}
\end{center}
\vspace{-0.3cm}
\caption{Test problem specifications. Columns (from left to right): test problem ID, test problem stencil type, symmetric (Y/N), positive definite (Y/N), matrix size $N = n_x\times n_y ~(\times~ n_z)$, matrix stencil, stencil parameter value, and preconditioner type.}
\label{tab:specs}
\end{table}

\vspace{0.2cm}

\textbf{Test problem 2: Unsymmetric positive definite 5-point 2D stencil.} [Results shown in Figs.~\ref{fig:unsymmetric}, \ref{fig:timings2} \& \ref{fig:symmetric_and_unsymmetric}.] The model problem for the second benchmark case is a perturbed 2D FDD of the Poisson matrix with slightly increased upper triangular entries. The matrix stencil for this test problem is given in Table \ref{tab:specs}.
The one-dimensional size of the problem is $n_x = n_y = 1,000$; hence the matrix size is $N = 1,000,000$. The corresponding matrix $A$ is positive definite but unsymmetric, and is normalized by its 2-norm as indicated above. The right-hand side $b$ and initial guess $\bar{x}_0$ are defined similarly to (TP1). 

\vspace{0.2cm}

\textbf{Test problem 3: Indefinite symmetric 5-point 2D stencil.} [Results shown in Figs.~\ref{fig:indefinite} \& \ref{fig:timings3}.] The third benchmark problem is a (normalized) FDD discretized Helmholtz operator of size $n_x = n_y = 500$ in each spatial dimension, such that the matrix size is $N = 250,000$. 
The matrix stencil for (TP3) is given in Table \ref{tab:specs}. This shifted Poisson matrix is indefinite and severely ill-conditioned; it has a small number of negative eigenvalues (namely the 15 smallest eigenvalues of the corresponding Laplace operator with $\varepsilon = 0$) and a large number of eigenvalues very close to the origin. As such, the condition number for this model problem is very large and solving the problem is particularly challenging from a Krylov subspace solver perspective. The right-hand side and initial guess are again defined analogously to (TP1). 

\vspace{0.2cm}

\begin{figure}[t]
\begin{center}
\includegraphics[width=0.49\textwidth]{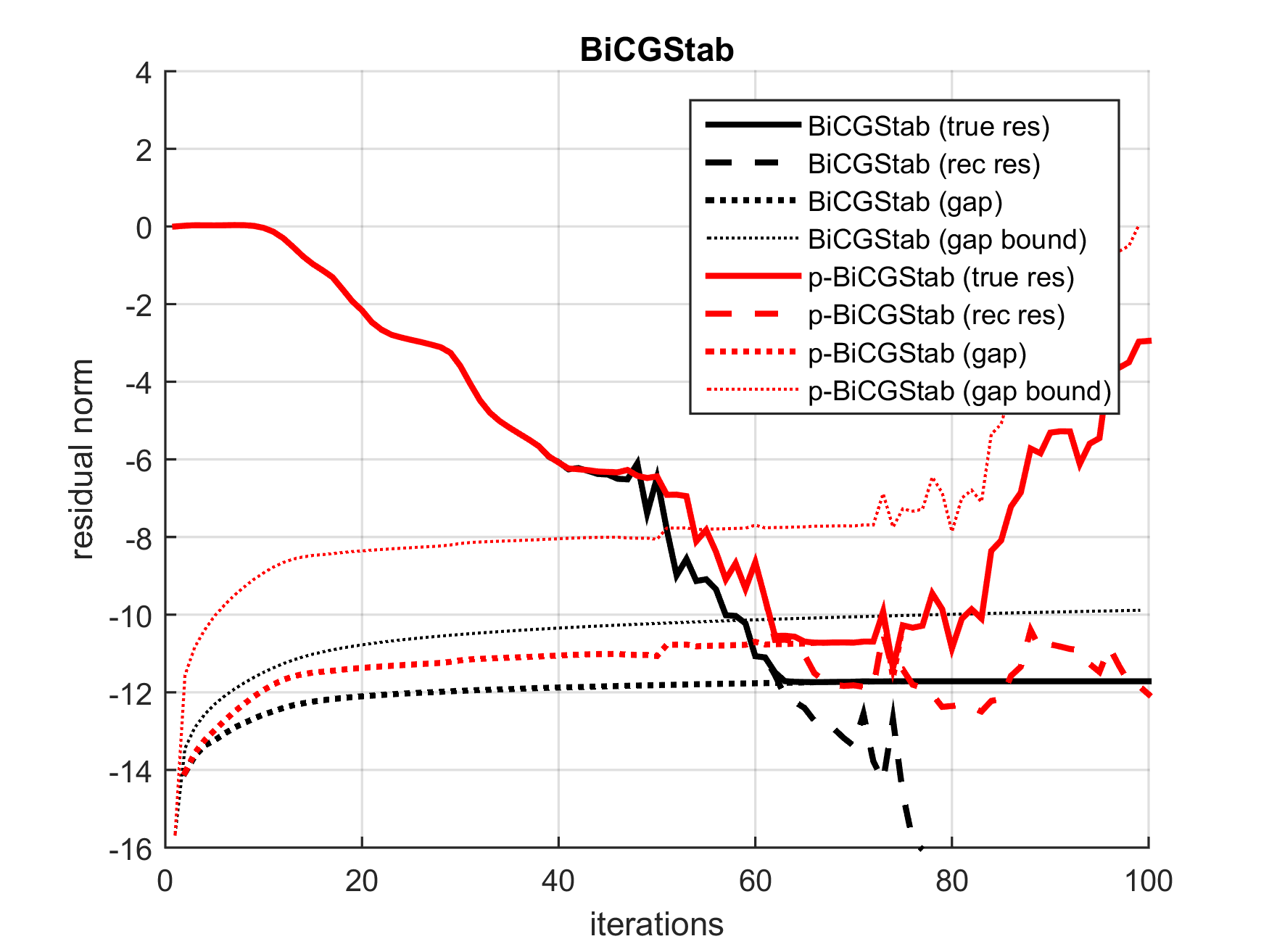}
\includegraphics[width=0.49\textwidth]{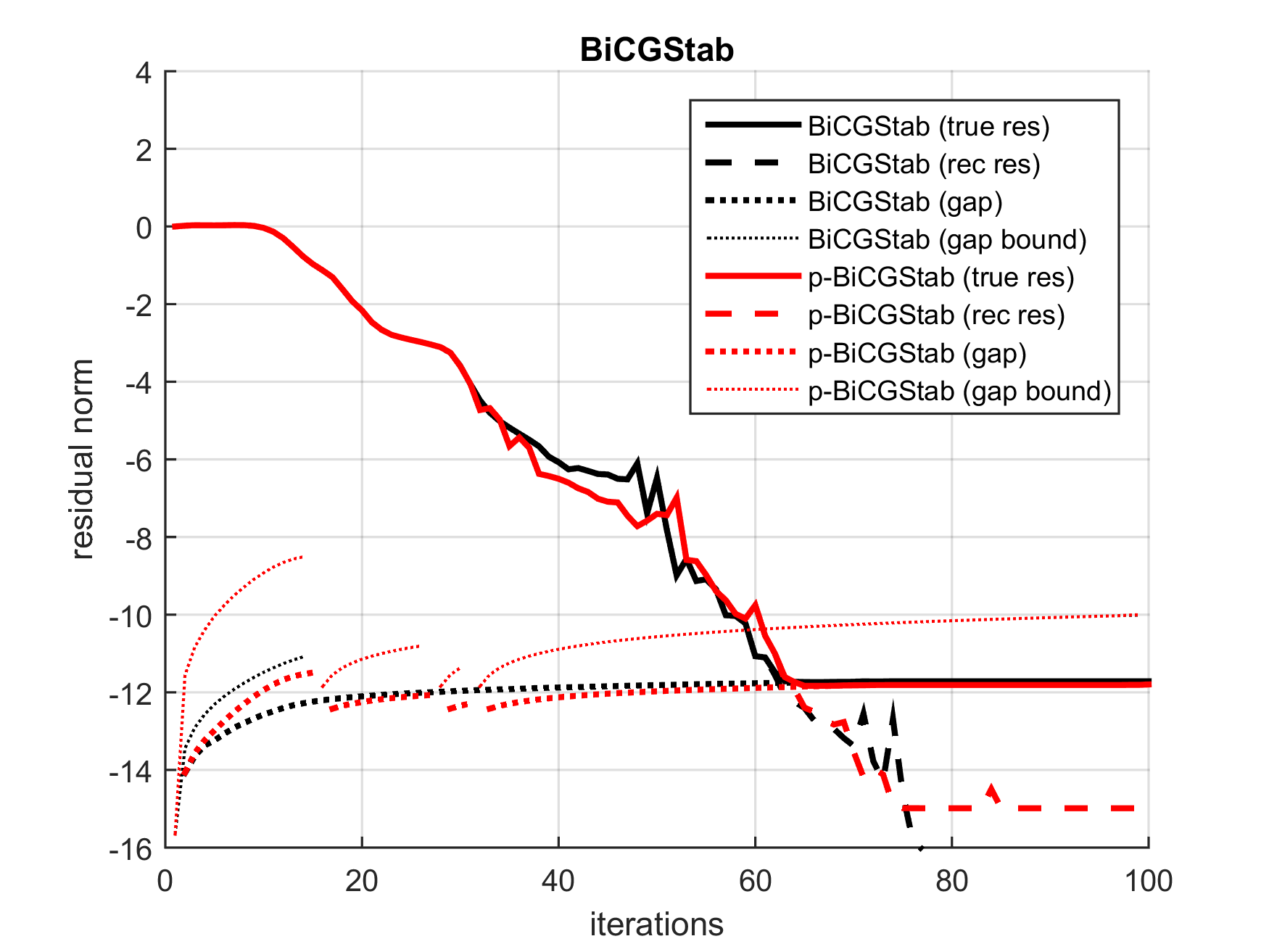} \\
\includegraphics[width=0.49\textwidth]{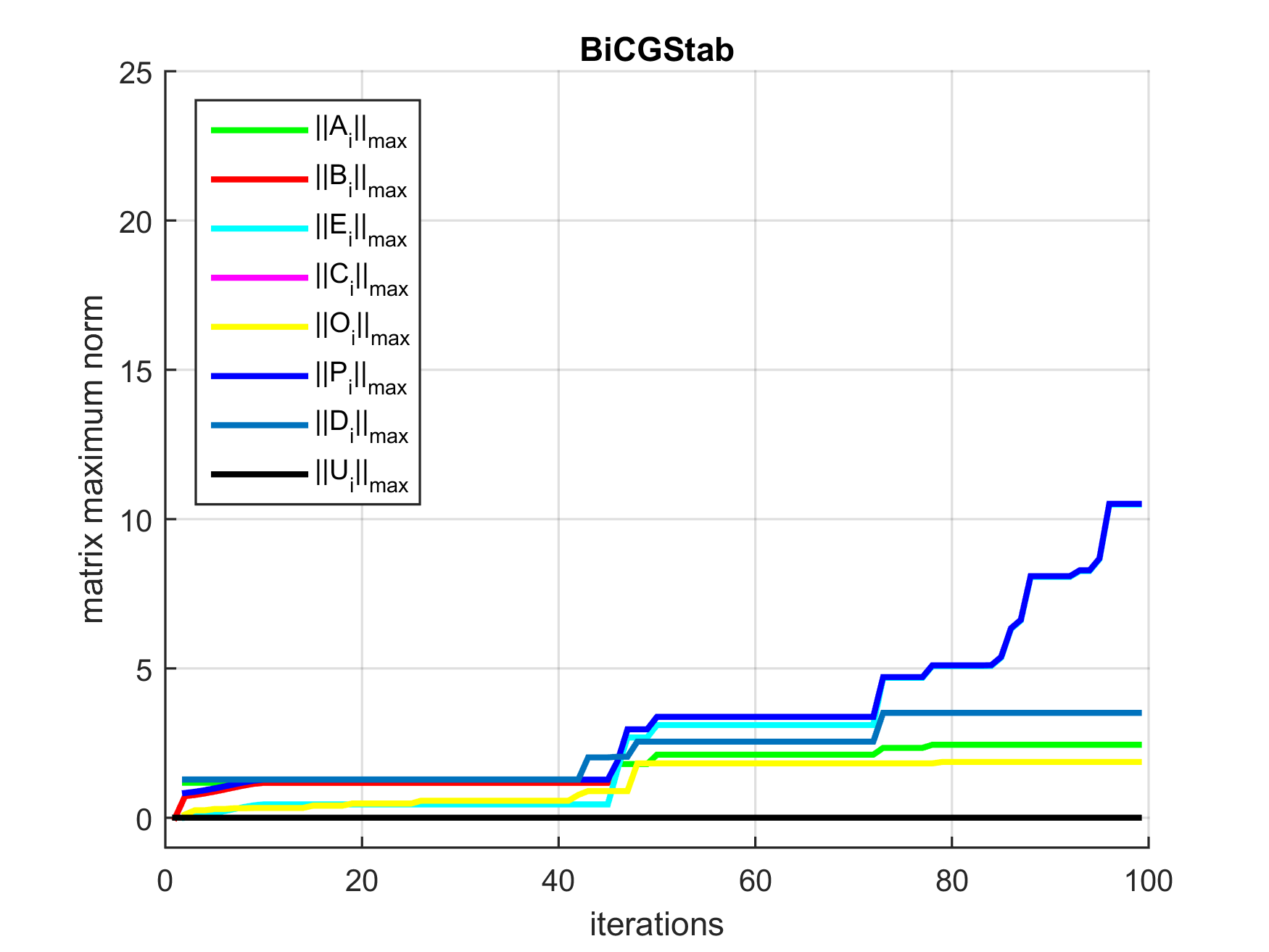}
\includegraphics[width=0.49\textwidth]{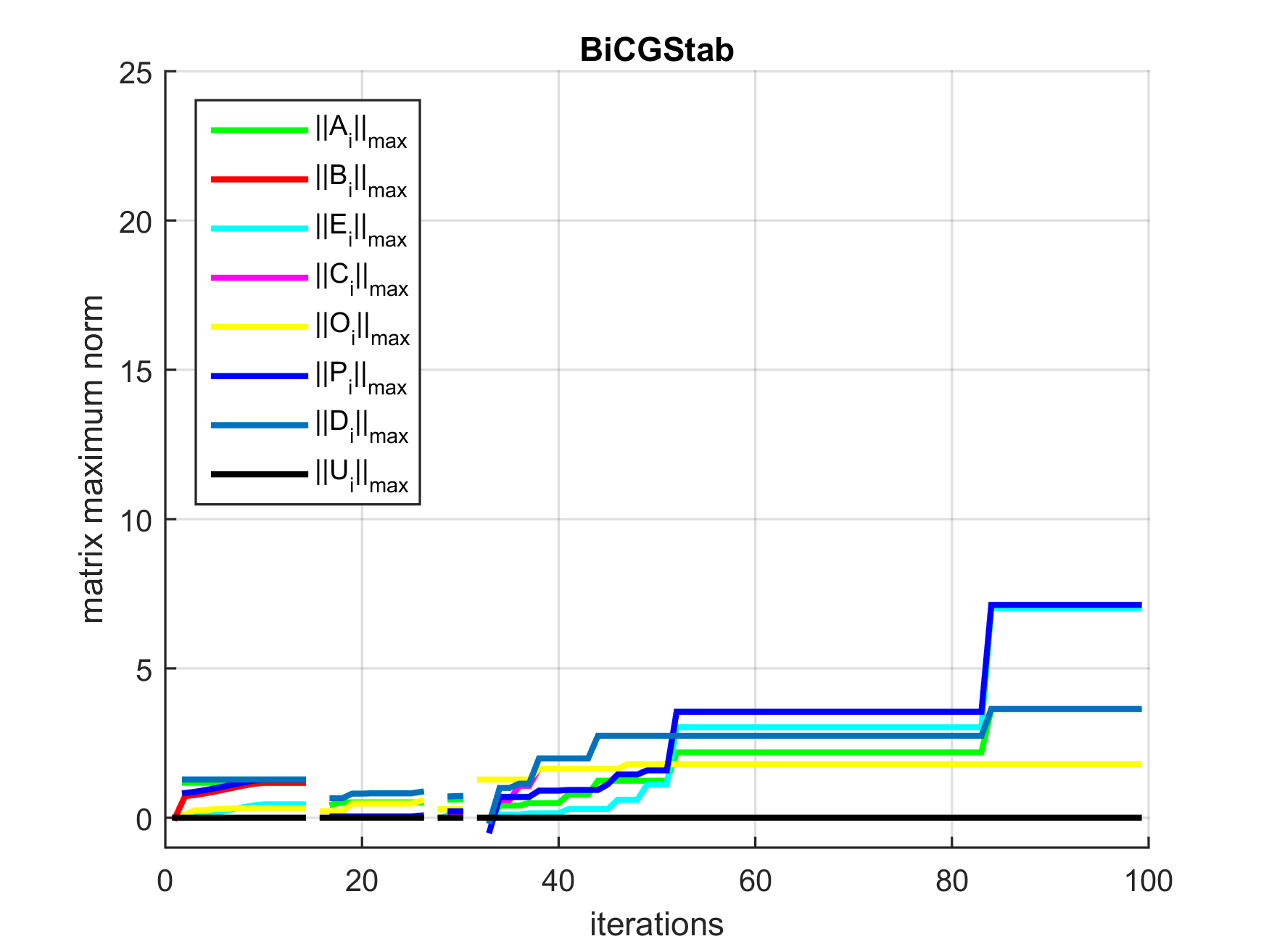} \\
\includegraphics[width=0.49\textwidth]{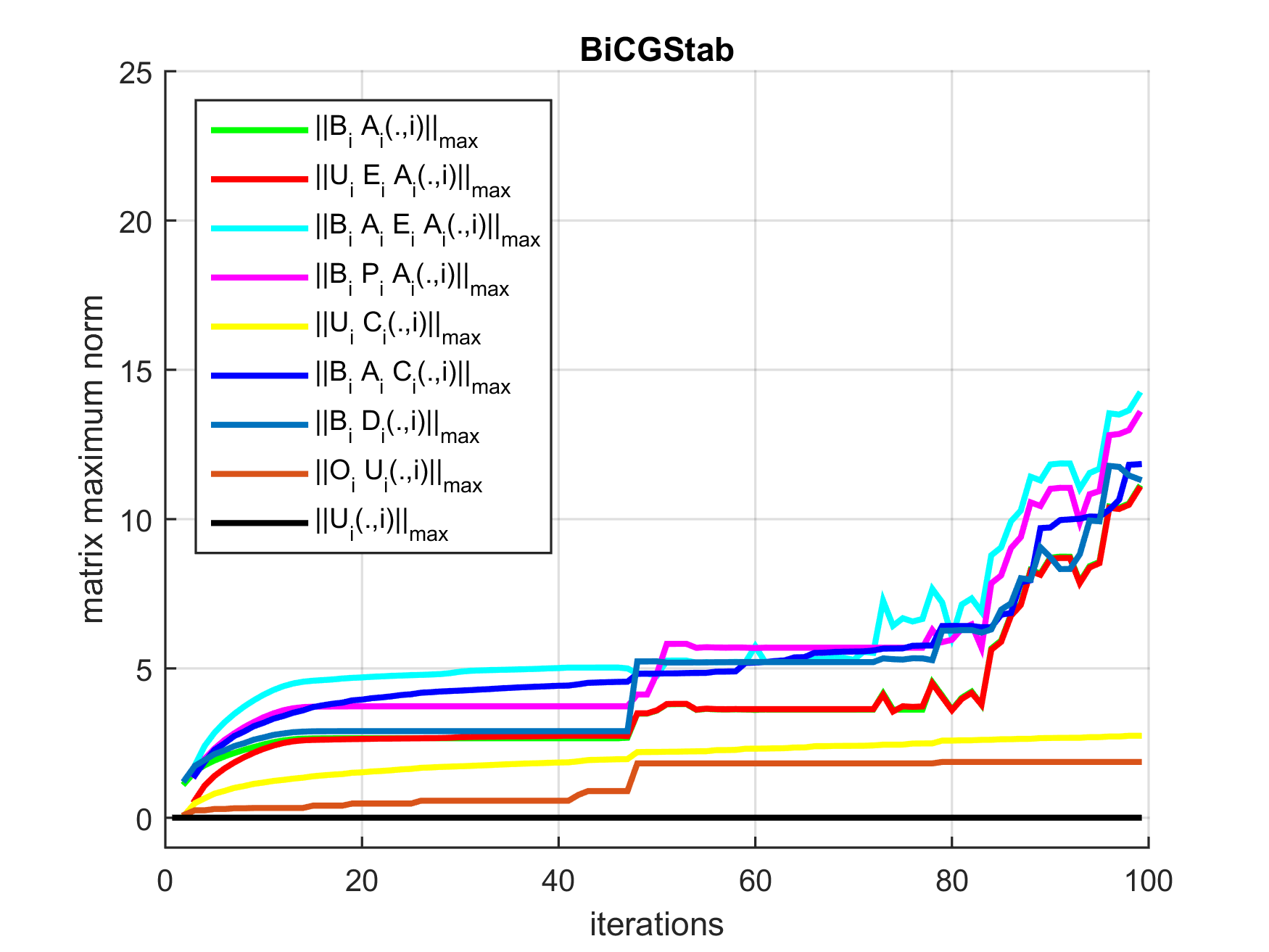}
\includegraphics[width=0.49\textwidth]{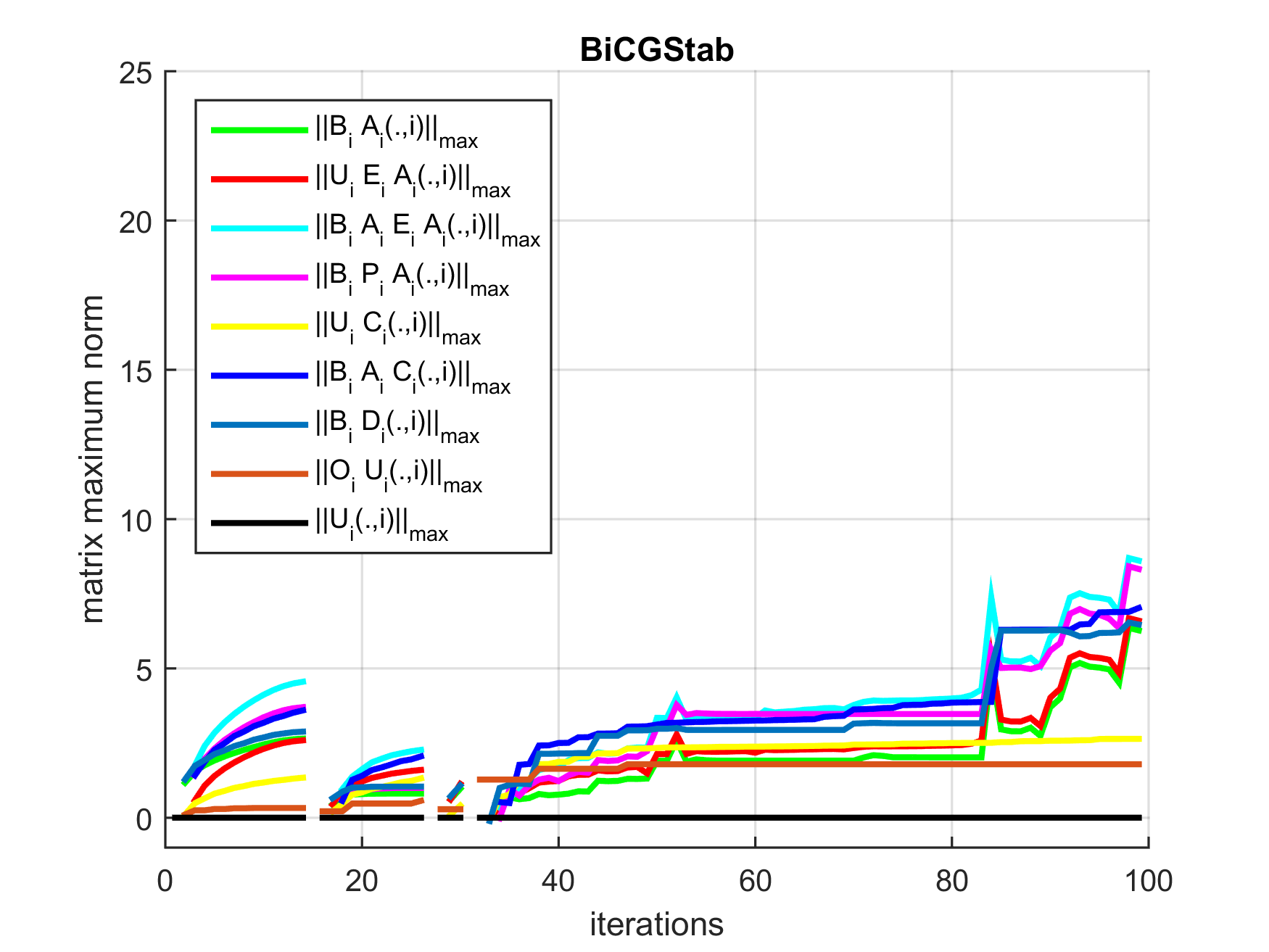} 
\end{center}
\caption{\textbf{(TP4)} Comparison BiCGStab (left) vs.~BiCGStab with fully automated residual replacements (right). \textbf{Top:} residual norm history $\|r_i\|_2$ for ICC(0) preconditioned BiCGStab/p-BiCGStab as a function of iterations. Dotted lines denote the residual gaps $\|(b-A\bar{x}_i) - \bar{r}_i\|_2$ and their computed upper bounds. \textbf{Middle:} maximum norms of various matrices occurring in the numerical stability analysis for p-BiCGStab, see \eqref{eq:matrix_expr}, as a function of iterations. \textbf{Bottom:} maximum norms of the $i$-th column of products of matrices occurring in the stability analysis \eqref{eq:all_error_matrices}. Vertical axis in $\log_{10}$ scale.}
\label{fig:9point}
\end{figure}

\textbf{Test problem 4:  Symmetric positive definite 9-point 2D stencil.} [Results shown in Fig.~\ref{fig:9point}.] Similarly to (TP1) this benchmark problem is a (normalized) FDD discretized Poisson operator of size $n_x = n_y = 200$ in each spatial dimension, such that the matrix size is $N = 40,000$. 
The nine-point matrix stencil for (TP4) is given in Table \ref{tab:specs}. The right-hand side and initial guess are defined in analogy to (TP1). 

\vspace{0.2cm}

\begin{figure}[t]
\begin{center}
\includegraphics[width=0.49\textwidth]{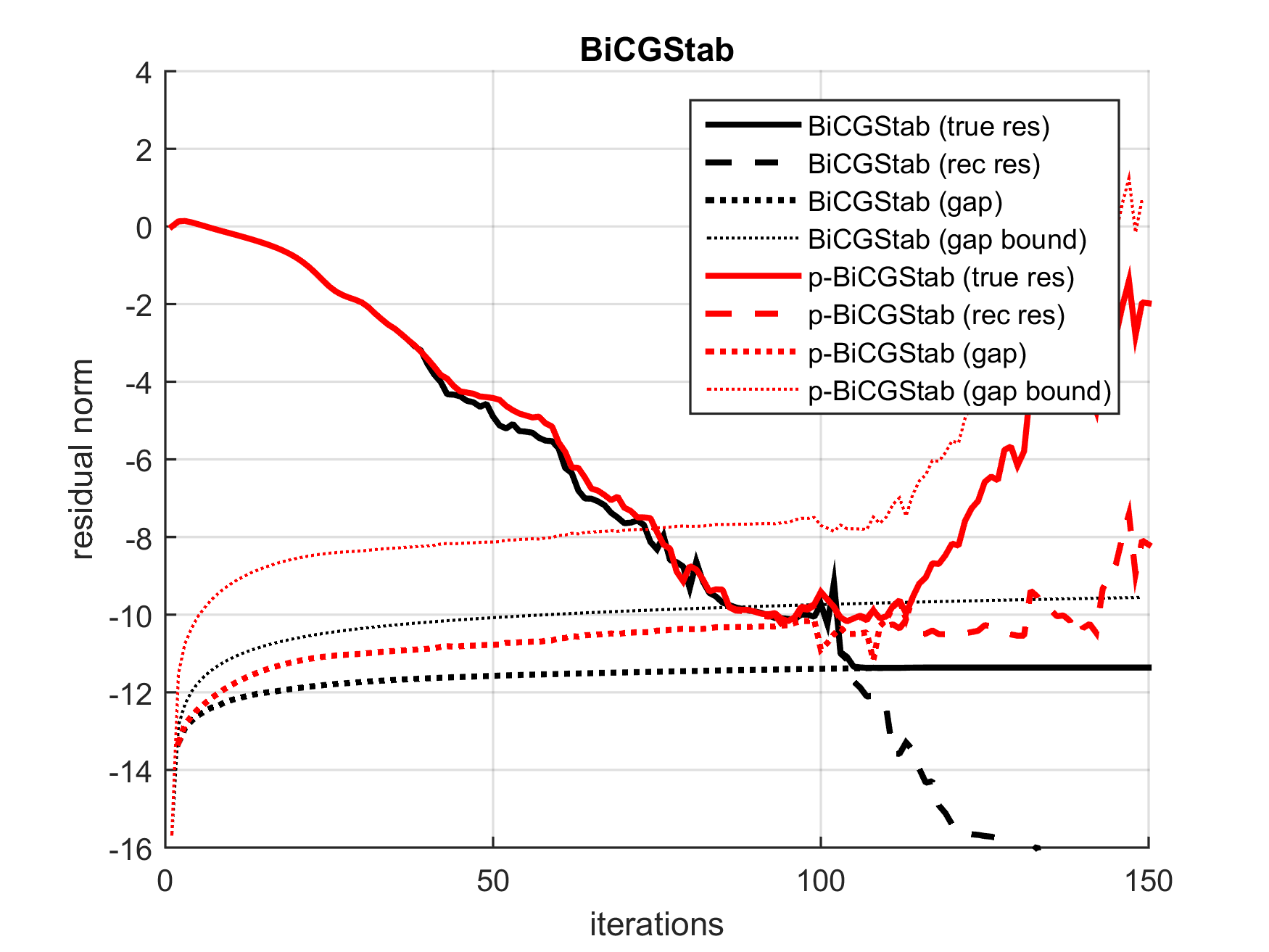}
\includegraphics[width=0.49\textwidth]{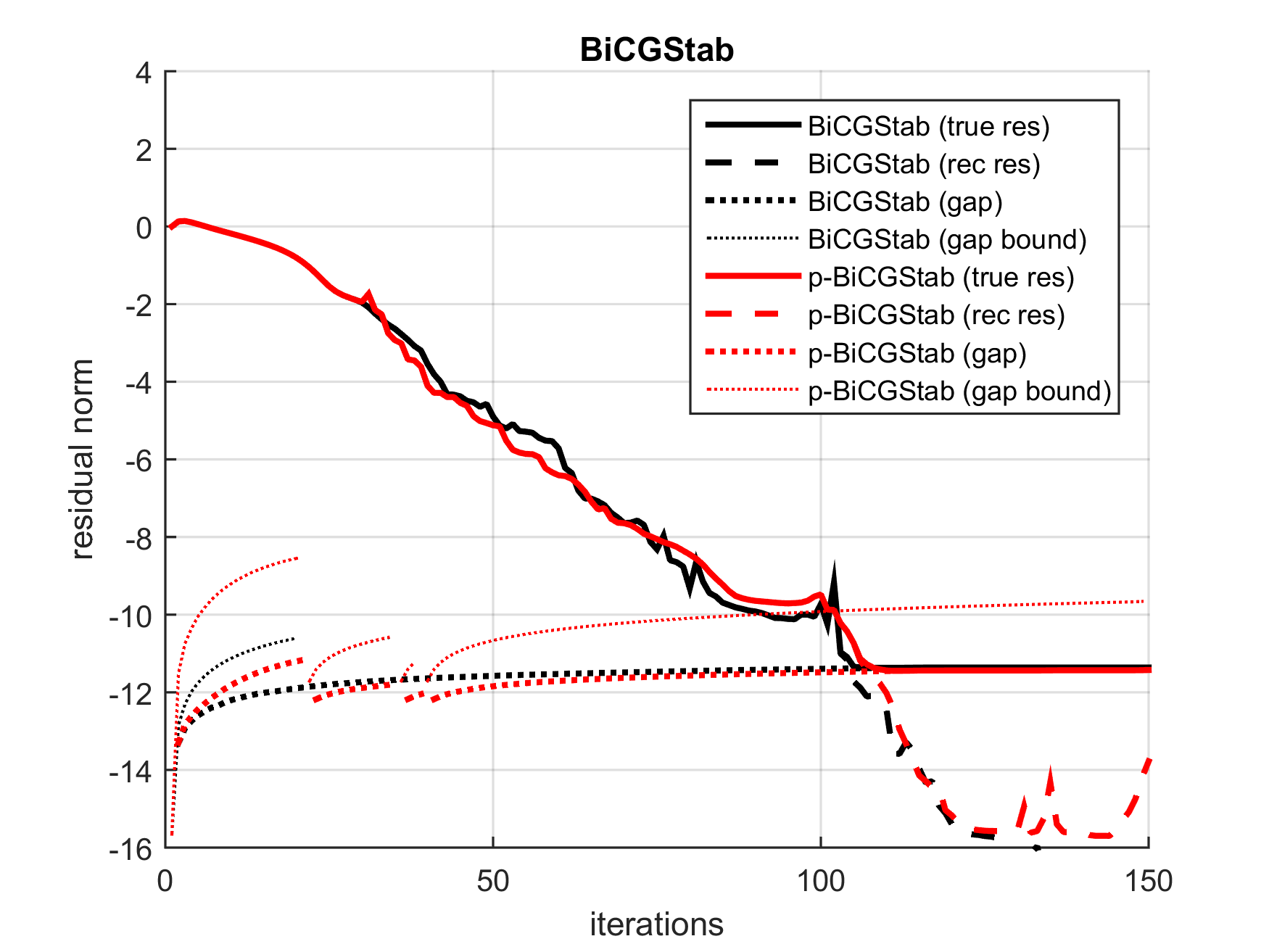} \\
\includegraphics[width=0.49\textwidth]{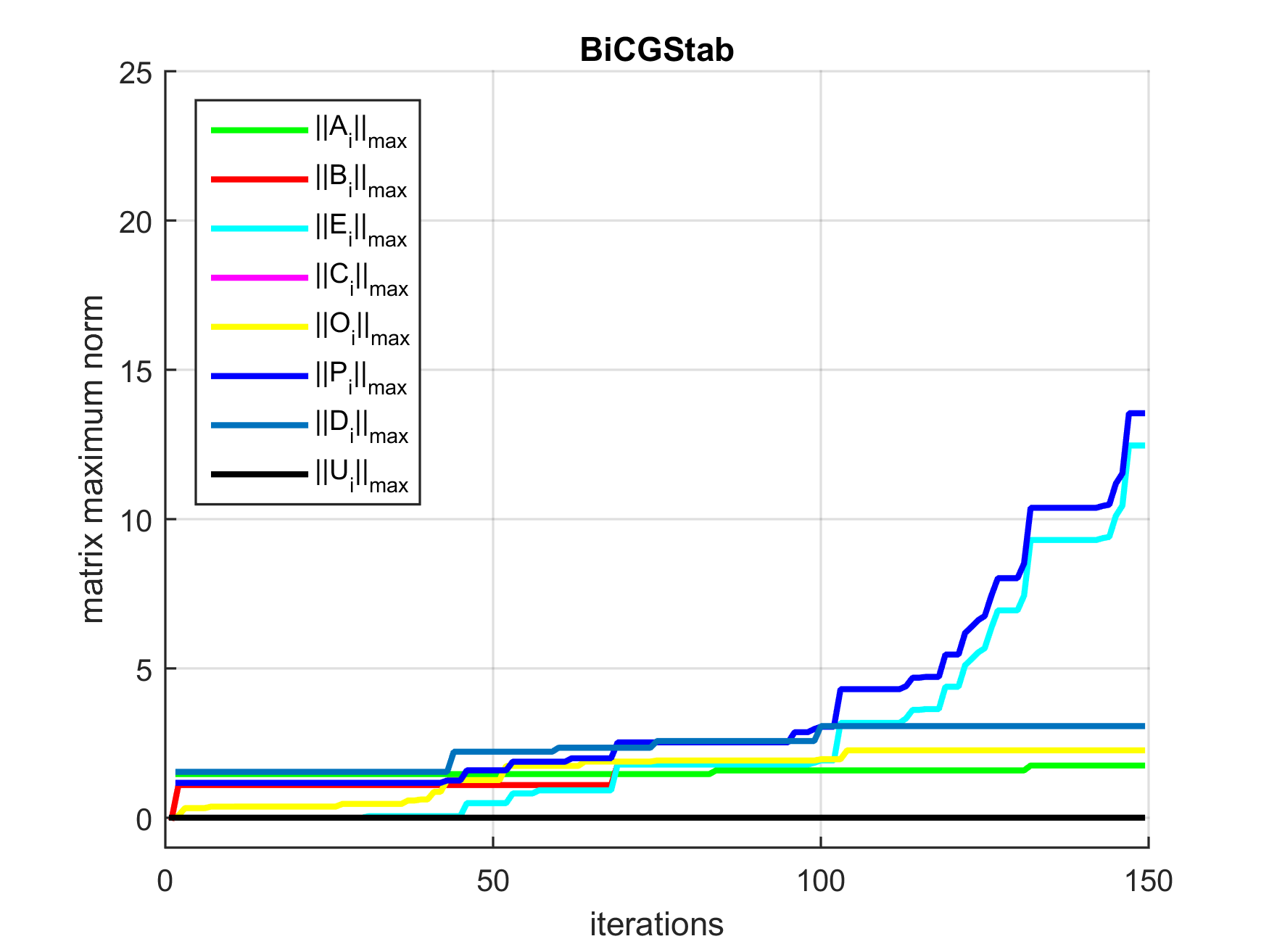}
\includegraphics[width=0.49\textwidth]{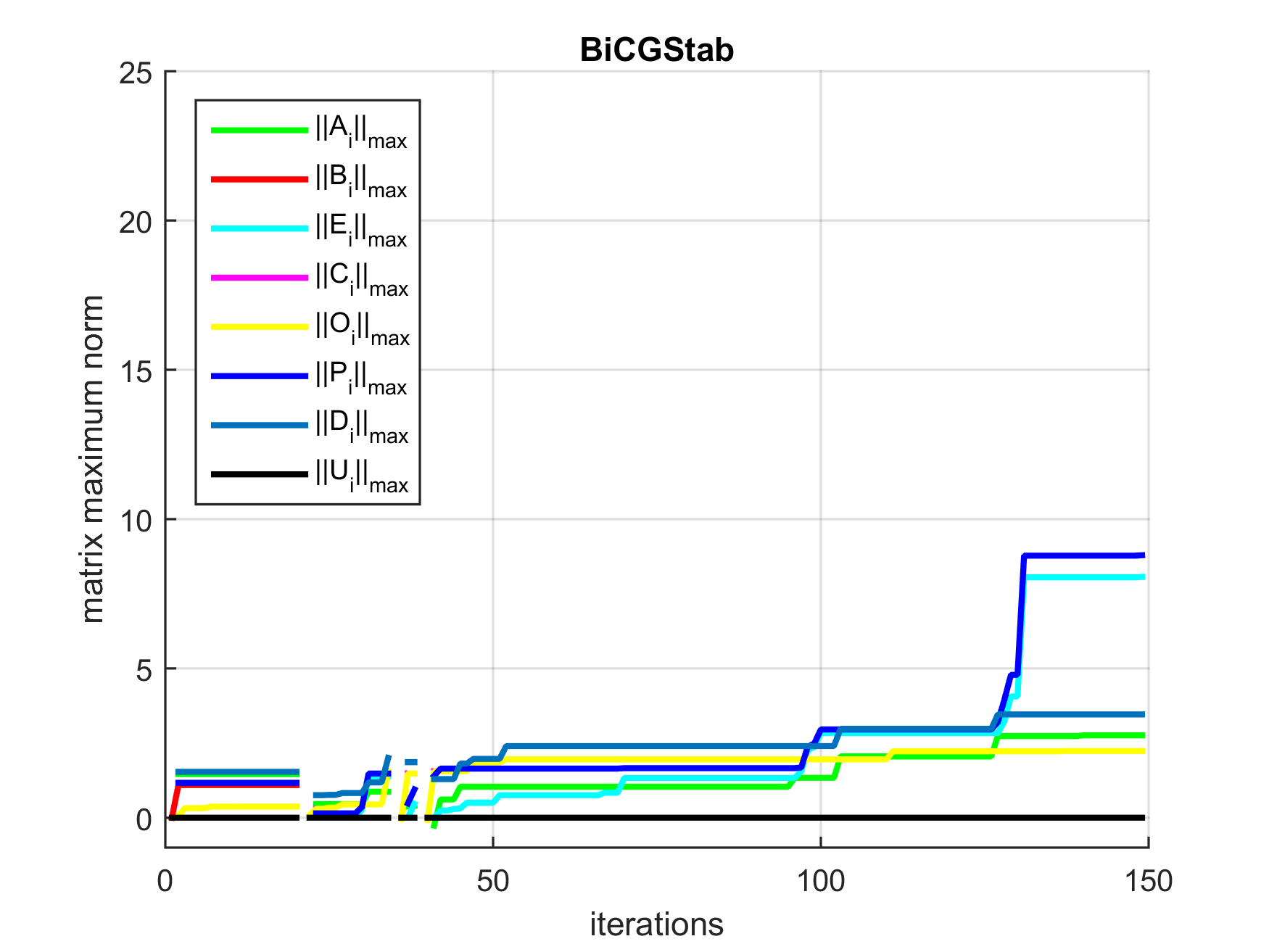} \\
\includegraphics[width=0.49\textwidth]{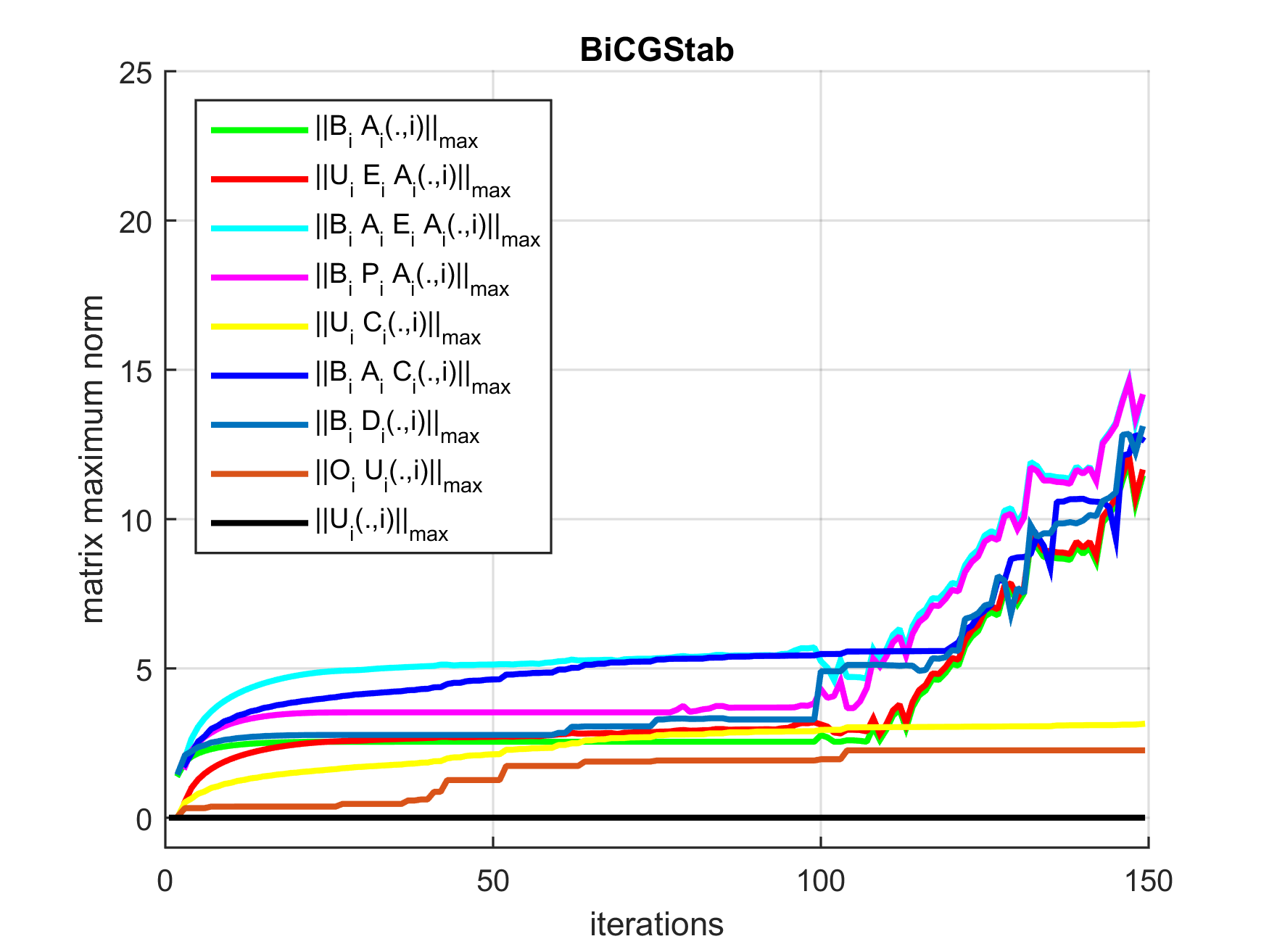}
\includegraphics[width=0.49\textwidth]{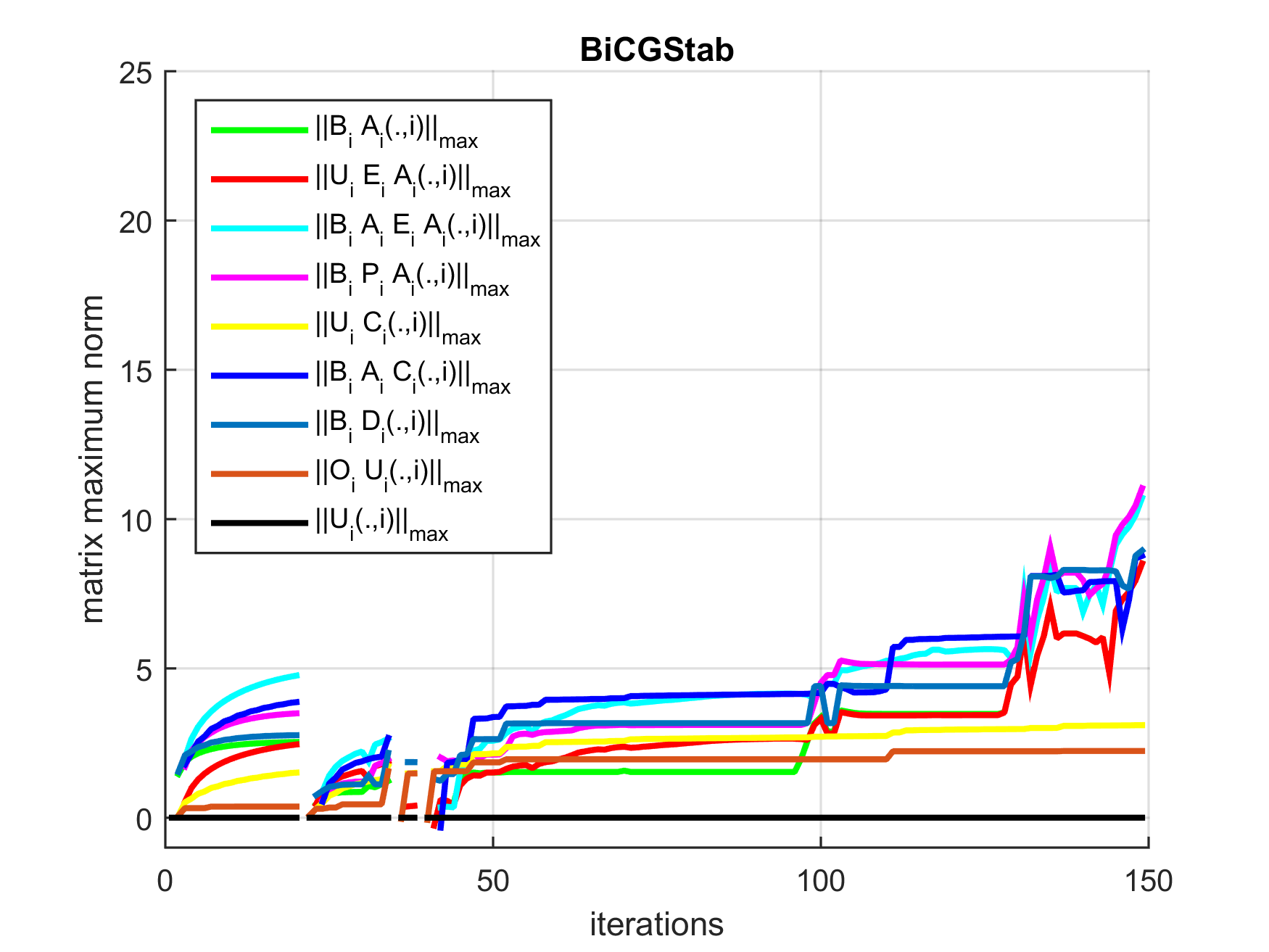} 
\end{center}
\caption{\textbf{(TP5)} Comparison BiCGStab (left) vs.~BiCGStab with fully automated residual replacements (right). \textbf{Top:} residual norm history $\|r_i\|_2$ for ICC(0) preconditioned BiCGStab/p-BiCGStab as a function of iterations. Dotted lines denote the residual gaps $\|(b-A\bar{x}_i) - \bar{r}_i\|_2$ and their computed upper bounds. \textbf{Middle:} maximum norms of various matrices occurring in the numerical stability analysis for p-BiCGStab, see \eqref{eq:matrix_expr}, as a function of iterations. \textbf{Bottom:} maximum norms of the $i$-th column of products of matrices occurring in the stability analysis \eqref{eq:all_error_matrices}. Vertical axis in $\log_{10}$ scale.}
\label{fig:7point}
\end{figure}

\textbf{Test problem 5:  Symmetric indefinite 7-point 3D stencil.} [Results shown in Fig.~\ref{fig:7point}.] The final benchmark problem is a (normalized) FDD discretized 3D Helmholtz operator of size $n_x = n_y = n_z = 50$ in each spatial dimension, yielding a matrix of size $N = 125,000$. 
Table \ref{tab:specs} shows the corresponding seven-point stencil for (TP5). The right-hand side and initial guess are identical to (TP1).

\subsection{Hardware and software specifications for parallel performance experiments} \label{sec:hardware}

We illustrate the parallel performance of the pipelined BiCGStab method, Alg.~\ref{algo::pipebicgstab}, by performing strong scaling and accuracy experiments on a small cluster with up to 15 compute nodes, consisting of two 6-core Intel Xeon X5660 Nehalem 2.80 GHz processors each (12 cores per node). 
These nodes are connected by $4\,\times\,$QDR InfiniBand technology which provides 32 Gb/s of point-to-point bandwidth for message passing and I/O.
To fully exploit parallelism on the machine 12 MPI processes per node are used.
The MPI library used for this experiment is MPICH-3.1.3\footnote{\url{http://www.mpich.org/}}. The environment variables 
\texttt{MPICH\_ASYNC\_PROGRESS=1} and \texttt{MPICH\_MAX\_THREAD\_SAFETY=multiple} are set to ensure optimal 
parallelism
; the first variable enables asynchronous non-blocking reductions, while the second allows to have multiple threads
calling MPI functions simultaneously.

The pipelined BiCGStab Alg.~\ref{algo::pipebicgstab} is available in the open-source PETSc library \cite{petsc-web-page} since v.\,3.8.4. The method is implemented as a modified version of the (flexible) BiCGStab implementation \texttt{fbcgs} and can be found in the PETSc Krylov solvers (KSP) folder.

The inclusion of a preconditioner in pipelined algorithms is generally straightforward from the algorithmic point of view. However, to efficiently overlap the preconditioner application with the global communication phase, the preconditioner should not be bottlenecked by communication. Block preconditioners or (non-overlapping) domain decomposition methods \cite{smith2004domain}, for example, are well-suited for this purpose, whereas the inclusion of advanced preconditioning schemes like parallel ILU \cite{chow2015fine} or specialized physics-based preconditioners or multigrid techniques for the Helmholtz equation (TP3) \cite{engquist2011sweeping,erlangga2004class} require more careful and specific treatment. For simplicity no preconditioner is included in the performance experiments reported in Figs.~\ref{fig:timings2}-\ref{fig:timings3}.  

\subsection{Numerical comparison between pipelined CG and pipelined BiCGStab} \label{sec:comparison}

\begin{figure}[t]
\begin{center}
\includegraphics[width=0.49\textwidth]{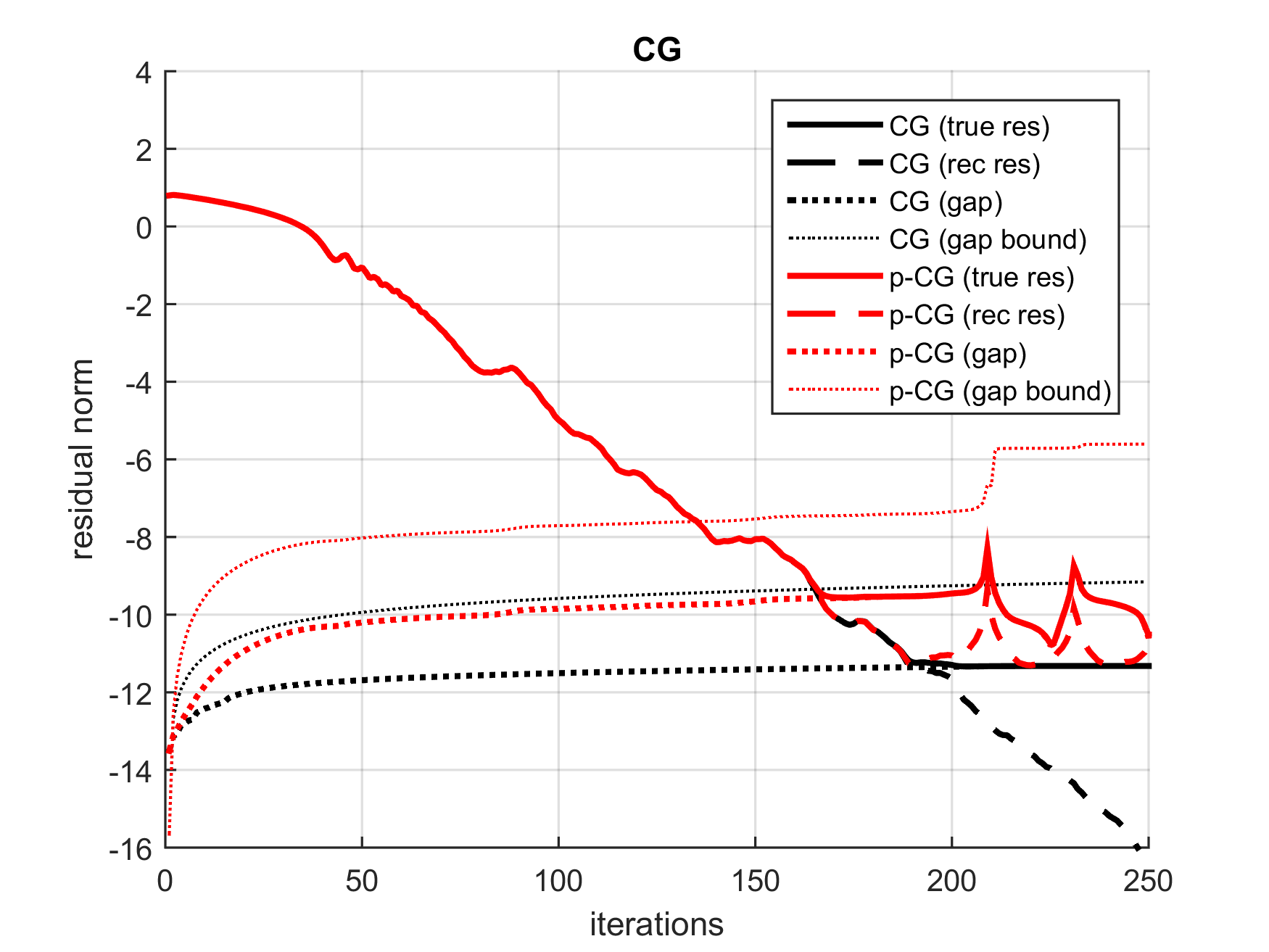}
\includegraphics[width=0.49\textwidth]{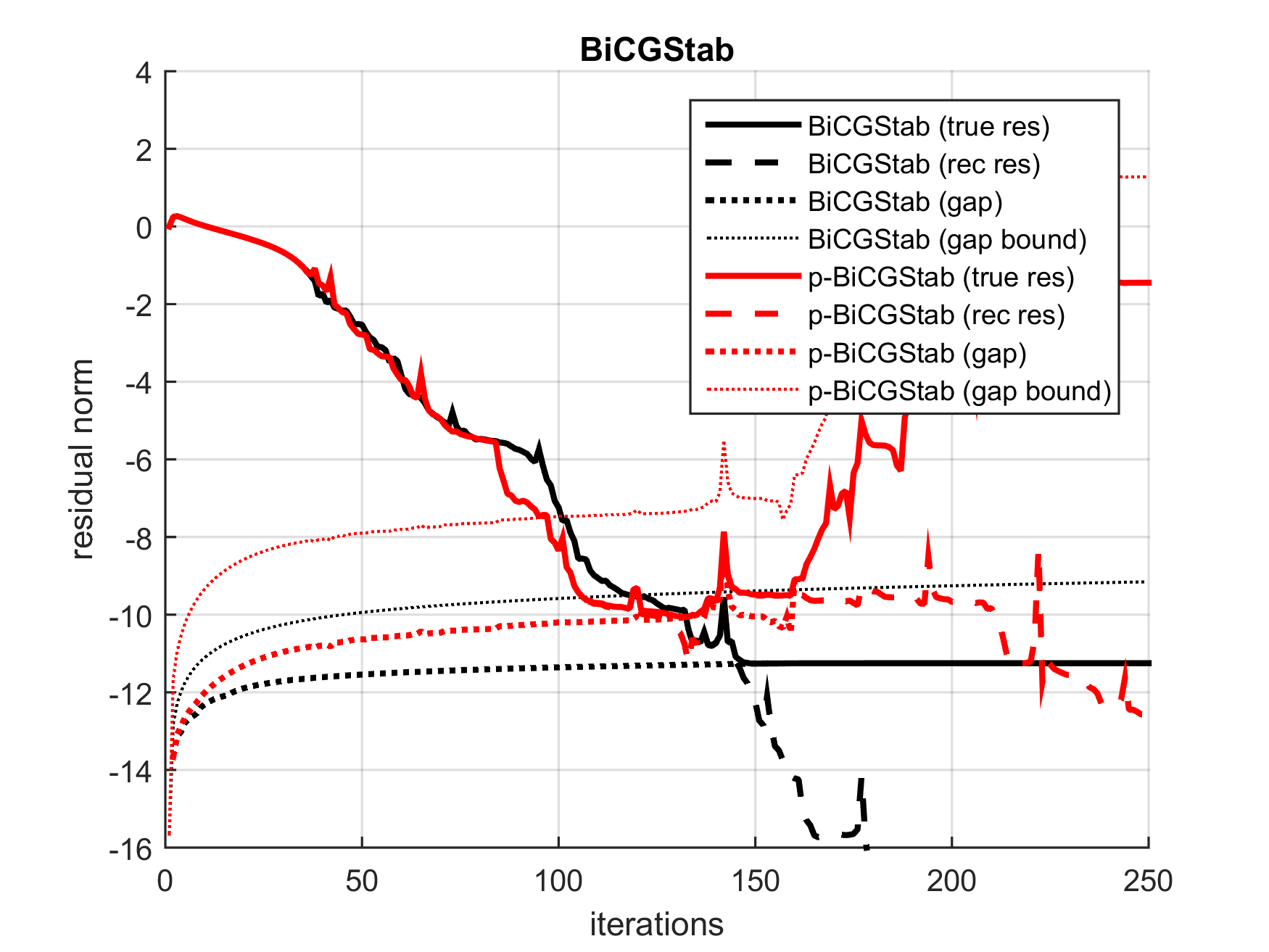} \\
\includegraphics[width=0.49\textwidth]{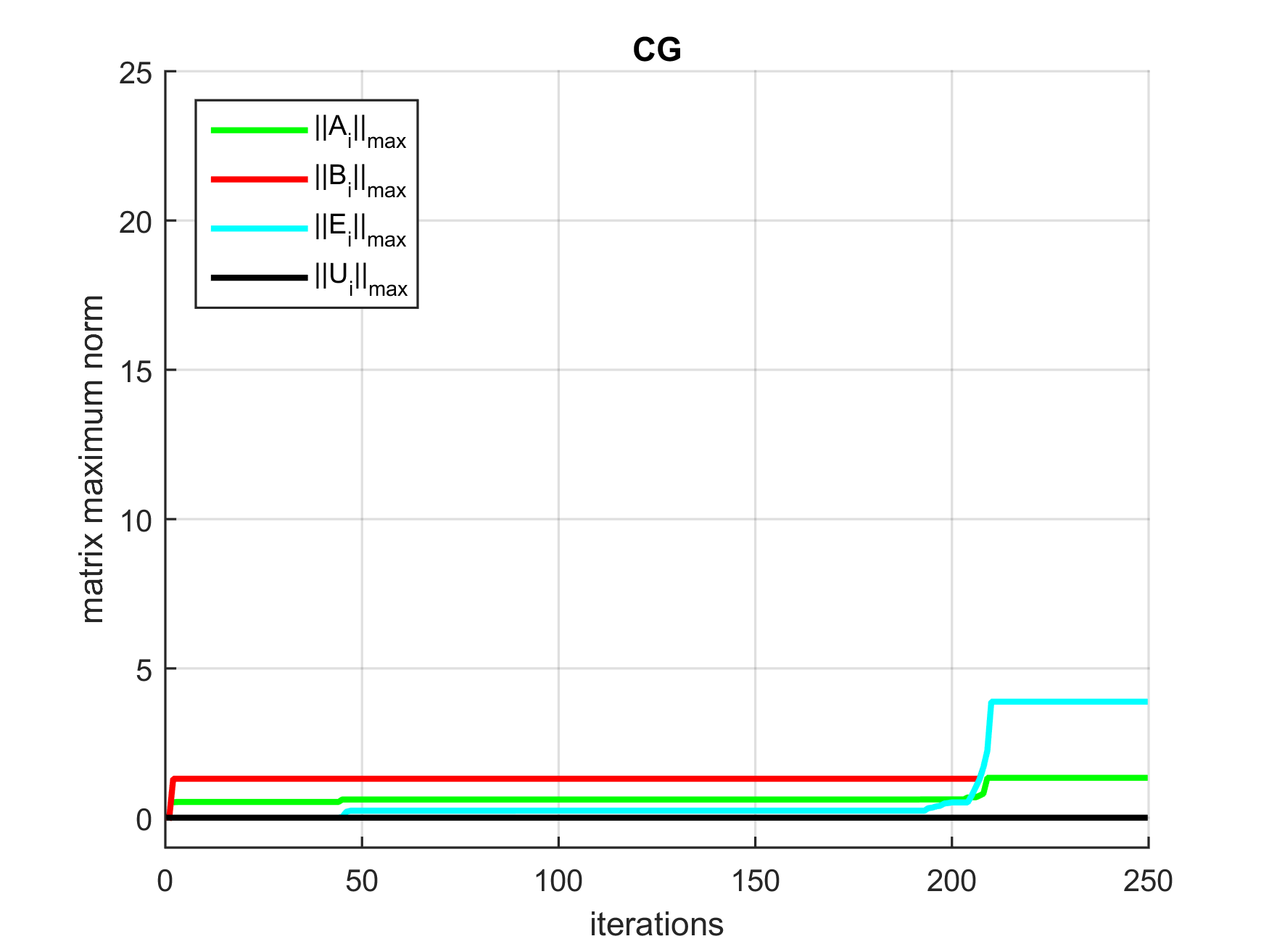}
\includegraphics[width=0.49\textwidth]{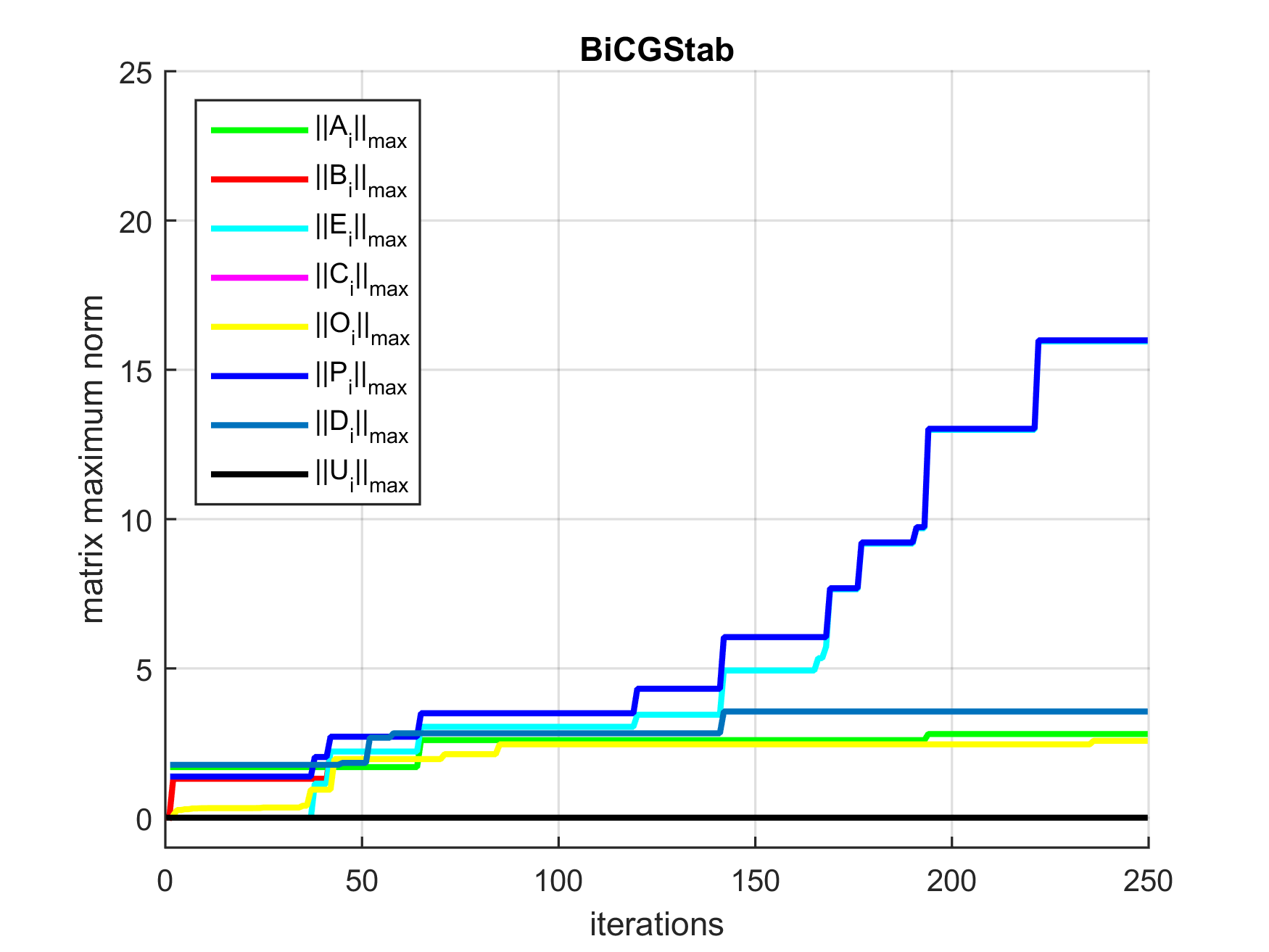} \\
\includegraphics[width=0.49\textwidth]{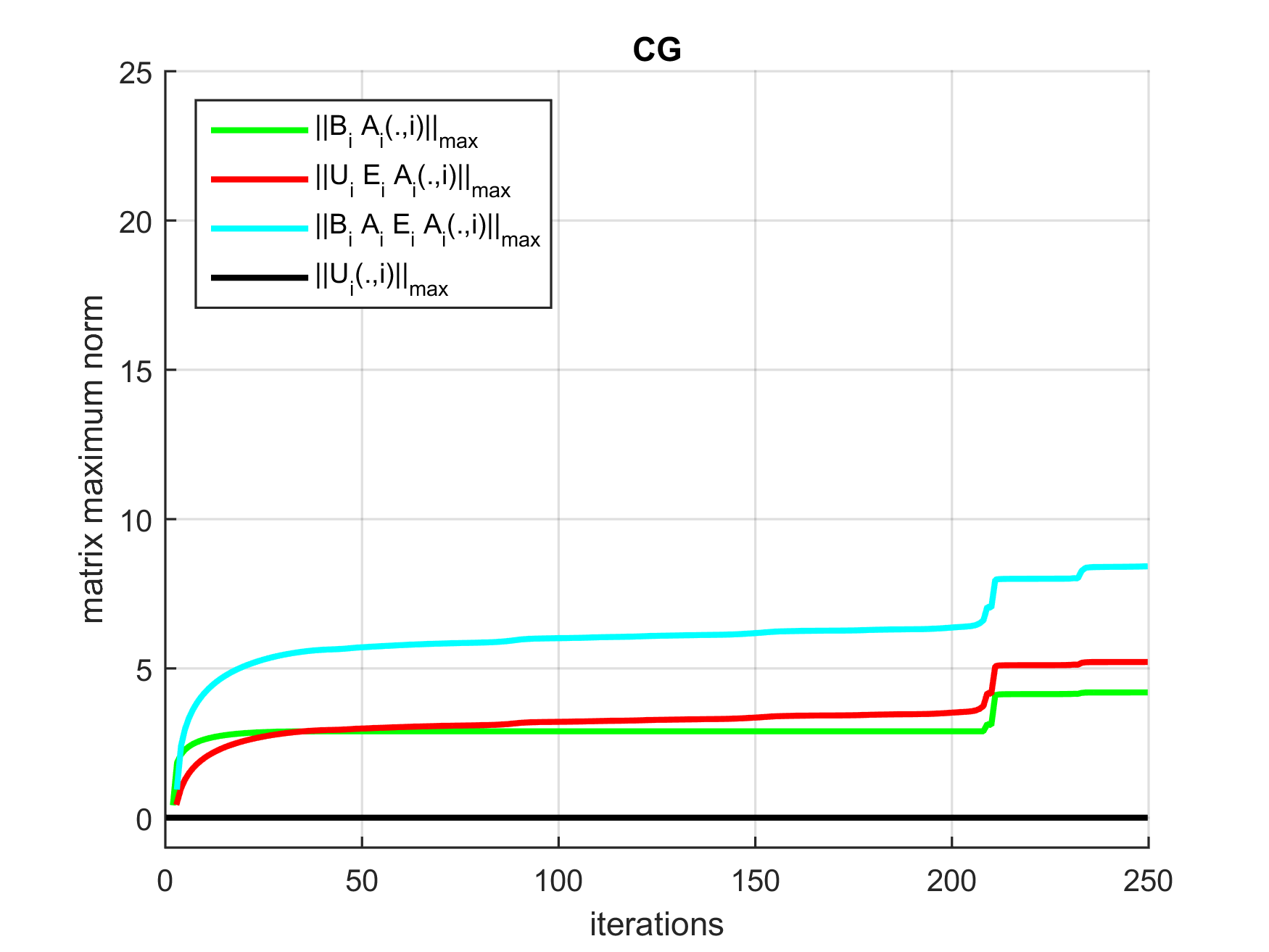}
\includegraphics[width=0.49\textwidth]{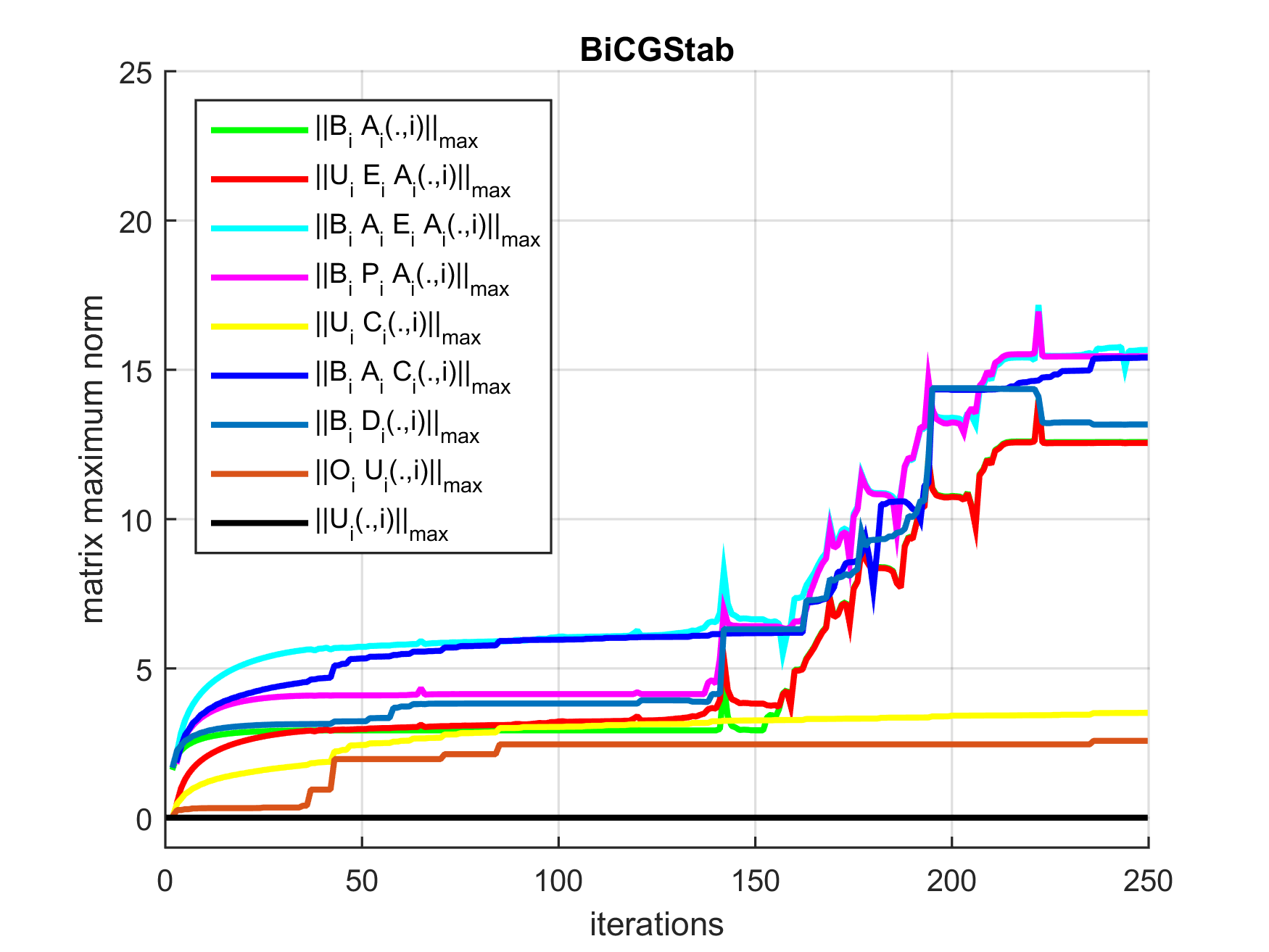} 
\end{center}
\caption{\textbf{(TP1)} Comparison CG (left) vs.~BiCGStab (right). \textbf{Top:} residual norm history $\|r_i\|_2$ for ICC(0) preconditioned CG/p-CG and BiCGStab/p-BiCGStab as a function of iterations. Dotted lines denote the residual gaps $\|(b-A\bar{x}_i) - \bar{r}_i\|_2$ and their computed upper bounds. \textbf{Middle:} maximum norms of various matrices occurring in the numerical stability analysis for p-CG (left) and p-BiCGStab (right), see \eqref{eq:matrix_expr}, as a function of iterations. \textbf{Bottom:} maximum norms of the $i$-th column of products of matrices occurring in the stability analysis \eqref{eq:all_error_matrices}. Vertical axis in $\log_{10}$ scale.}
\label{fig:comparison_ichol}
\end{figure}

Figure \ref{fig:comparison_ichol} (top) shows the convergence histories of CG/p-CG (left) and BiCGStab/p-BiCGStab (right) for solving the SPD benchmark problem (TP1) with an Incomplete Cholesky preconditioner with zero fill-in (ICC(0)). The norms of the true residuals $b-A\bar{x}_i$ (full line), the recursively computed residuals $\bar{r}_i$ (dashed line) and the residual gaps $(b-A\bar{x}_i)-\bar{r}_i$ (dotted line) are displayed as a function of iterations. The true residuals stagnate when the norm of the residual gap becomes larger than the actual residual norm. The recursively computed residuals start to deviate from the true residuals at this point and keep decreasing although the actual precision on the solution no longer increases.
We point out that one SpMV is computed in each iteration of the CG algorithms, whereas in BiCGStab (Alg.~\ref{algo::bicgstab}-\ref{algo::pipebicgstab}) a residual is computed for each two SpMVs. BiCGStab thus generally converges slower than CG when compared on the same symmetric model problem, although the `number of iterations' reported on the horizontal axes in the figures may appear to be smaller for BiCGStab than for CG.

Numerically unstable behavior can be observed from Fig.~\ref{fig:comparison_ichol} (top) for the p-BiCGStab method beyond the stagnation point. These numerical instabilities were also reported in the original publication on pipelined BiCGStab \cite{cools2017communication}. For (TP1) they are observed after iteration 140 at which point the true residual norms suddenly increase. The p-CG method also shows unstable behavior (from iteration 200 onward) although the effect on convergence is generally less pronounced. 
Figure \ref{fig:comparison_ichol} (top) furthermore presents the upper bound $f^r_i$ for the residual gap (thin dotted line) based on the analysis in Section \ref{sec:expressing}, cf.~expression \eqref{eq:DDr_bound}. 
As shown in the figure the resulting bound is not guaranteed to be sharp; 
however, the bound captures the general behavior of the gap norm quite well. 

Figure \ref{fig:comparison_ichol} (middle) shows the maximum norms of the matrices that occur in the stability analysis, cf.\,\eqref{eq:matrix_expr}. This figure indicates that p-BiCGStab is significantly more prone to rounding error amplification compared to p-CG, leading to the unstable behavior observed in Fig.~\ref{fig:comparison_ichol} (top).

Figure \ref{fig:comparison_ichol} (bottom) presents the maximum norms of the $i$-th column of the products of matrices 
from \eqref{eq:all_error_matrices}. These columns are multiplied by the respective local rounding errors from expression \eqref{eq:matrix_expr} to form the gap $\Delta^r_i$ on the residual $\bar{r}_i$ in iteration $i$. The matrix and vector norms shown in the middle and bottom panels do not account for a full quantitative characterization of the respective gaps, since the numerical value of these gaps also depends on the size of the local rounding error contributions. However, when the norm of the error propagation matrices grows large, local rounding errors can be expected to be amplified, resulting in a possibly dramatic increase of the residual gap. On the other hand, when the matrix norms do not increase rounding errors do not propagate and the residual gap is not expected to increase dramatically, although a moderate increase due to local rounding error accumulation may still be observed.  


In summary, it is clear from Figure \ref{fig:comparison_ichol} that the attainable accuracy for p-CG and p-BiCGStab on (TP1) is comparable, but the minimal attainable residual norm for the pipelined methods lies several orders of magnitude above the accuracy attainable by the standard Krylov subspace algorithms. This observation is substantiated by the norms shown in the middle and bottom panels of Fig.~\ref{fig:comparison_ichol}.

\subsection{Pipelined BiCGStab with periodic residual replacements} \label{sec:residual}

\begin{figure}[t]
\begin{center}
\includegraphics[width=0.49\textwidth]{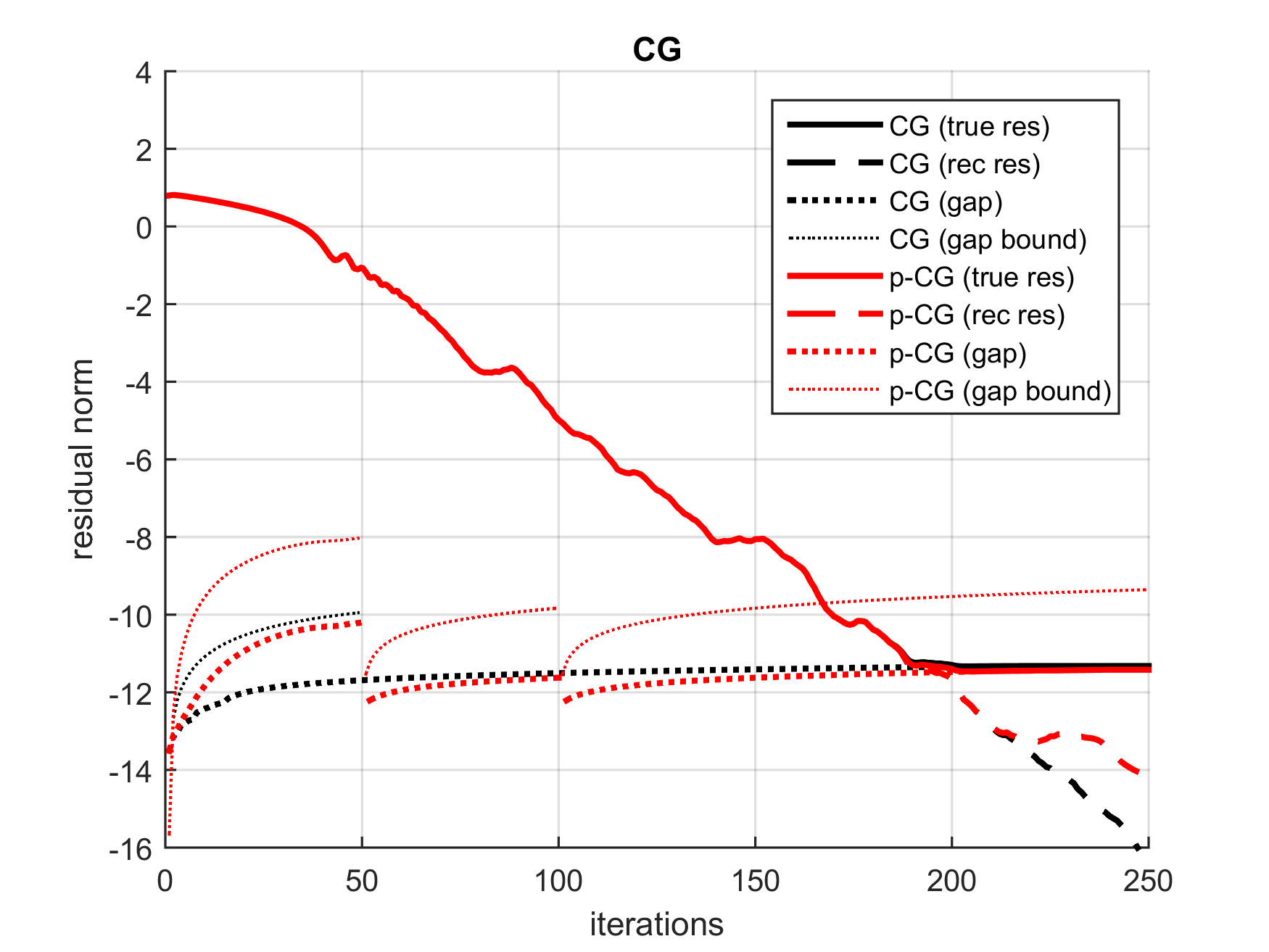}
\includegraphics[width=0.49\textwidth]{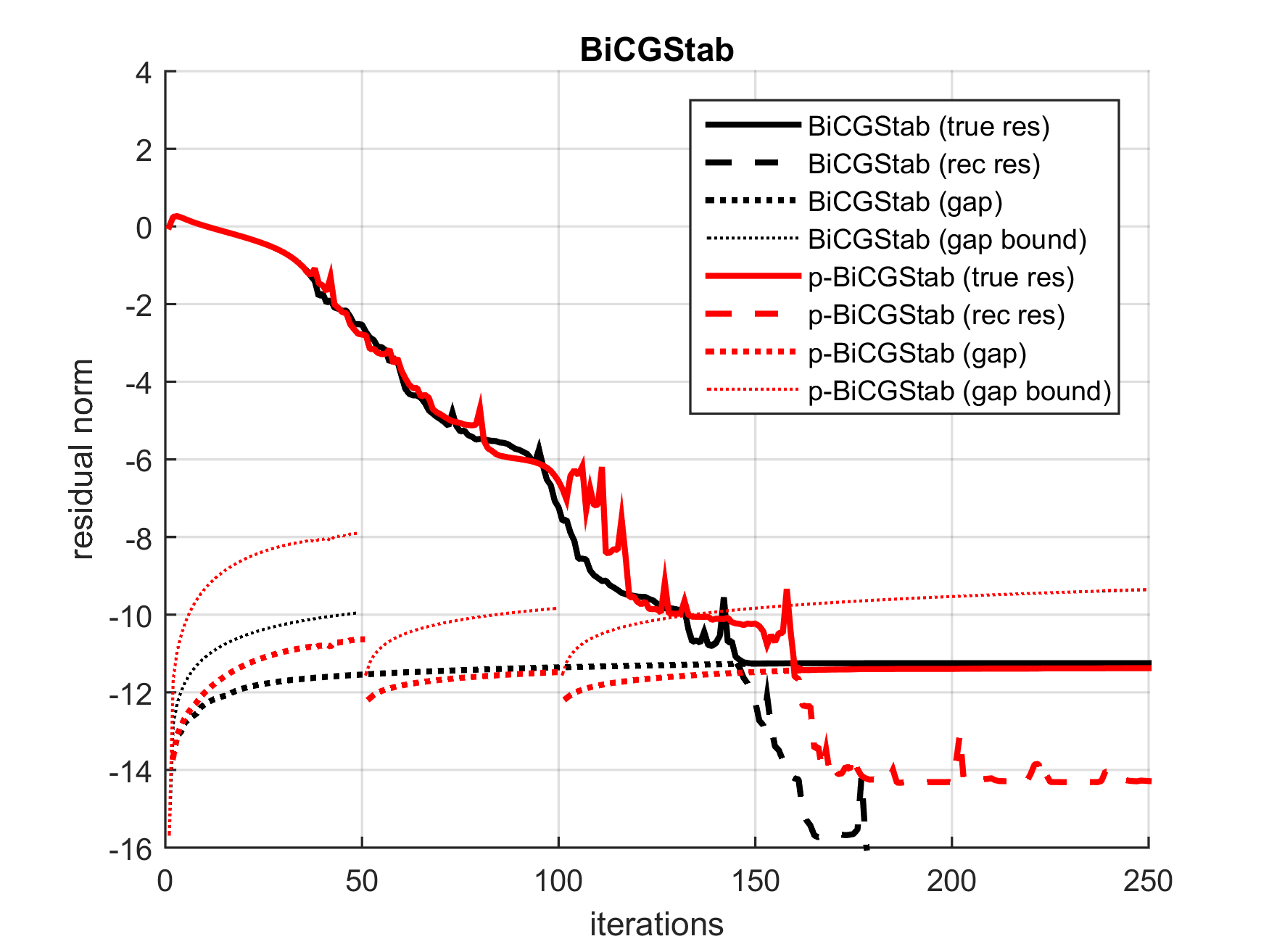} \\
\includegraphics[width=0.49\textwidth]{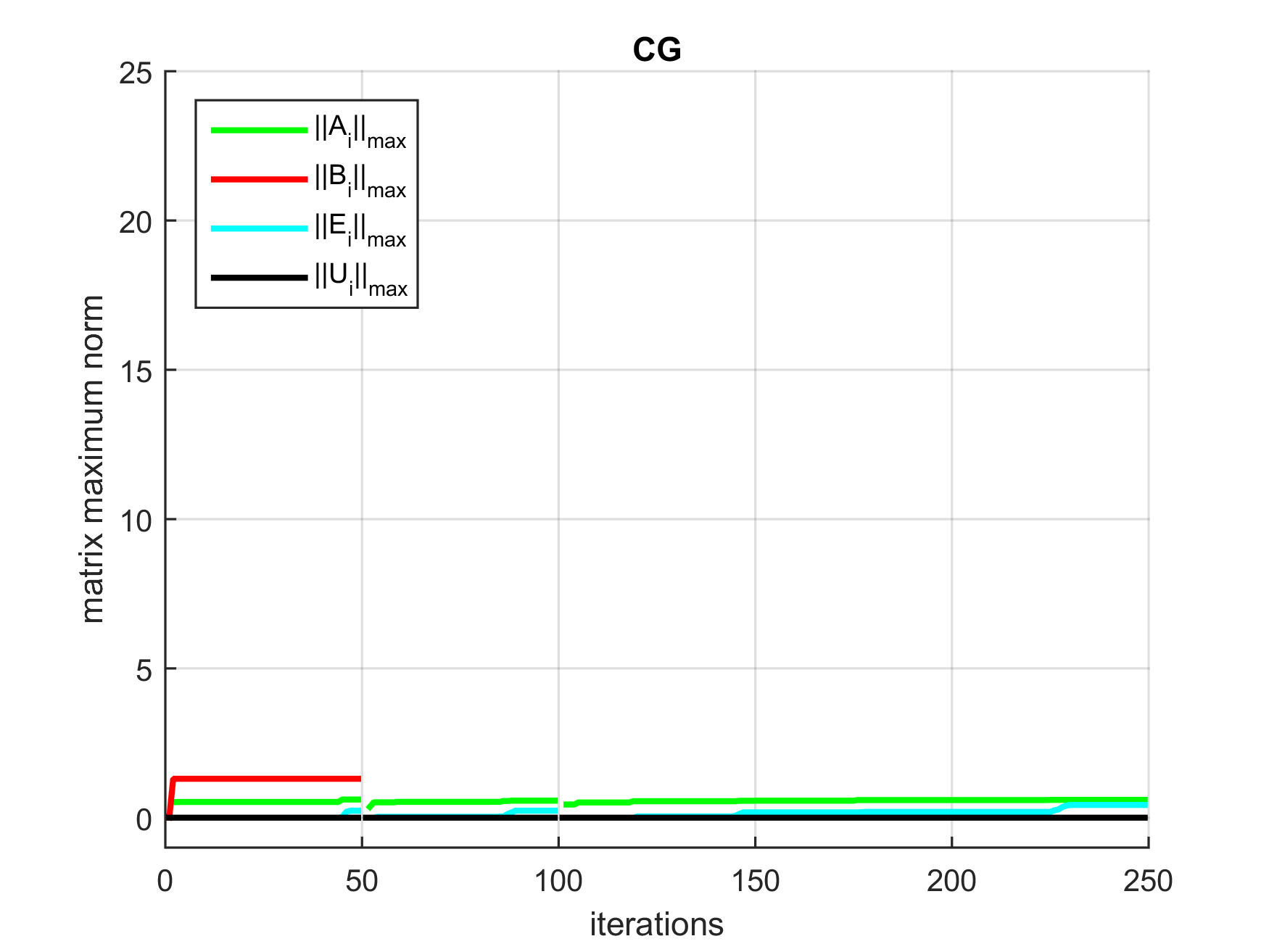}
\includegraphics[width=0.49\textwidth]{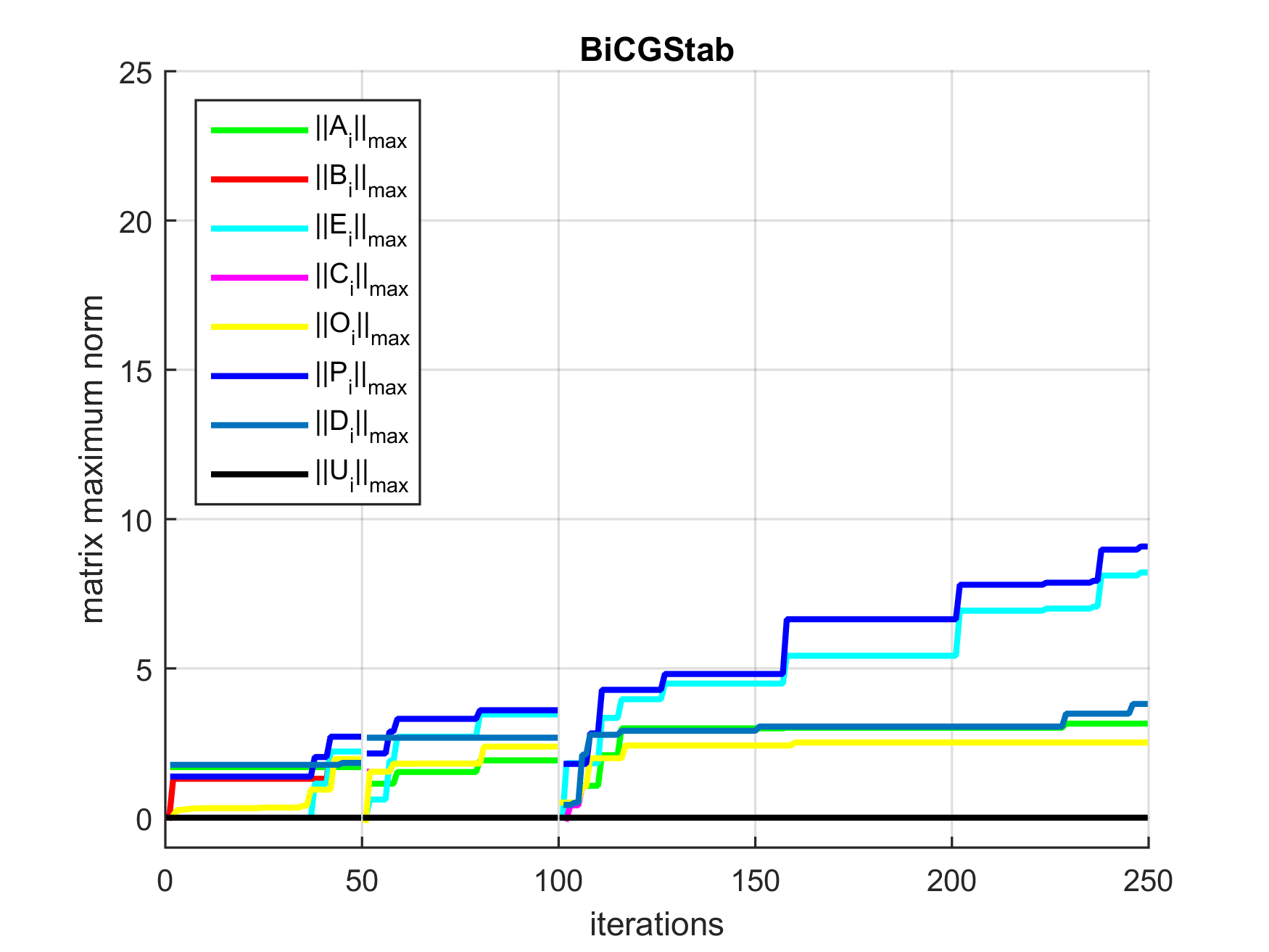} \\
\includegraphics[width=0.49\textwidth]{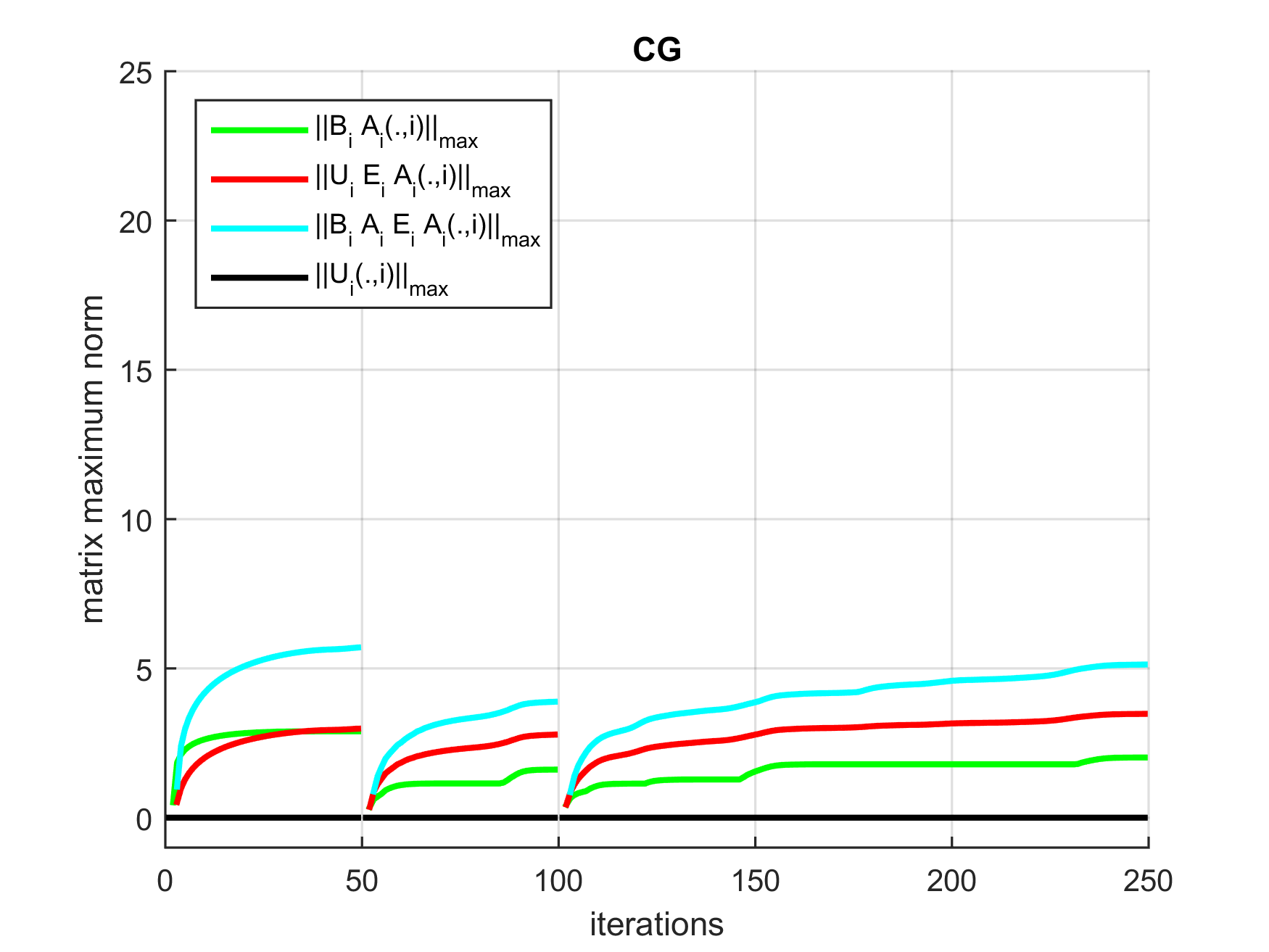}
\includegraphics[width=0.49\textwidth]{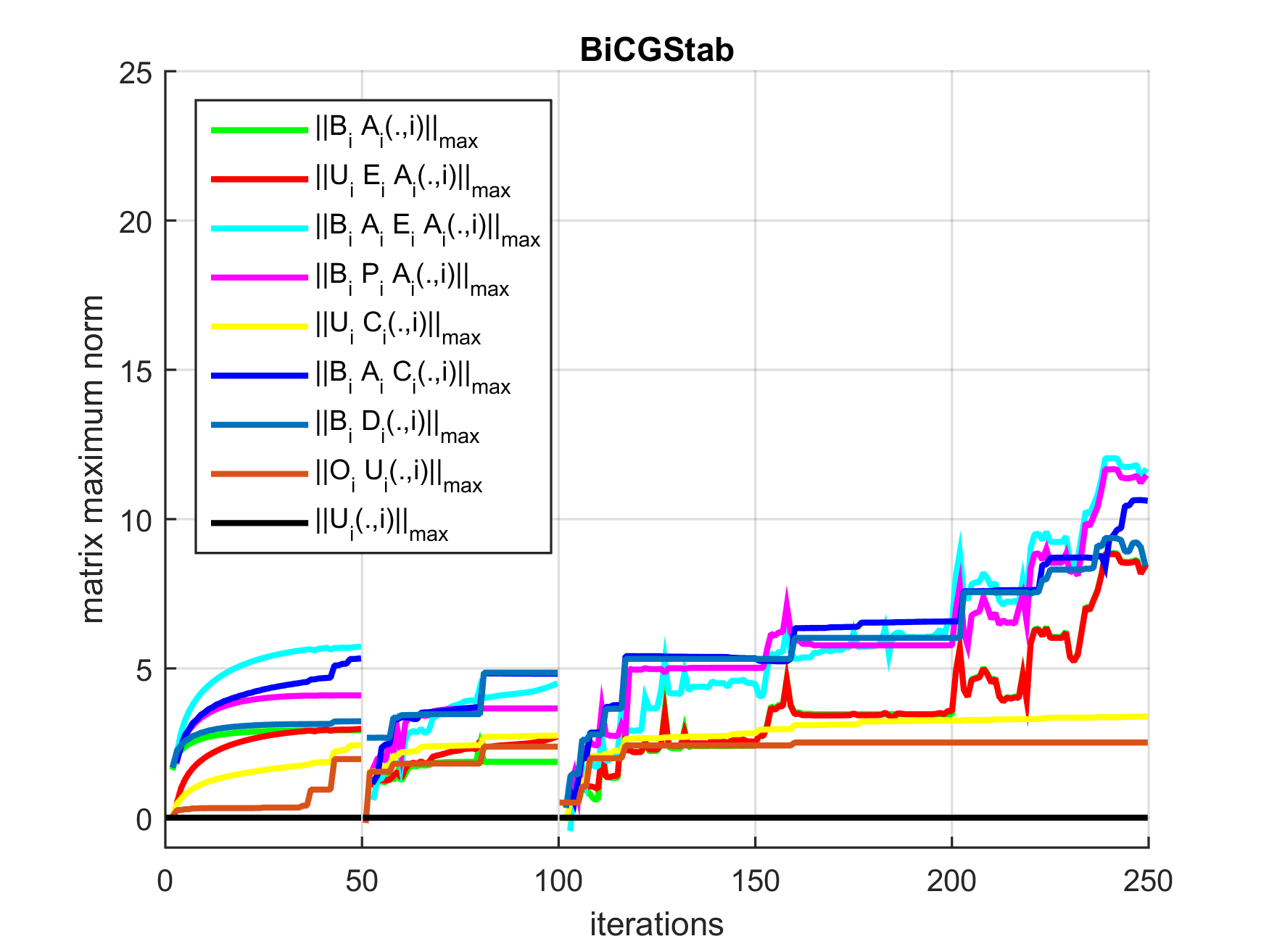} 
\end{center}
\caption{\textbf{(TP1)} Comparison CG (left) vs.~BiCGStab (right). Stabilized version of Fig.~\ref{fig:comparison_ichol} with residual replacement for ICC(0) preconditioned p-CG and p-BiCGStab. Detailed description: see Fig.~\ref{fig:comparison_ichol}. For pipelined CG replacement of the residual $\bar{r}_i$ and auxiliary variables $\bar{u}_i, \bar{w}_i, \bar{s}_i, \bar{q}_i, \bar{z}_i$ by their true values is performed in iterations 50 and 100 (fixed replacement period). For pipelined BiCGStab replacement of the residual $\bar{r}_i$ and auxiliary variables $\bar{s}_i, \bar{\ell}_i, \bar{z}_i, \bar{k}_i, \bar{w}_i$ by their true values, see \eqref{eq:rr}, is performed in iterations 50 and 100 (fixed replacement period). Vertical axis in $\log_{10}$ scale.}
\label{fig:comparison_ichol_rr}
\end{figure}

In Figure \ref{fig:comparison_ichol_rr} (top) benchmark problem (TP1) is again solved, but a residual replacement technique is included to improve the numerical stability of the pipelined algorithms. In the experiment reported in Fig.~\ref{fig:comparison_ichol_rr} replacements are performed periodically every 100 iterations until the criterion $\|\bar{r}_i\|_2 < \sqrt{\epsilon} \, \|\bar{r}_0\|_2$ is satisfied. 
By periodically resetting the accumulated gaps to zero the replacement technique improves numerical stability, resulting in better maximal attainable accuracy for the pipelined algorithms. 

The middle and bottom panels of Figure \ref{fig:comparison_ichol_rr} show the corresponding matrix norms for this experiment. When residual replacement takes place in iteration $i_r$ the accumulated and possibly amplified local rounding errors from previous iterations are eliminated as the residual gap is reset to zero. This implies that the entries in the first $i_r$ rows of the propagation matrices \eqref{eq:all_error_matrices} are set to zero whenever replacements are performed, since the corresponding local rounding errors do not contribute to the gap(s) from iteration $i_r$ onward. 

We remark that no replacements are performed in the classic CG and BiCGStab methods, since replacing the residual only marginally improves the maximal attainable accuracy for these methods. Instead, we use the accuracy attainable by the classic Krylov methods to benchmark the effectiveness of the residual replacement technique for improving the attainable accuracy of the pipelined methods.

By comparing Figure \ref{fig:comparison_ichol_rr} to Figure \ref{fig:comparison_ichol} one observes that the periodic reset of the auxiliary variables 
reduces the size of the maximum norms of the error propagation matrices, and hence improves the final attainable accuracy of the solution. Moreover, the unstable behavior of p-BiCGStab  observed in Figure \ref{fig:comparison_ichol} is reduced, as indicated by the error propagation matrix norms in Figure \ref{fig:comparison_ichol_rr}.

\subsection{Improving numerical stability of p-BiCGStab for unsymmetric and indefinite problems} \label{sec:improving}

\begin{figure}[t]
\begin{center}
\includegraphics[width=0.49\textwidth]{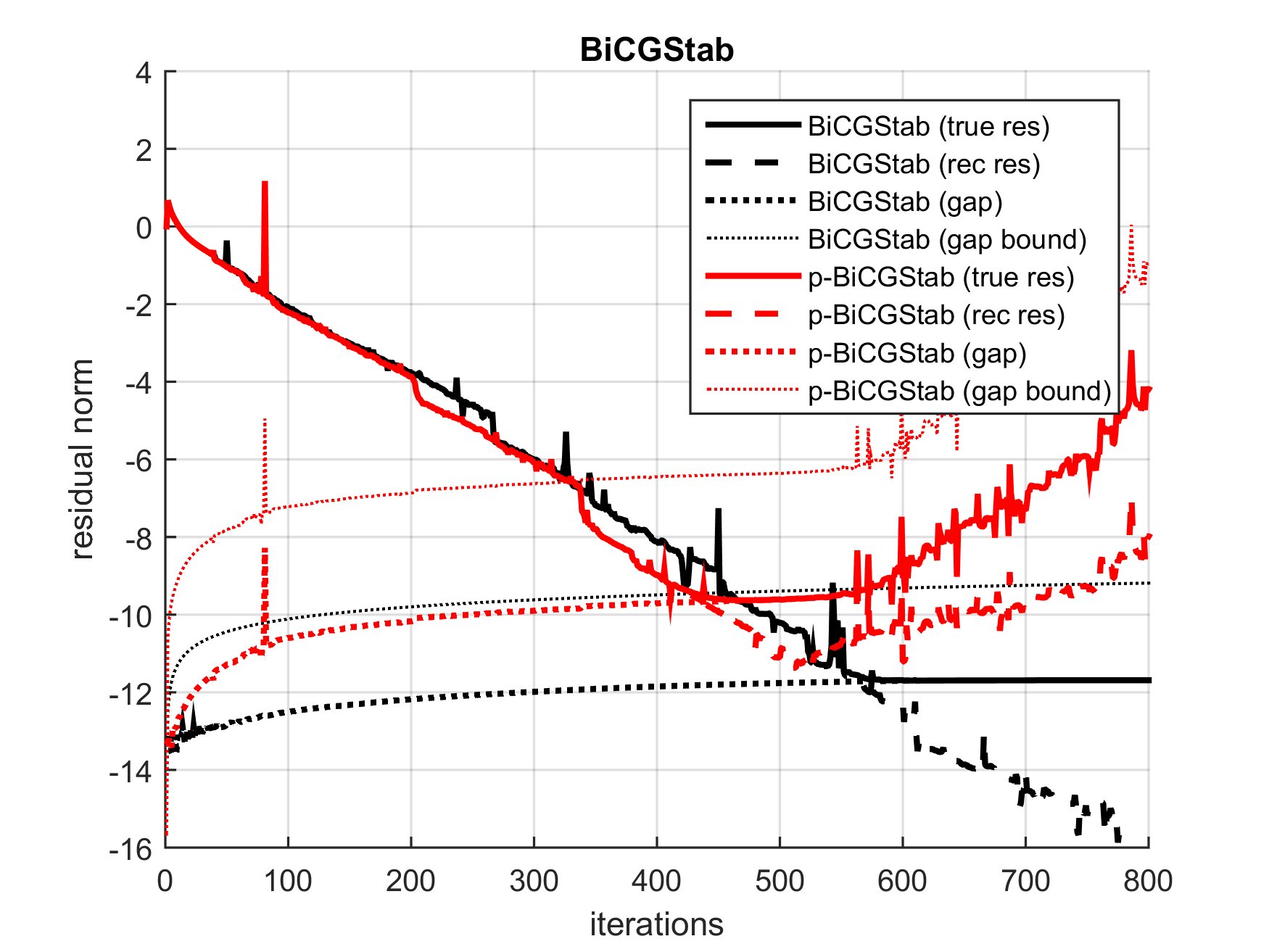}
\includegraphics[width=0.49\textwidth]{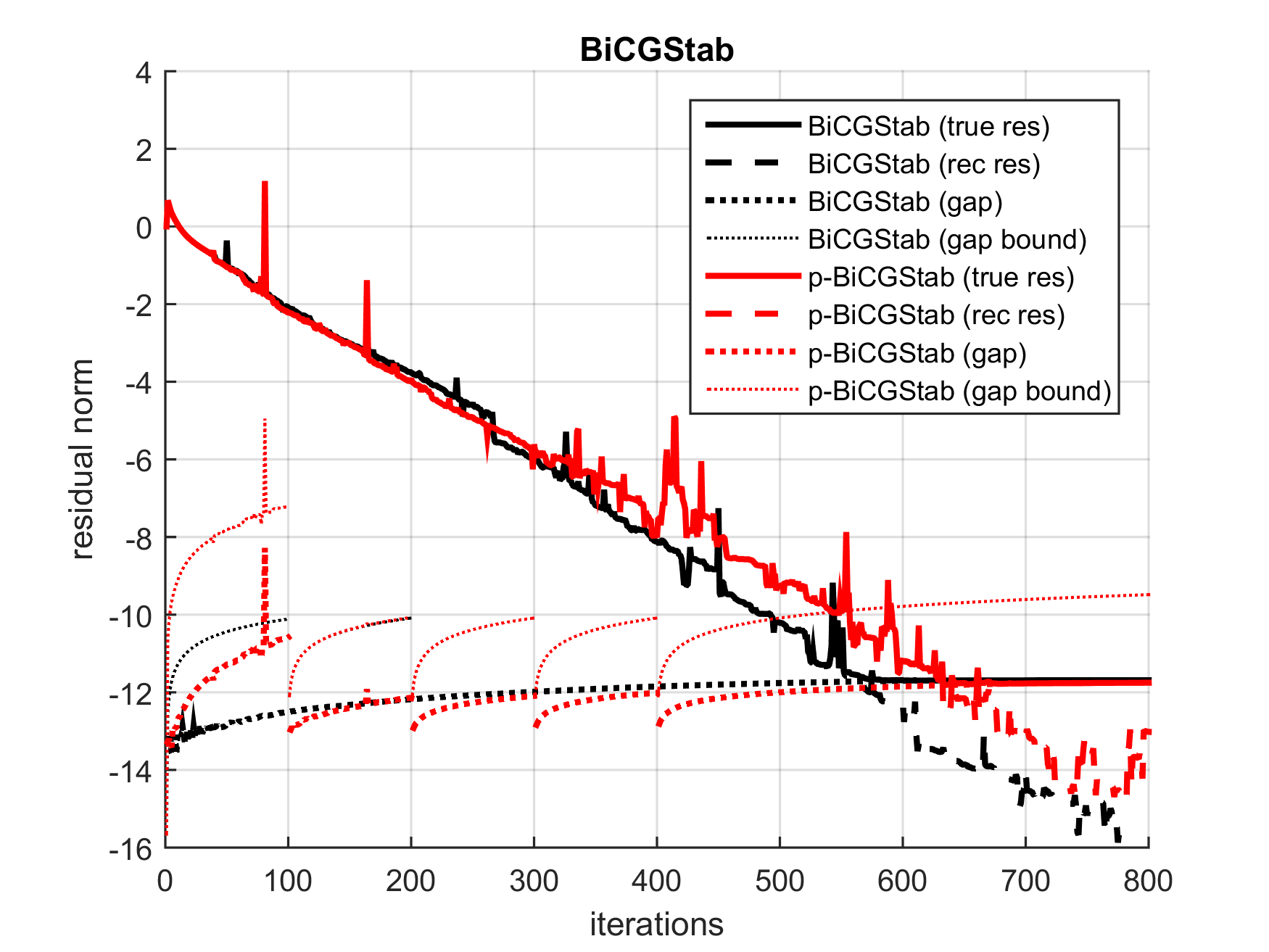} \\
\includegraphics[width=0.49\textwidth]{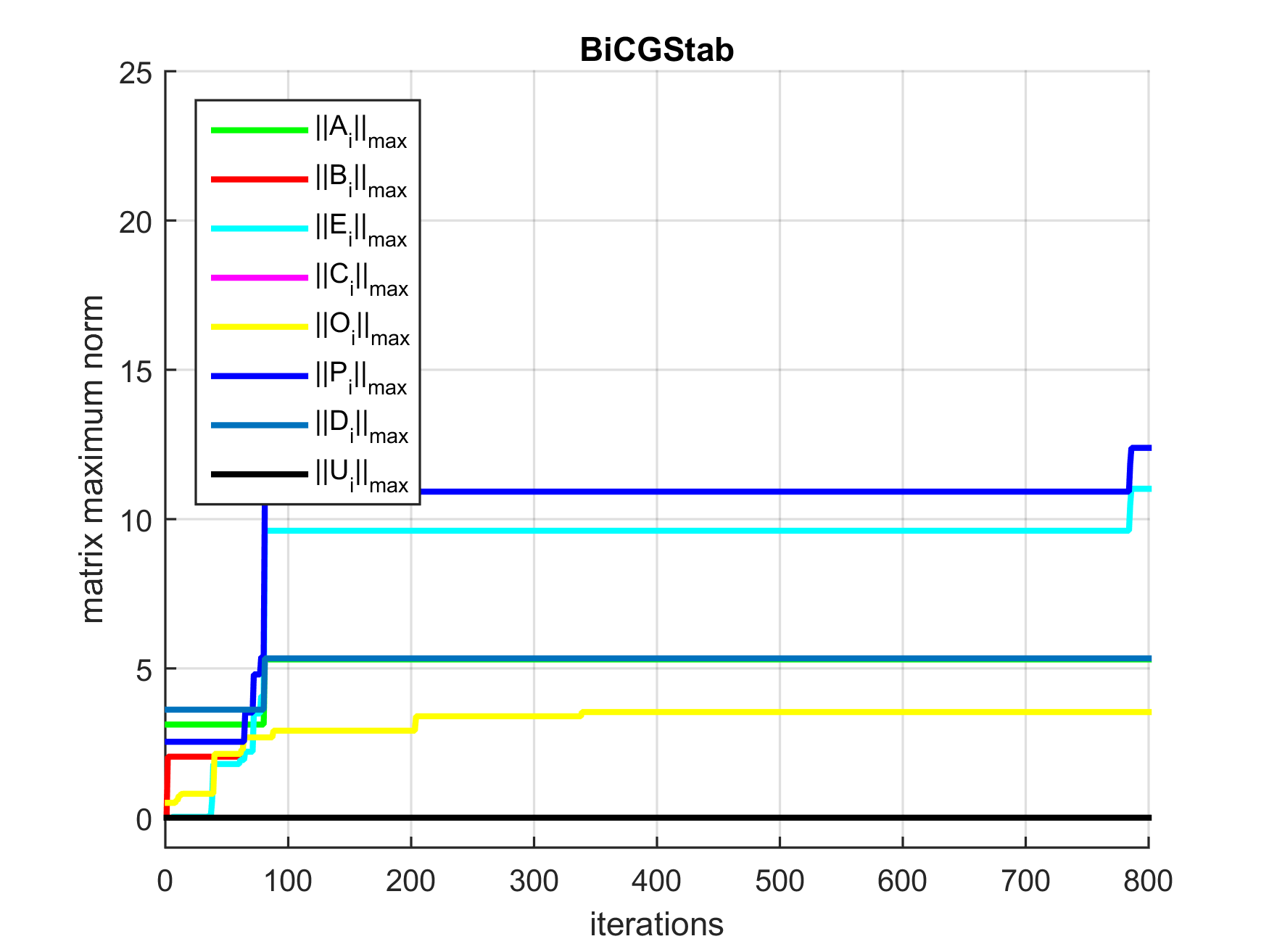}
\includegraphics[width=0.49\textwidth]{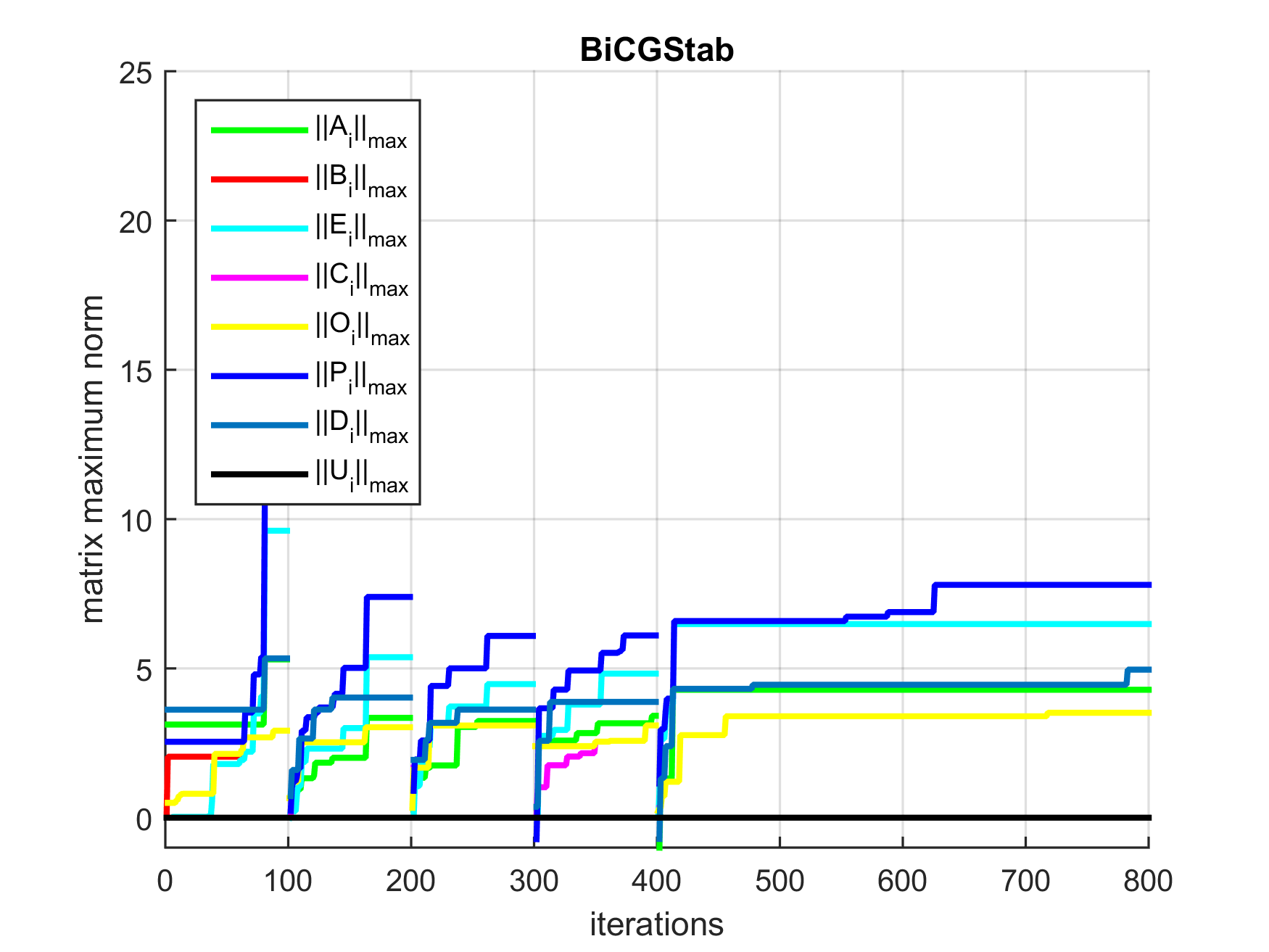} \\
\includegraphics[width=0.49\textwidth]{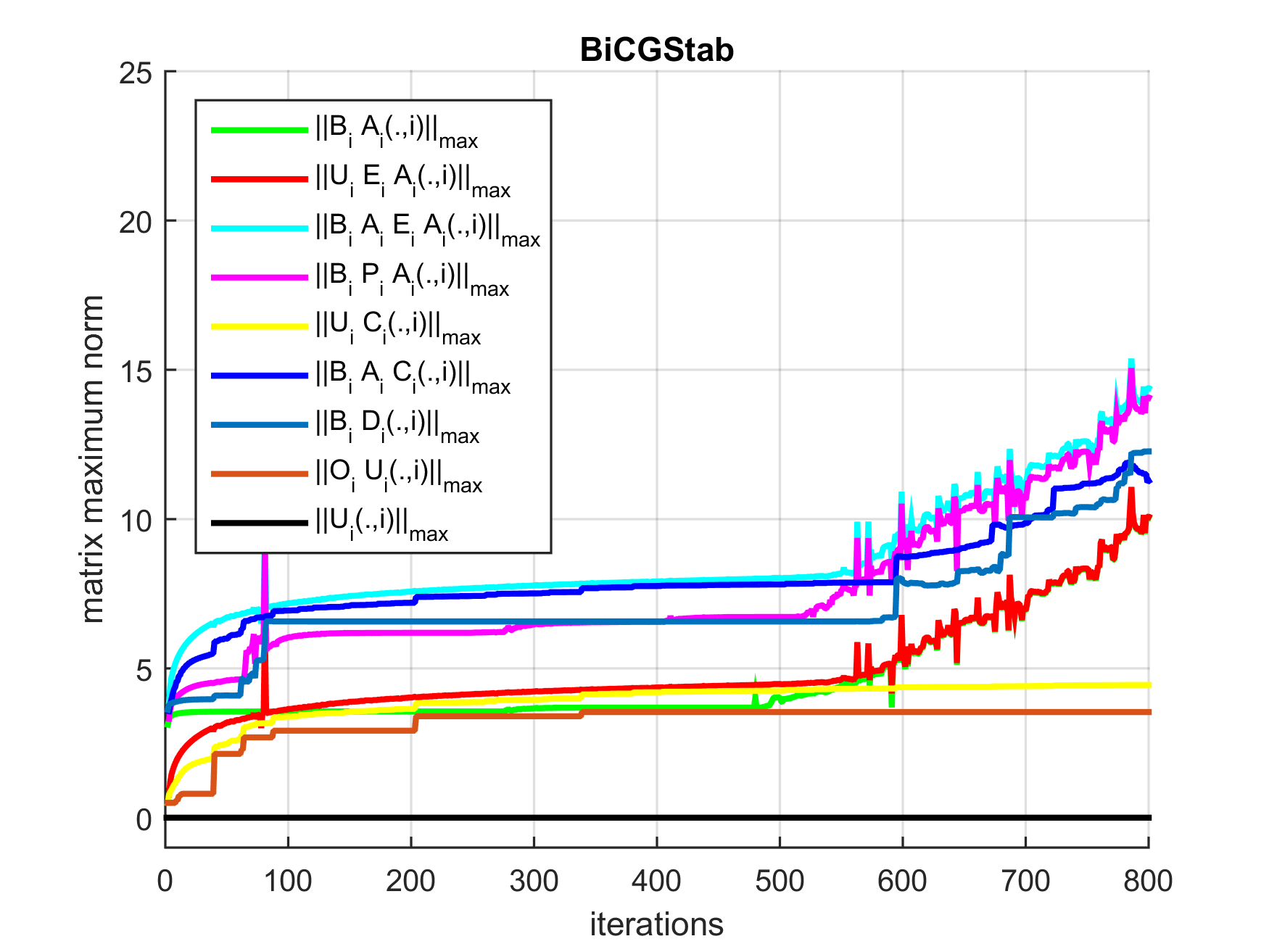}
\includegraphics[width=0.49\textwidth]{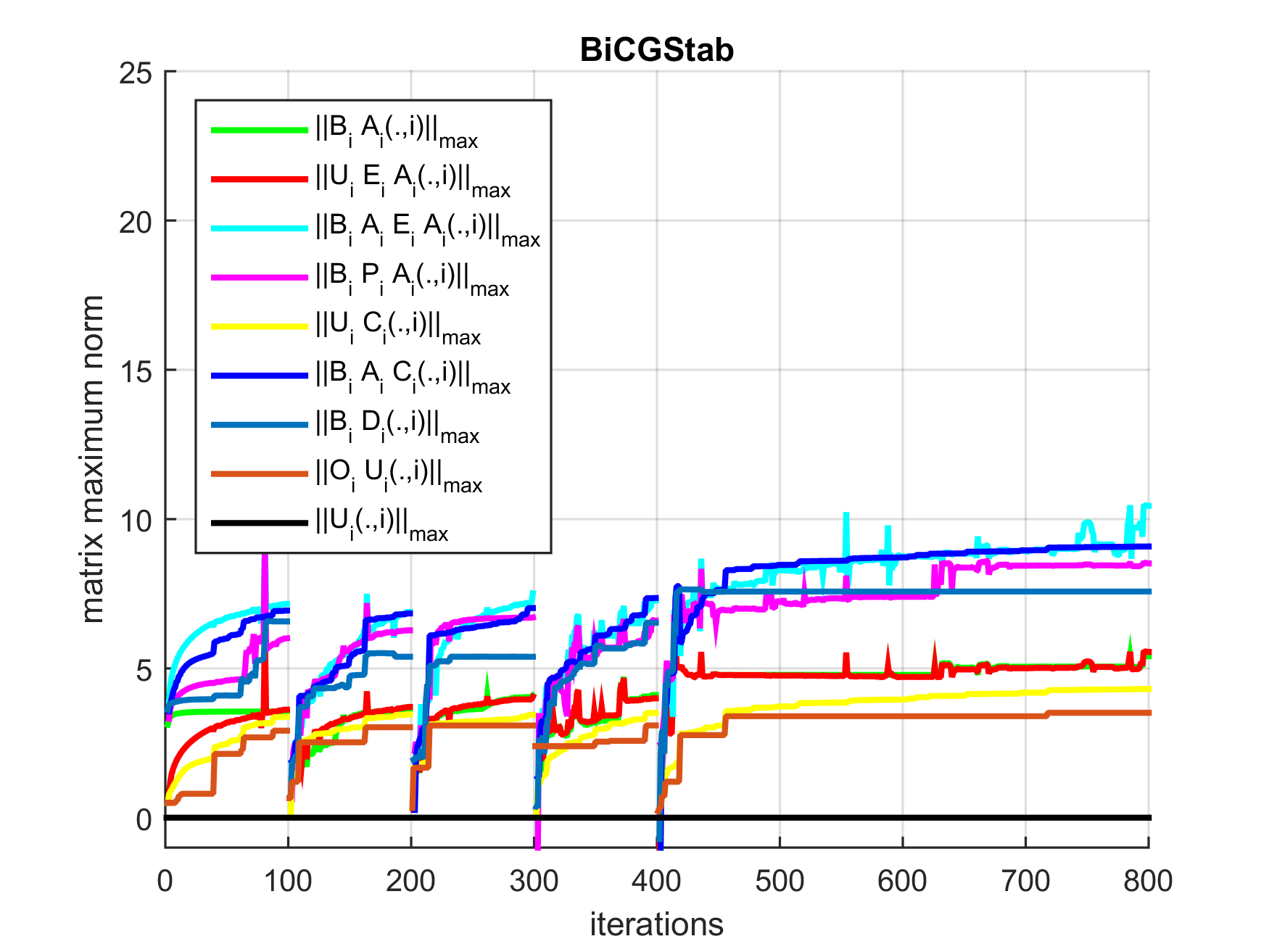} 
\end{center}
\caption{\textbf{(TP2)} Comparison BiCGStab (left) vs.~BiCGStab with residual replacements (right). Fixed replacement period (100 iterations). \textbf{Top:} residual norm history $\|r_i\|_2$ for BiCGStab/p-BiCGStab as a function of iterations. Dotted lines denote the residual gaps $\|(b-A\bar{x}_i) - \bar{r}_i\|_2$ and their computed upper bounds. \textbf{Middle:} maximum norms of various matrices occurring in the numerical stability analysis for p-BiCGStab, see \eqref{eq:matrix_expr}, as a function of iterations. \textbf{Bottom:} maximum norms of the $i$-th column of products of matrices occurring in the stability analysis \eqref{eq:all_error_matrices}. Vertical axis in $\log_{10}$ scale.}
\label{fig:unsymmetric}
\end{figure}

\begin{figure}[t]
\begin{center}
\includegraphics[width=0.49\textwidth]{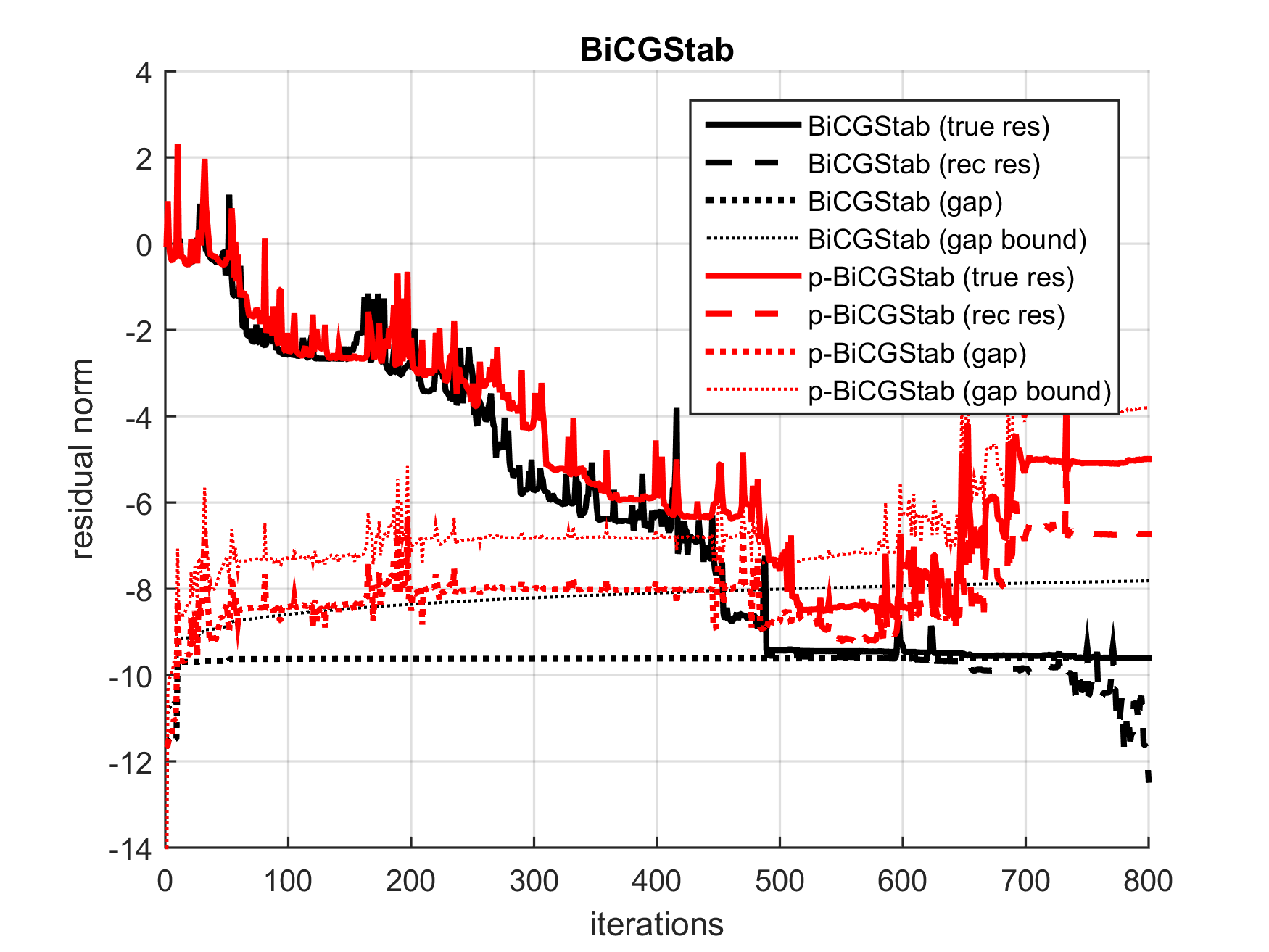}
\includegraphics[width=0.49\textwidth]{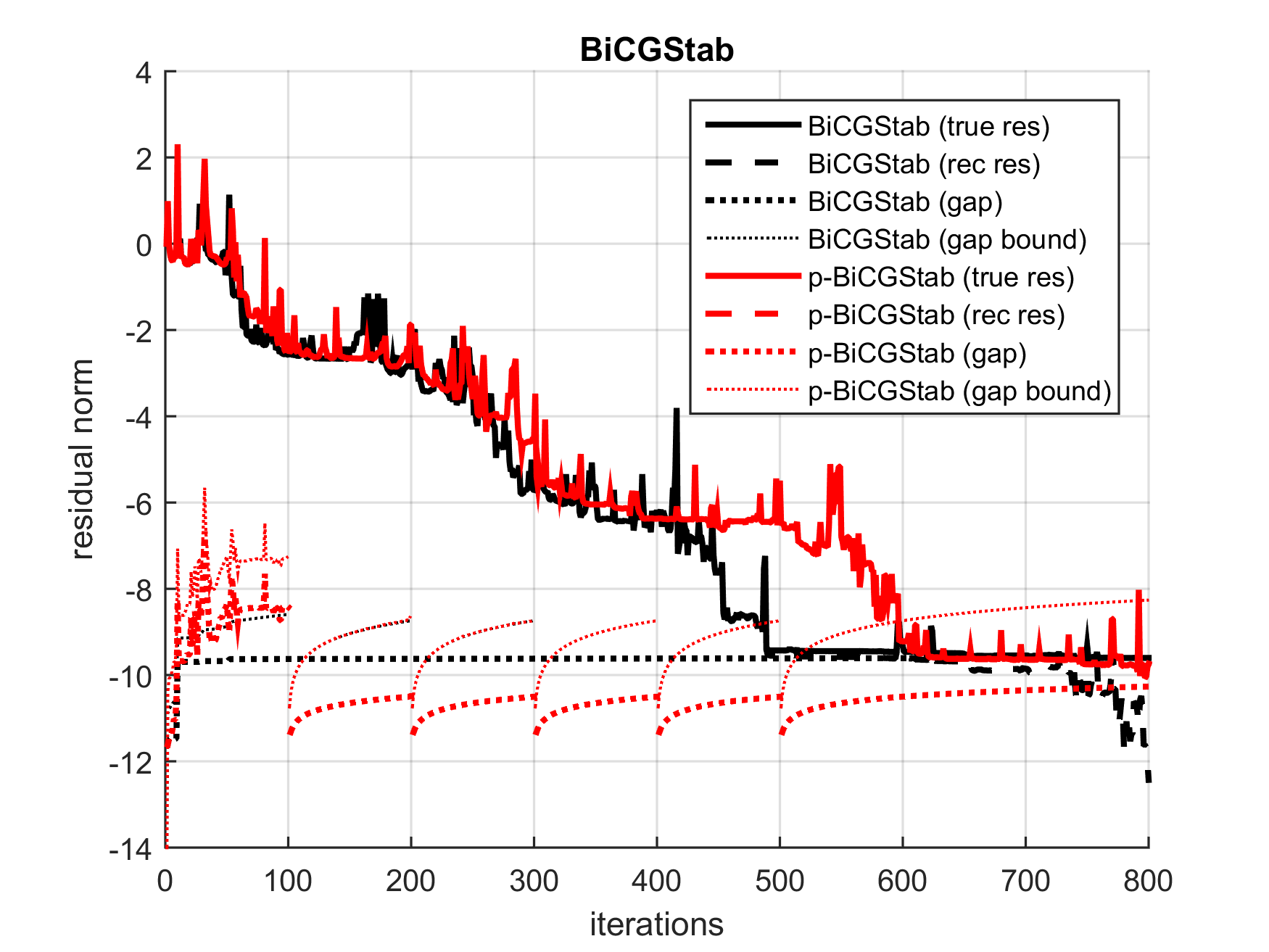} \\
\includegraphics[width=0.49\textwidth]{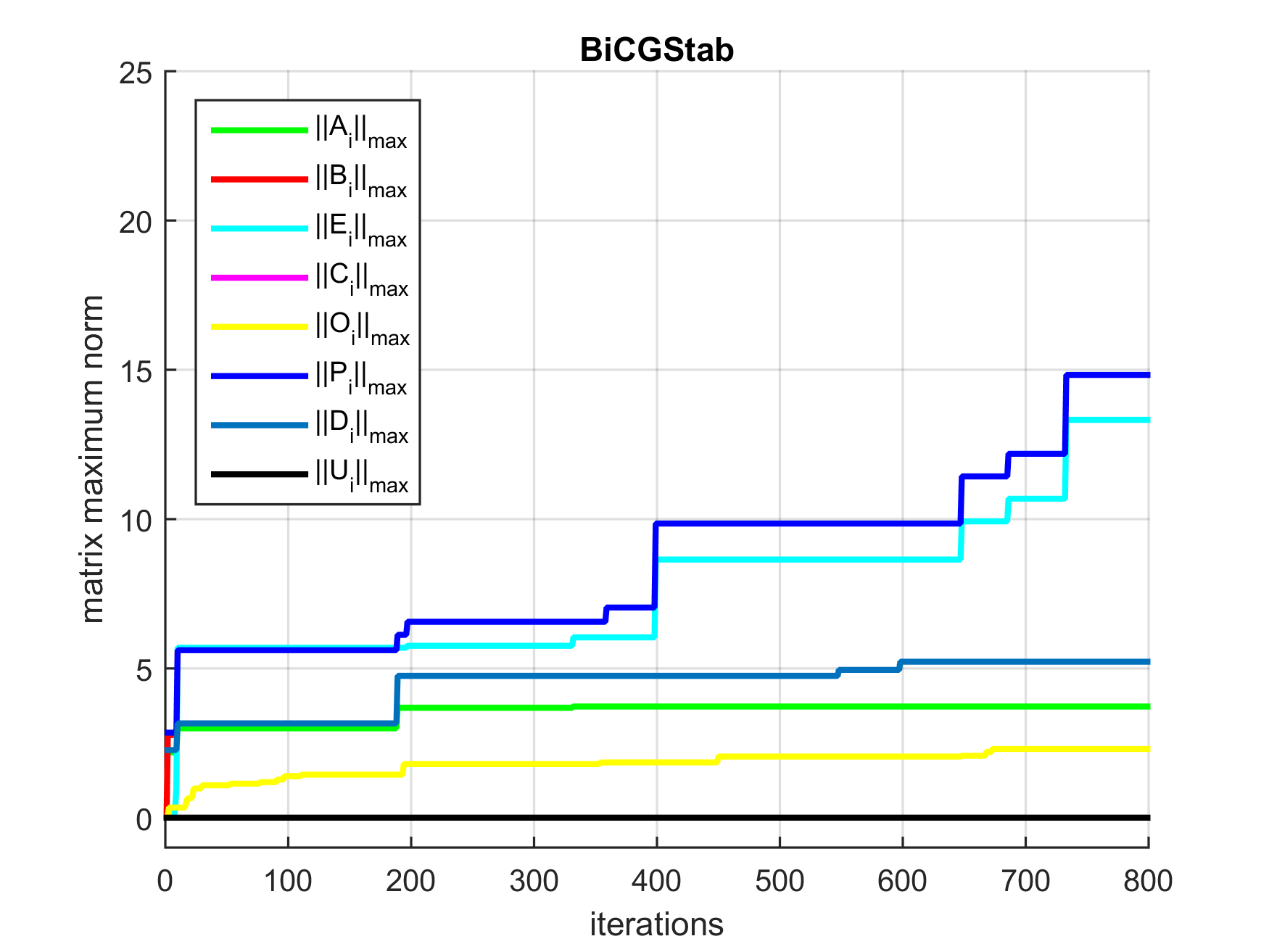}
\includegraphics[width=0.49\textwidth]{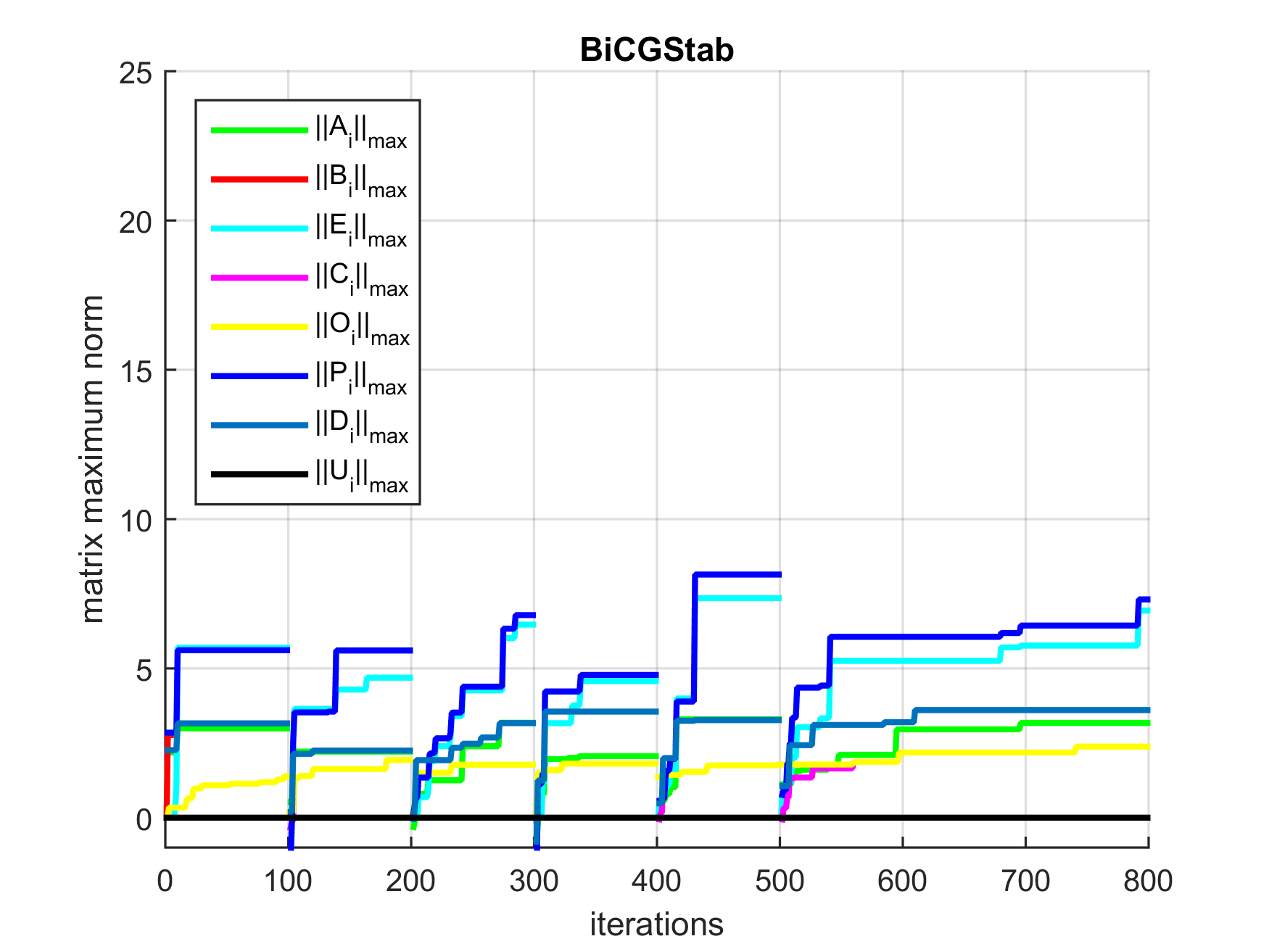} \\
\includegraphics[width=0.49\textwidth]{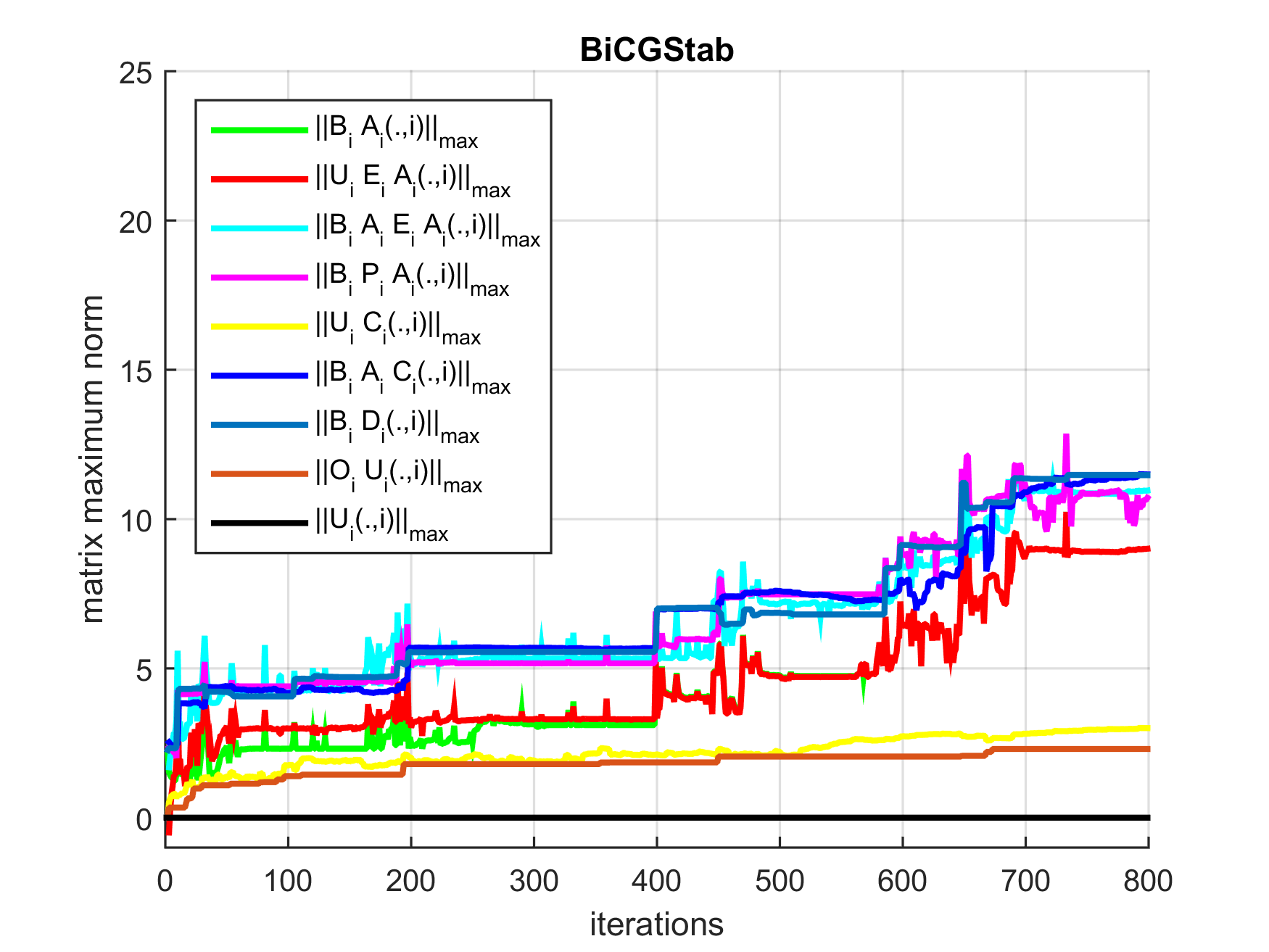}
\includegraphics[width=0.49\textwidth]{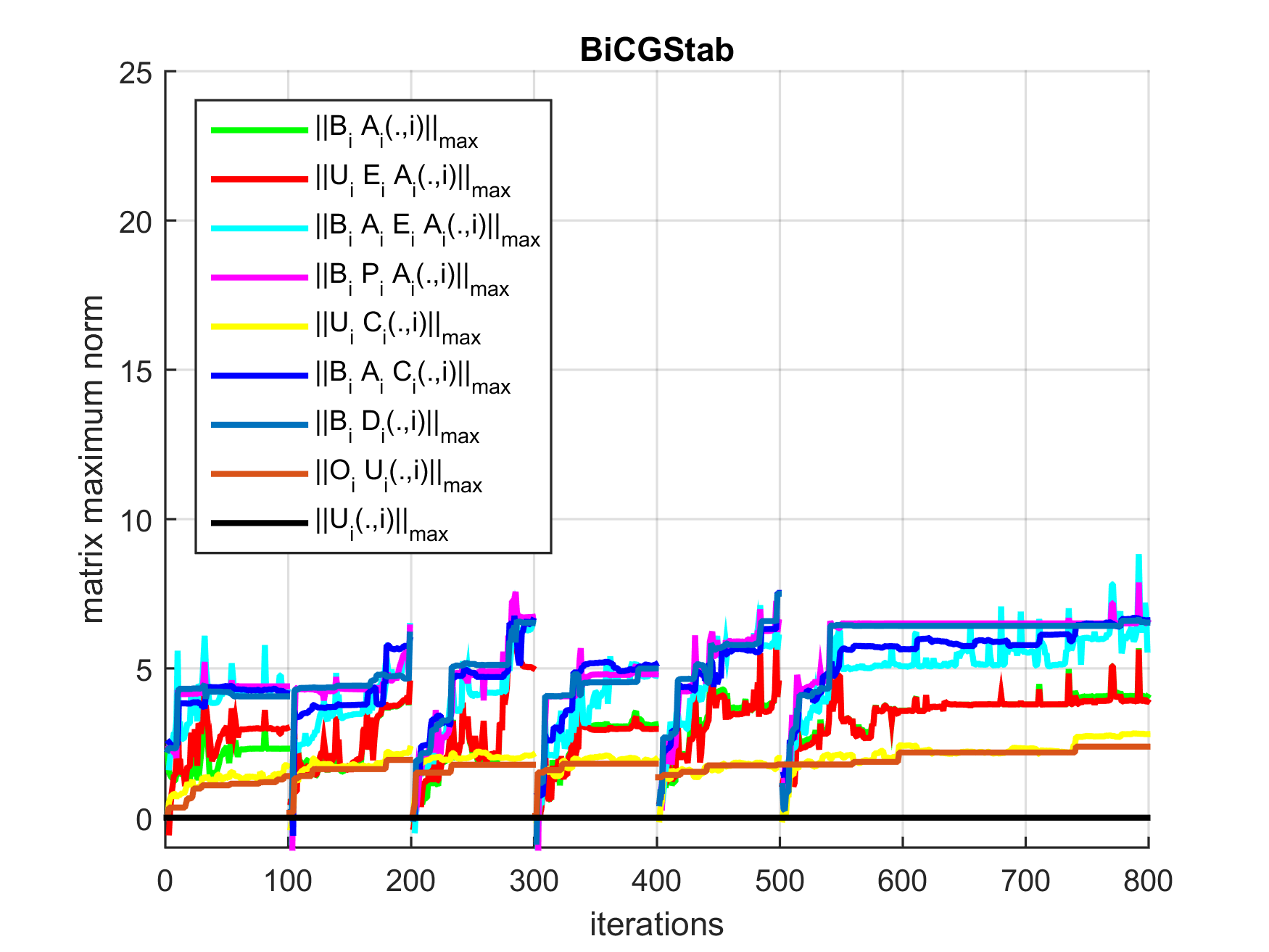} 
\end{center}
\caption{\textbf{(TP3)} Comparison BiCGStab (left) vs.~BiCGStab with residual replacements (right).  Fixed replacement period (100 iterations). \textbf{Top:} residual norm history $\|r_i\|_2$ for ICC(0) preconditioned BiCGStab/p-BiCGStab as a function of iterations. Dotted lines denote the residual gaps $\|(b-A\bar{x}_i) - \bar{r}_i\|_2$ and their computed upper bounds. \textbf{Middle:} maximum norms of various matrices occurring in the numerical stability analysis for p-BiCGStab, see \eqref{eq:matrix_expr}, as a function of iterations. \textbf{Bottom:} maximum norms of the $i$-th column of products of matrices occurring in the stability analysis \eqref{eq:all_error_matrices}.  Vertical axis in $\log_{10}$ scale.}
\label{fig:indefinite}
\end{figure}

\begin{figure}[t]
\begin{center}
\includegraphics[width=0.49\textwidth]{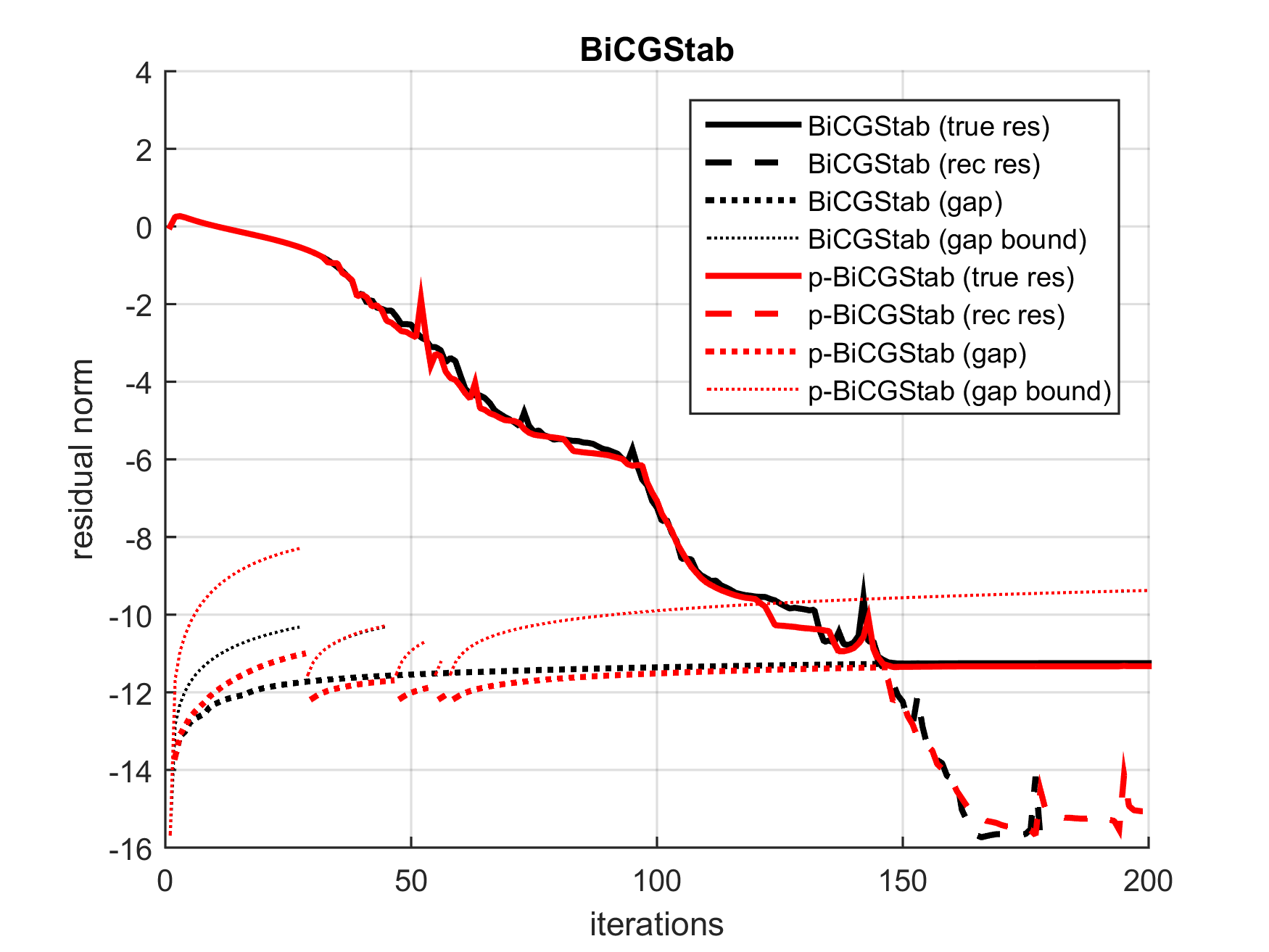}
\includegraphics[width=0.49\textwidth]{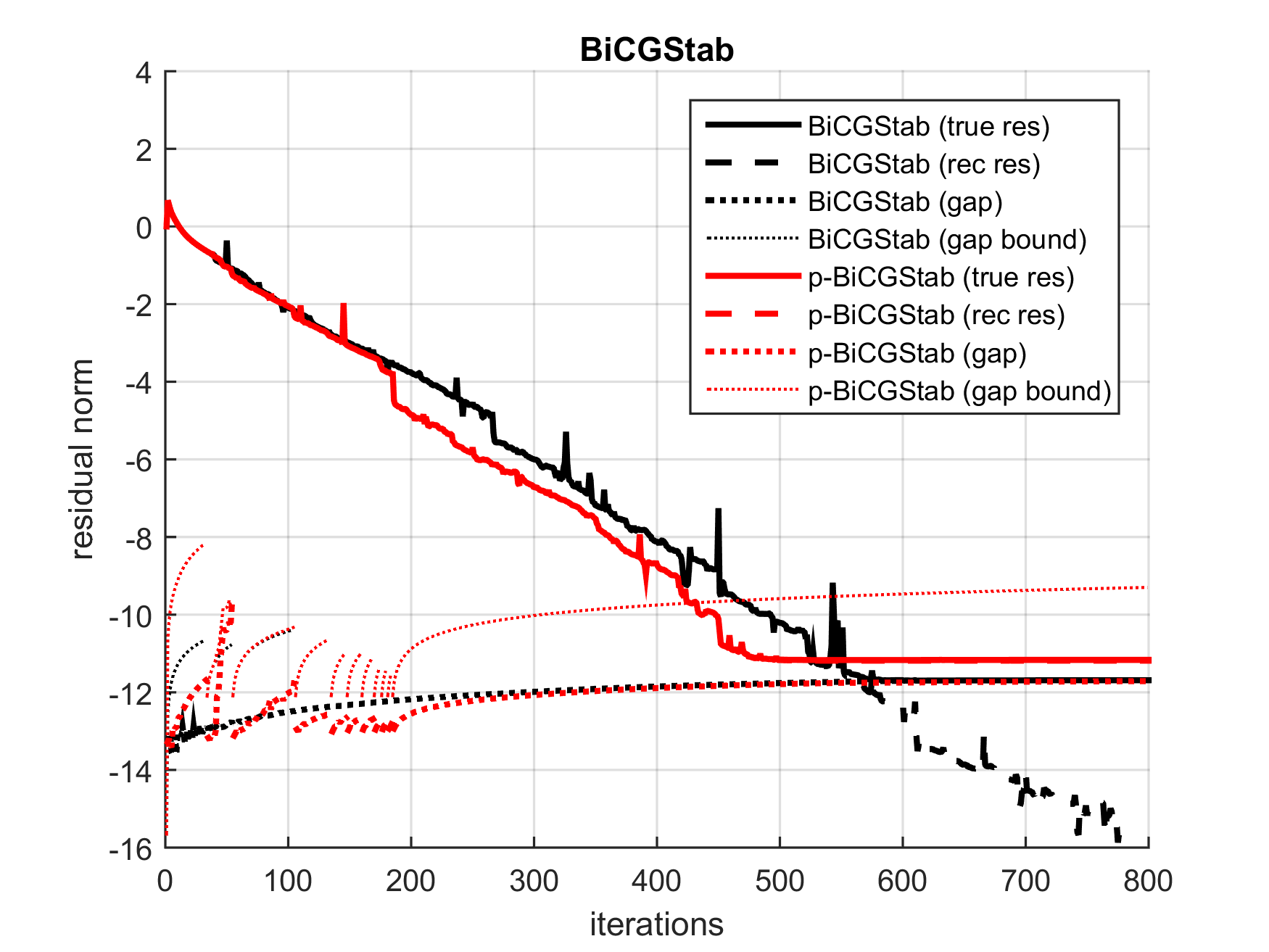} \\
\includegraphics[width=0.49\textwidth]{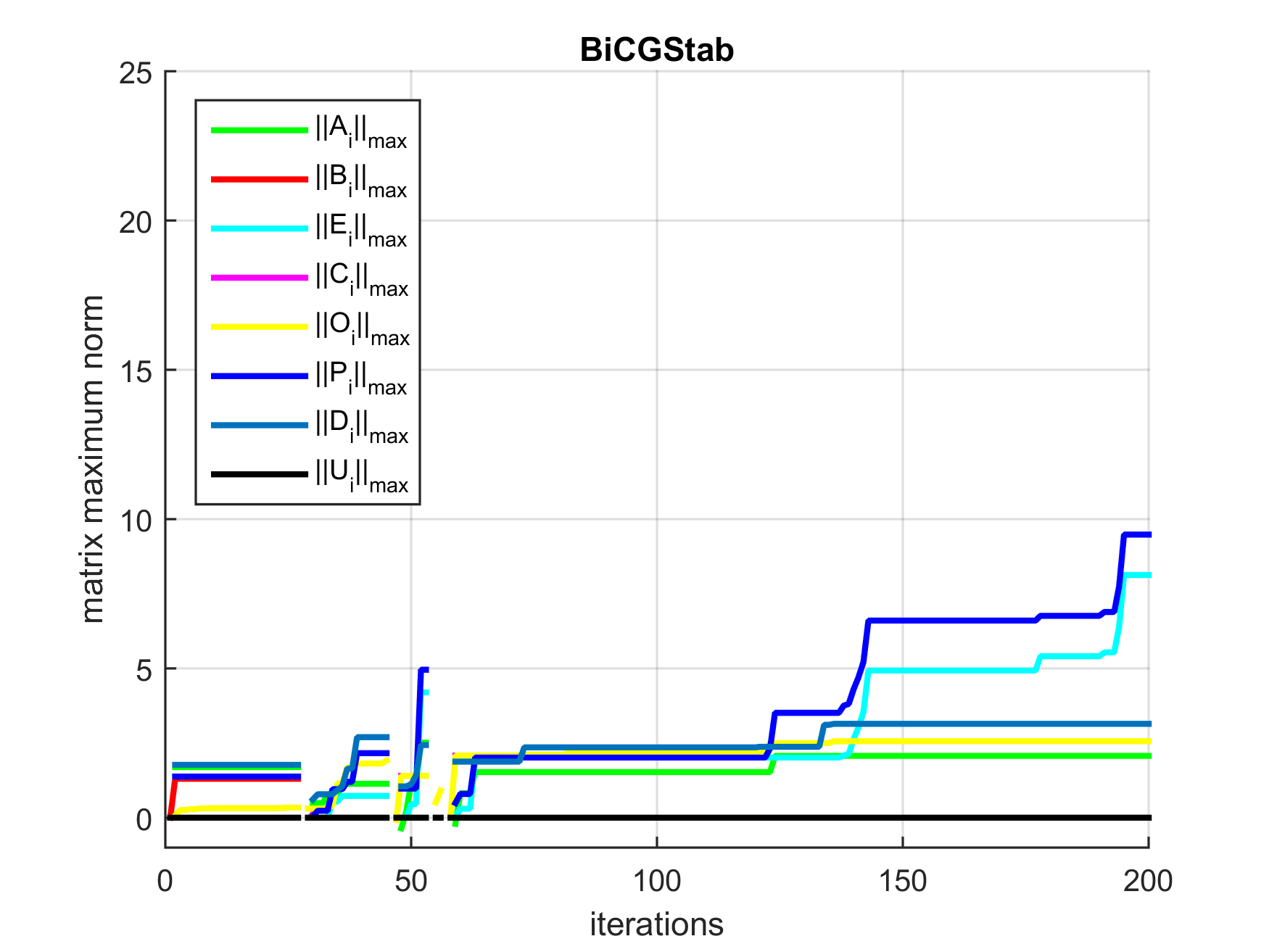}
\includegraphics[width=0.49\textwidth]{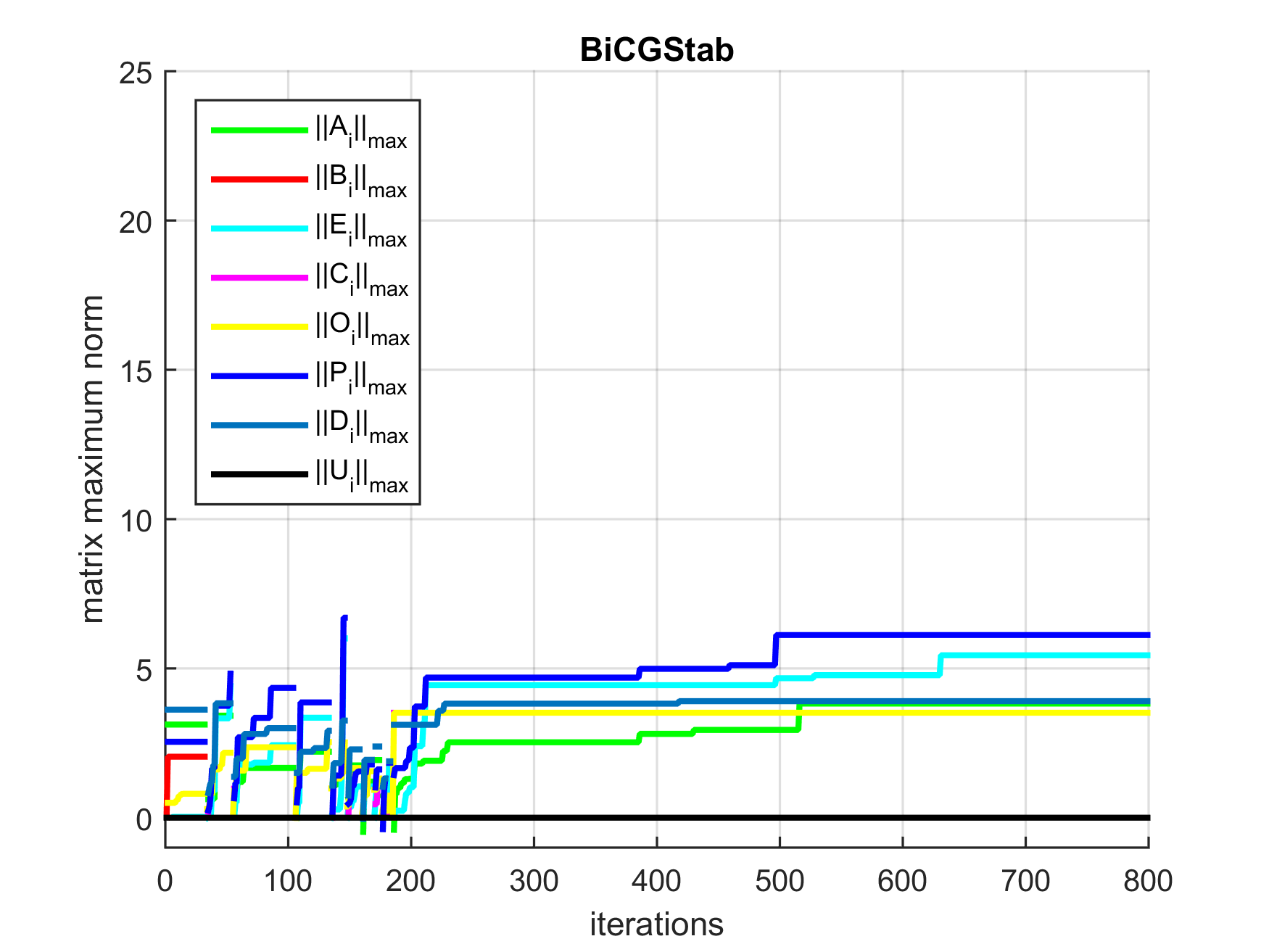} \\
\includegraphics[width=0.49\textwidth]{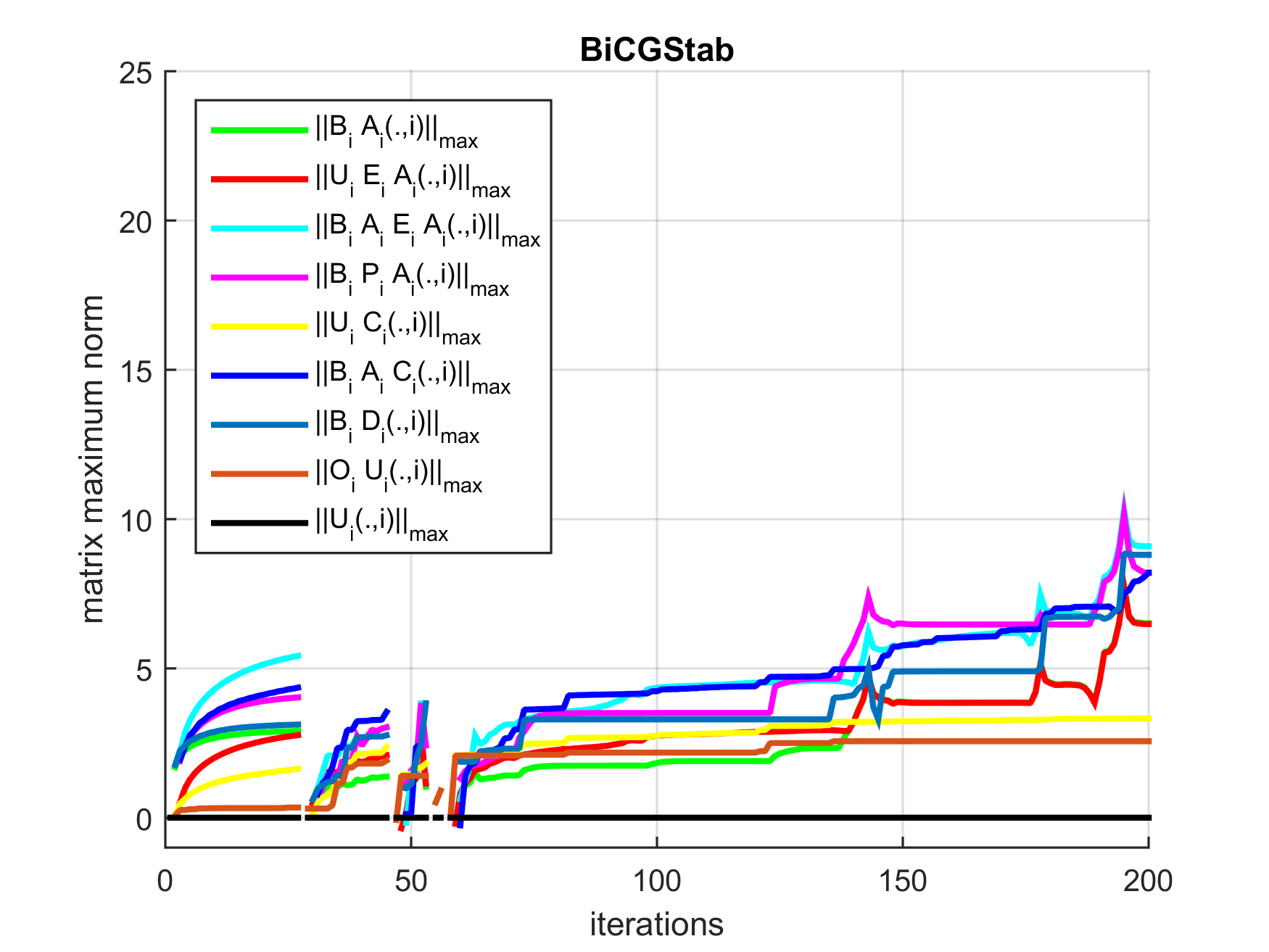}
\includegraphics[width=0.49\textwidth]{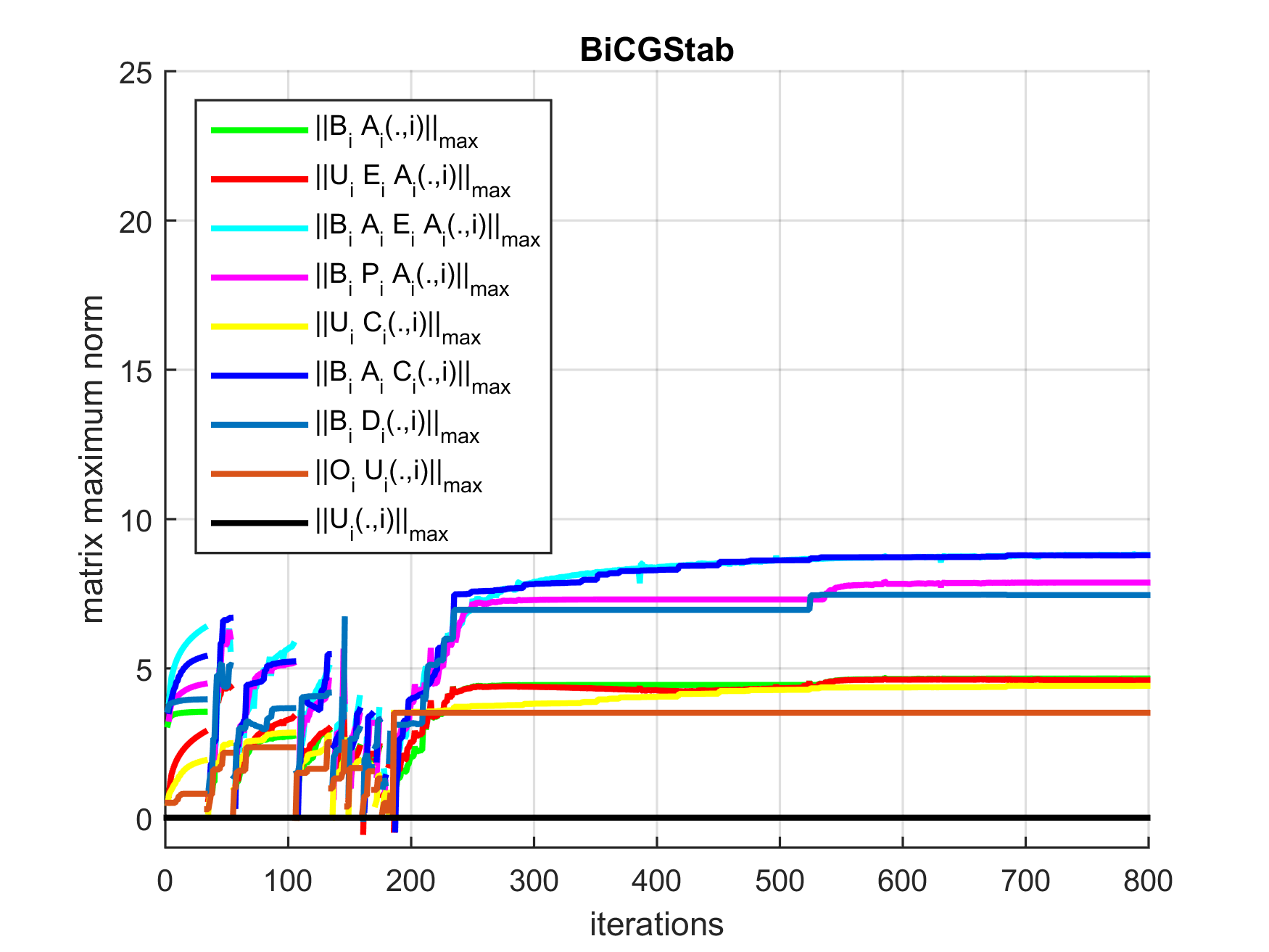} 
\end{center}
\caption{\textbf{(TP1)} (left) \& \textbf{(TP2)} (right). BiCGStab with fully automated residual replacements. Compare to Figs.~\ref{fig:comparison_ichol}-\ref{fig:comparison_ichol_rr} for \textbf{(TP1)} and Fig.\ref{fig:unsymmetric} for \textbf{(TP2)} respectively. \textbf{Top:} residual norm history $\|r_i\|_2$ for BiCGStab/p-BiCGStab as a function of iterations. Dotted lines denote the residual gaps $\|(b-A\bar{x}_i) - \bar{r}_i\|_2$ and their computed upper bounds. \textbf{Middle:} maximum norms of various matrices occurring in the numerical stability analysis for p-BiCGStab, see \eqref{eq:matrix_expr}, as a function of iterations. \textbf{Bottom:} maximum norms of the $i$-th column of products of matrices occurring in the stability analysis \eqref{eq:all_error_matrices}. Vertical axis in $\log_{10}$ scale.}
\label{fig:symmetric_and_unsymmetric}
\end{figure}

Figures \ref{fig:unsymmetric} and \ref{fig:indefinite} show the residual history and propagation matrix norms for the BiCGStab and p-BiCGStab methods on the unsymmetric test problem (TP2) and the indefnite test problem (TP3) respectively. The standard BiCGStab and p-BiCGStab algorithms are presented in the left panels, whereas the right panels show stabilized versions of the methods using the residual replacement strategy. 
Unstable behavior and reduced maximal attainable accuracy are again observed for the p-BiCGStab method compared to standard BiCGStab. Indeed, beyond the p-BiCGStab stagnation point the norms of the error propagation matrices grow exponentially. The numerical analysis from Section \ref{sec:analysis_pipe} provides more insight into the convergence of the pipelined method, as shown by the matrix norms in the bottom panels in Figures \ref{fig:unsymmetric}-\ref{fig:indefinite}. The norms of the matrices $\mathcal{B}_i$ $\mathcal{E}_i$ and $\mathcal{P}_i$, which all contain products of the coefficients $\bar{\beta}_i$, tend to be particularly large. A rapid growth of the product matrices can be observed for p-BiCGStab beyond the stagnation point.

The residual replacement strategy presented in the right panels of Figures \ref{fig:unsymmetric}-\ref{fig:indefinite} reduces the size of the propagation matrix norms, decreasing the detrimental amplification of local rounding errors throughout the algorithm. The true residual norms stagnate several orders of magnitude below the accuracy achievable by p-BiCGStab. Note that the residual replacement strategy may induce a slight delay of convergence due to loss of basis orthogonality, cf.~the discussion in Section \ref{sec:residual}, implying more iterations will sometimes have to be performed to reach a comparable maximal accuracy compared to standard BiCGStab. This can be observed in particular in Fig.~\ref{fig:indefinite} (top). The use of residual replacements thus induces a tradeoff between numerical solution precision and computational cost.

\subsection{A fully automated residual replacement strategy for p-BiCGStab}

As explained in Section \ref{sec:rr} an automated replacement strategy was implemented based on the rounding error analysis from Sections \ref{sec:gap}--\ref{sec:expressing}. Figure \ref{fig:symmetric_and_unsymmetric} presents results for (TP1) (left panel) and (TP2) (right panel) with 
automated residual replacements for the pipelined BiCGStab method. Note that replacements are typically performed early in the iteration
to prevent local rounding errors from propagating. We point out that the number of replacements performed is generally very small (4 and 10 for (TP1) and (TP2) respectively) compared to the total number of iterations required to attain maximal attainable accuracy.
For both test cases the maximal accuracy attained by p-BiCGStab with automated replacements is comparable to that of classic BiCGStab,
while the residual norms appear to stay close to the ones computed by classic BiCGStab.
Furthermore, maximal attainable accuracy is similar to the precision obtained by performing fixed periodic replacements, see Figs.~\ref{fig:comparison_ichol}-\ref{fig:comparison_ichol_rr}.

Figs.~\ref{fig:9point}-\ref{fig:7point} demonstrate the practical useability of the automated replacement strategy based on the bound \eqref{eq:DDr_bound} and the criterion \eqref{eq:criterion} on two additional test problems, (TP4) and (TP5). The left panels show the convergence histories and the error propagation matrices for classic BiCGStab and pipelined BiCGStab without replacements, while the right panels detail the same information for p-BiCGStab with automated residual replacement. For both benchmark problems three replacement steps are performed, causing only very limited computational overhead. The fully automated residual replacement strategy improves the maximal attainable accuracy and stability of the residuals after maximal accuracy has been attained. In both cases nearly identical convergence behavior to classic BiCGStab is obtained by p-BiCGStab when using automated replacements.

\subsection{Performance of p-BiCGStab for practical benchmark cases} \label{sec:performance}

\begin{figure}[t]
\begin{center}
\begin{tabular}{cc}
\includegraphics[width=0.48\textwidth]{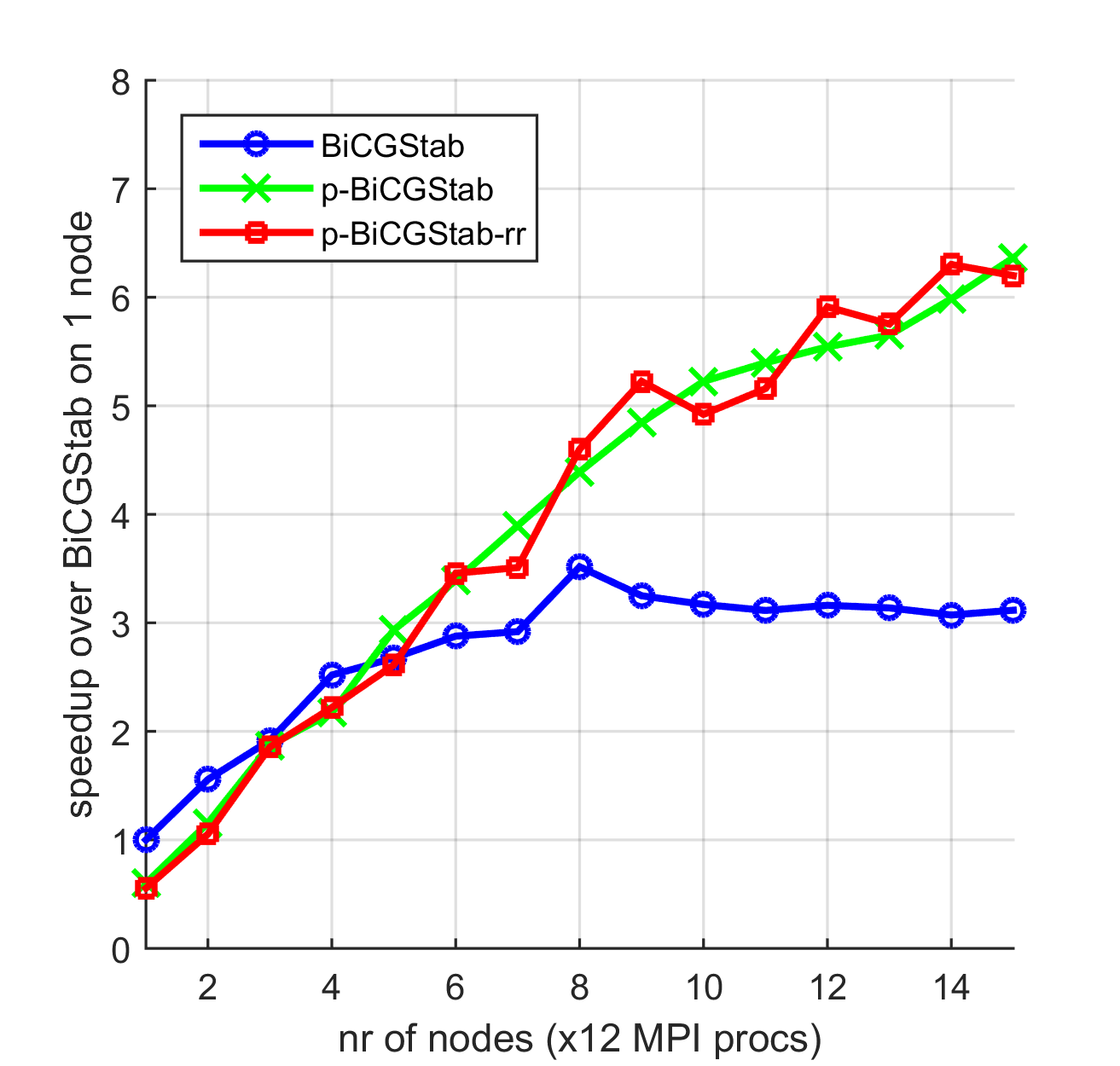} & 
\includegraphics[width=0.48\textwidth]{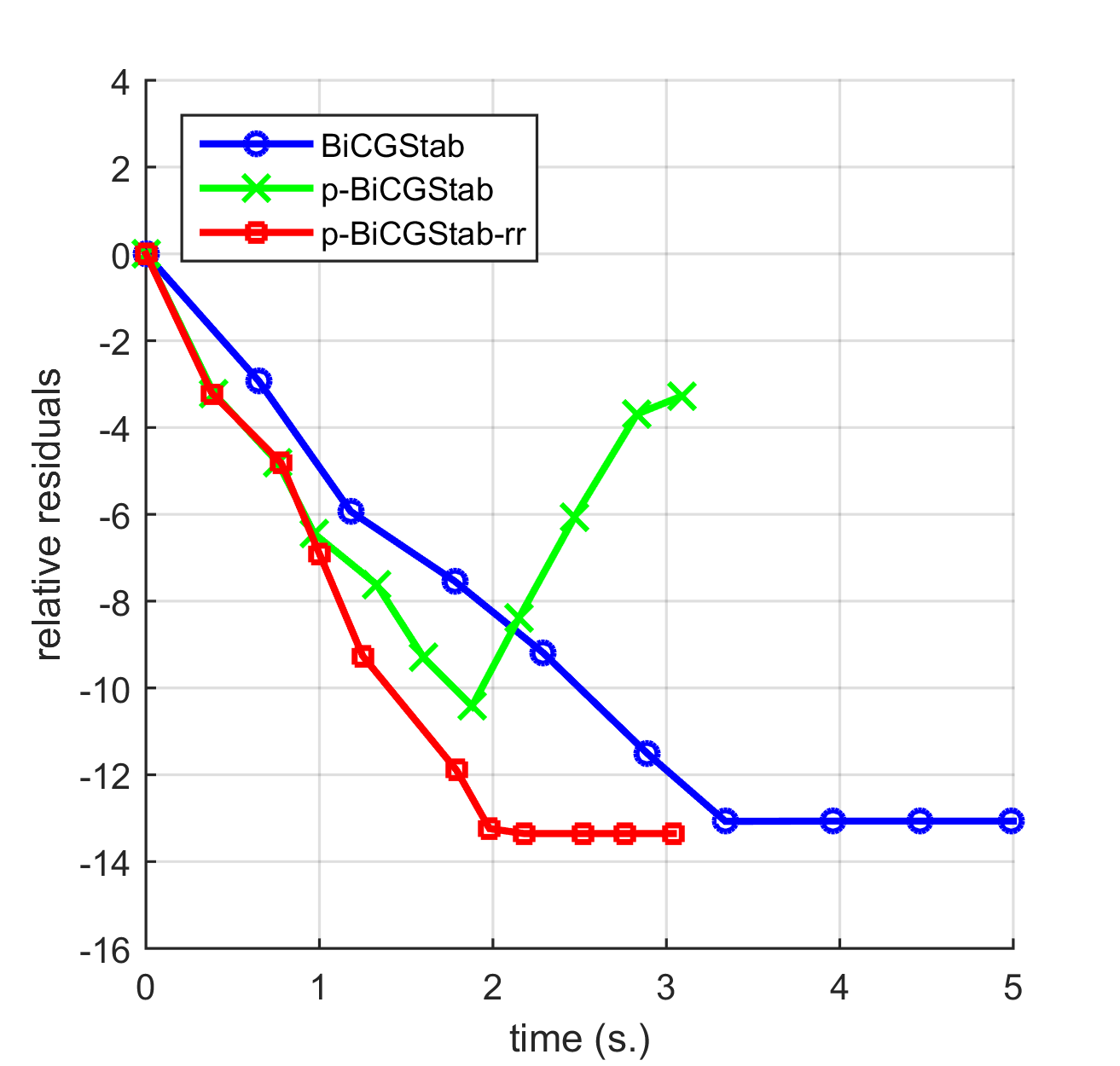} 
\end{tabular}
\end{center}
\caption{\textbf{(TP2)} Strong scaling experiment on up to 15 nodes (180 cores). 
\textbf{Left:} Speedup (based on average time per iteration; i.e.~total time divided by the number of iterations) over standard BiCGStab on a single node. All methods converged to a scaled residual tolerance of $10^{-6}$, which was reached on average in 249 iterations (205 min./282 max.). 
\textbf{Right:} Relative true residual norm (\texttt{log10} scale) as a function of total time spent by the algorithm on 10 nodes. Measurement points are based on 100 iteration intervals.
The p-BiCGStab-rr algorithm performs a replacement step every 100 iterations.  Vertical axis in $\log_{10}$ scale.
\label{fig:timings2}}
\end{figure}

\begin{figure}[t]
\begin{center}
\begin{tabular}{cc}
\includegraphics[width=0.48\textwidth]{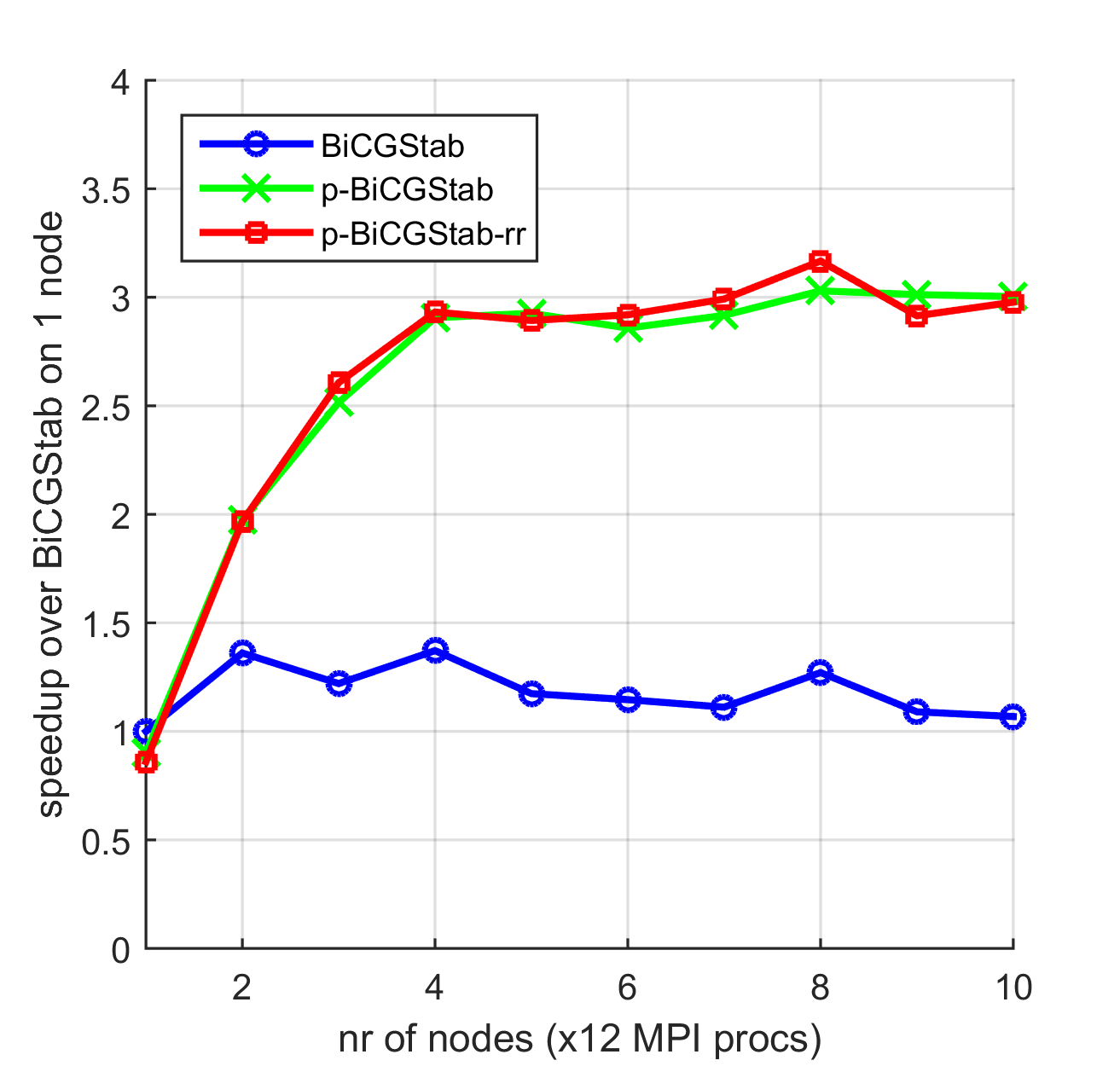} & 
\includegraphics[width=0.48\textwidth]{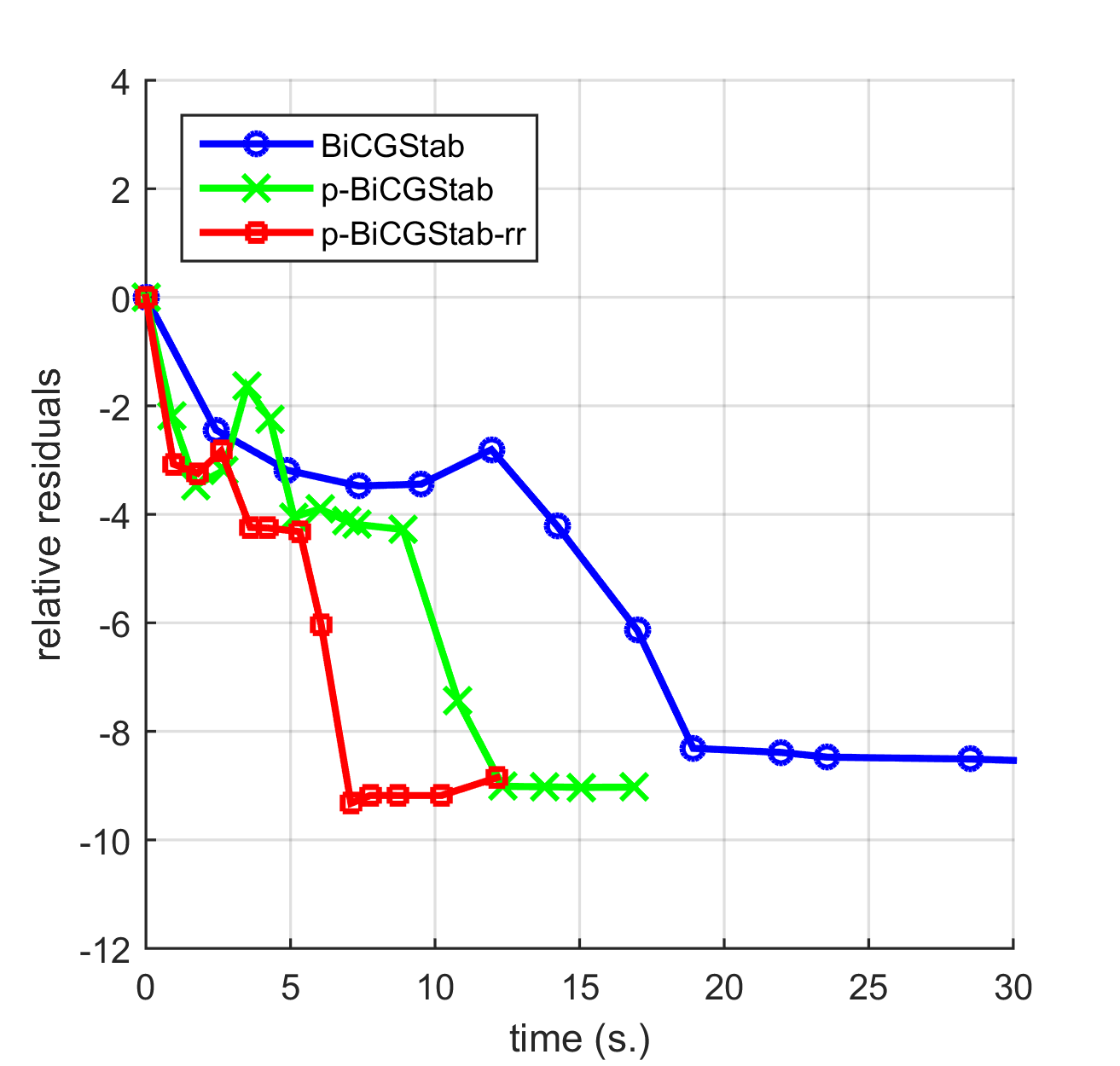} 
\end{tabular}
\end{center}
\caption{\textbf{(TP3)} Strong scaling experiment on up to 10 nodes (120 cores). 
\textbf{Left:} Speedup (based on average time per iteration; i.e.~total time divided by the number of iterations) over standard BiCGStab on a single node. All methods converged to a scaled residual tolerance of $10^{-6}$, which was reached on average in 4481 iterations (2784 min./7559 max.). 
\textbf{Right:} Relative true residual norm (\texttt{log10} scale) as a function of total time spent by the algorithm on 10 nodes. Measurement points are based on 500/1.000 iteration intervals (before/after iteration 5.000 respectively).
The p-BiCGStab-rr algorithm performs a replacement step every 200 iterations.  Vertical axis in $\log_{10}$ scale.
\label{fig:timings3}}
\end{figure}

Figures \ref{fig:timings2} and \ref{fig:timings3} show strong scaling results (left panel) and accuracy experiments (right panel) for the benchmark problems (TP2) and (TP3) respectively. The primary aim of these figures is to illustrate that the performance of p-BiCGStab is not negatively affected by performing replacements, whereas attainable accuracy is increased significantly. It is observed from Figure \ref{fig:timings2} (left) that for this specific problem setup and hardware configuration p-BiCGStab is able to achieve a speedup of roughly $2\times$ over classic BiCGStab when both are executed on 15 nodes. Note that classic BiCGStab stops scaling at around 8 nodes in this experiment whereas p-BiCGStab scales quite well on up to 15 nodes. Figure \ref{fig:timings3} (left) shows that for the smaller (and thus relatively more latency dominated) problem (TP3) classic BiCGStab does not scale at all for this setup, whereas p-BiCGStab scales well on up to 5 nodes, after which speedup stagnates. The p-BiCGStab method achieves a speedup of approximately $3\times$ over BiCGStab for (TP3) when both are executed on 10 nodes. 

The right panels of Figures \ref{fig:timings2}-\ref{fig:timings3} show the attainable accuracy (given by the relative true residual norm $\|b-A\bar{x}_i\|/\|b\|$) as a function of the total time spent by the algorithm on 10 nodes. For the unsymmetric test problem (TP2) shown in Figure \ref{fig:timings2} it is observed that classic BiCGStab is able to achieve a relative residual norm in the order of 1e-13 in 3.34 seconds ($\sim$600 iterations). Pipelined BiCGStab converges significantly faster, but is only able to reach an accuracy of the order of 1e-10 in 1.88 seconds ($\sim$600 iters), after which the residual norms start to increase. 
By using the residual replacement strategy for p-BiCGStab a relative residual norm of the order 1e-13 can be obtained in 1.98 seconds ($\sim$600 iters) for this test problem.

Figure \ref{fig:timings3} indicates that the time gained by using p-BiCGStab with residual replacements compared to classic BiCGStab is even more pronounced for (TP3). Classic BiCGStab is able to attain a relative residual norm in the order of 1e-8 in 18.93 s.~($\sim$4,000 iters), whereas p-BiCGStab reaches a relative residual norm in the order of 1e-9 in 12.34 s.~(requiring $\sim$7,000 iters) while after 6.93 s.~($\sim$4,000 iterations) a relative residual norm of only 8e-5 is reached by p-BiCGStab. The p-BiCGStab algorithm with residual replacements attains a similar precision to classic BiCGStab in only 7.09 s.~($\sim$4,000 iters). Note that for (TP3) convergence of the p-BiCGStab method stagnates around a residual norm of 1e-4 early on in the iteration, but succeeds in converging to a precision comparable to classic BiCGStab eventually. 
The residual replacement technique improves the speed of convergence of p-BiCGStab by 
removing the numerical instabilities which cause the early plateau.


\section{Conclusions} \label{sec:conclusions} 

This paper analyzes the effect of local rounding errors on the 
attainable accuracy of the communication hiding pipelined BiCGStab method \cite{cools2017communication} in finite precision arithmetic. In analogy to the work in \cite{cools2018analyzing} a theoretical framework is derived that characterizes the propagation of rounding errors stemming from the recurrence relations used in the algorithm. These coupled multi-term recurrence relations allow for improved parallel scalability on multi-node parallel hardware, yet significantly increase the algorithm's sensitivity to numerical round-off errors. 

The source of the error propagation in pipelined BiCGStab is designated by the error propagation model and --matrices introduced 
in this work. Depending on the magnitude of the (products of) scalar coefficients defined in each iteration of the algorithm, the amplification of local rounding errors throughout the pipelined BiCGStab algorithm may lead to significantly reduced precision on the solution and even cause highly unstable behavior of the corresponding residuals.

The propagation of local rounding errors in pipelined BiCGStab is compared to that of classic BiCGStab and to the related pipelined CG algorithm on a symmetric model problem. Furthermore, numerical experiments validate the numerical analysis performed in the first part of the paper, and substantiate theoretically the use of residual replacement type techniques that are proposed to improve the stability of the pipelined BiCGStab method. An automated replacement strategy is proposed based on the error analysis. Parallel performance and accuracy experiments illustrate the practical usefulness of the pipelined BiCGStab method and the automated residual replacement technique for solving large scale linear systems on parallel distributed memory hardware.

The analysis in this work focuses on understanding the impact of local rounding errors in the multi-term recurrence relation pipelined BiCGStab algorithm on the attainable precision of the iterative solution. Other sources of errors, such as rounding errors due to loss of basis orthogonality \cite{liesen2012krylov}, system noise related errors \cite{morgan2016stochastic} or hard faults/soft errors \cite{agullo2016hard} are not considered in this study. These topics will be considered in the context of pipelined Krylov methods as part of future work.

\section{Acknowledgments} 

\noindent S.\,Cools is funded by the Research Foundation Flanders (FWO) under research grant number 12H4617N. The author would like to thank Jeffrey Cornelis and Wim Vanroose (University of Antwerp) for valuable discussions on the contents of this manuscript.

\setlength{\bibsep}{0pt plus -1.5ex}
{\footnotesize
\bibliographystyle{plain}
\bibliography{refs3}
}

\end{document}